\documentclass[english]{smfbook}
\usepackage{amscd,amsfonts,amssymb,amscd,amsmath,latexsym,amsthm}
\usepackage{hyperref}
\usepackage[all]{xy}
\makeatletter
\renewcommand{\subsection}{\@startsection
{subsection}
{1}
{0mm}
{0mm}
{0mm}
{\normalfont\normalsize\itshape}}
\makeatother

\usepackage{rotating}
\usepackage{mathrsfs}

\newtheorem{theorem}{Theorem}[section] 
\newtheorem{prop}[theorem]{Proposition}
\newtheorem{lem}[theorem]{Lemma}
\newtheorem{lemma}[theorem]{Lemma}
\newtheorem{ddd}[theorem]{Definition}
\newtheorem{definition}[theorem]{Definition}
\newtheorem{kor}[theorem]{Corollary}
\newtheorem{corollary}[theorem]{Corollary}

\newcommand{\SSSS}{\mathbb{S}} 
\newcommand{\uQ}{\underline{\mathbb{Q}}}
\newcommand{\uF}{\underline{F}}

\newcommand{\Mor}{\mathrm{Mor}}

\newcommand{\uR}{\underline{\R}}

\newcommand{\colim}{\mathrm{colim}}

\newcommand{\Mod}{{\tt Mod}}

\newcommand{\bbA}{\mathbb{A}}

\newcommand{\naturals}{\mathbb{N}}

\newcommand{\tensor}{\otimes}
\newcommand{\into}{\hookrightarrow}

\newcommand{\iso}{\cong}

\newcommand{\eins}{{\mathbf{1}}}

\newcommand{\forget}[1]{}

 {\catcode`@=11\global\let\c@equation=\c@theorem}
 \renewcommand{\theequation}{\thetheorem}

\renewcommand{\theenumi}{(\arabic{enumi})}
\renewcommand{\labelenumi}{\theenumi}

\newcommand{\bU}{\mathbf{U}}
\newcommand{\bLE}{\mathbf{LE}}
\newcommand{\bLB}{\mathbf{LB}}
\newcommand{\uC}{\underline{\mathbb{C}}}

\newcommand{\Z}{\mathbb{Z}}

\newcommand{\bE}{{\bf E}}

\newcommand{\diag}{{\tt diag}}

\newcommand{\R}{\mathbb{R}}

\newcommand{\Q}{\mathbb{Q}}

\newcommand{\Tor}{{\tt Tor}}

\newcommand{\bB}{\mathbf{B}}

\newcommand{\C}{\mathbb{C}}

\newcommand{\Fl}{{\mathcal Fl}}

\newcommand{\Aut}{{\tt Aut}}

\newcommand{\cZ}{\mathcal{Z}}
\newcommand{\cH}{\mathcal{H}}

\newcommand{\cA}{\mathcal{A}}

\newcommand{\Hom}{{\tt Hom}}

\newcommand{\bLX}{\mathbf{LX}}
\newcommand{\End}{{\tt End}}

\newcommand{\im}{{\tt im}}

\newcommand{\cF}{\mathcal{F}}
\newcommand{\cD}{\mathcal{D}}

\newcommand{\gerbes}{{\tt gerbes}}

\newcommand{\coker}{{\tt coker}}
\newcommand{\id}{{\tt id}}

\newcommand{\nat}{\mathbb{N}}
\newcommand{\supp}{{\tt supp}}

\newcommand{\cB}{\mathcal{B}}
\def\imath{{i}}

\def\hB{\hspace*{\fill}$\Box$ \newline\noindent}

\newcommand{\cS}{\mathcal{S}}

\def\hB{\hspace*{\fill}$\Box$ \\[0cm]\noindent}
\newcommand{\cL}{\mathcal{L}}
 
\newcommand{\bH}{\mathbf{H}}

\newcommand{\pr}{{\tt pr}}
\newcommand{\bX}{\mathbf{X}}
\newcommand{\bY}{\mathbf{Y}}
\newcommand{\bZ}{\mathbf{Z}}

\newcommand{\bV}{\mathbf{V}}

\newcommand{\uZ}{\underline{\Z}}
\newcommand{\uRR}{\underline{R}}
\newcommand{\holim}{{\tt holim\:}}
\newcommand{\Top}{{\tt Top}}

\newcommand{\ori}{{\tt or}} 

\newcommand{\cov}{{\tt cov}}
\newcommand{\Ab}{{\tt Ab}}
\renewcommand{\Pr}{{\tt Pr}}
\newcommand{\Sh}{{\tt Sh}}
\newcommand{\Site}{{\tt Site}}
\newcommand{\Mf}{{\tt Mf}}
\newcommand{\bG}{\mathbf{G}}

\author{Ulrich Bunke}
\address{Mathematische Fakult\"at\\ Universit\"at Regensburg\\ 93040
  Regensburg, Germany}
\email{ulrich.bunke@mathematik.uni-regensburg.de}
\urladdr{http://www.mathematik.uni-regensburg.de/Bunke/index.html}
\author{Thomas Schick}
\address{Mathematisches Institut\\ Georg-August-Universit{\"a}t G{\"o}ttingen\\
Bunsenstr. 3-5\\ 37073 G{\"o}ttingen, Germany}   
\email{schick@uni-math.gwdg.de } 
\urladdr{http://www.uni-math.gwdg.de/schick}
\author{Markus Spitzweck}
\address{Mathematische Fakult\"at\\ Universit\"at Regensburg\\ 93040
  Regensburg, GErmany} \email{markus.spitzweck@mathematik.uni-regensburg.de}

\title{Periodic twisted cohomology and T-duality}
\alttitle{Cohomologie p\'eriodique tordue et T-dualit\'e}

\begin{document}

\frontmatter
\begin{abstract}
  Using the differentiable structure, twisted $2$-periodic de Rham cohomoology
  is well known, and showing up as the target of Chern characters for twisted
  K-theory. 
The main motivation of this work is  a topological interpretation
of two-periodic twisted  de Rham cohomology which is generalizable to
arbitrary topological spaces and at the same time to
arbitrary coefficients.

To this end we develop a sheaf theory in the context of locally compact
topological stacks with emphasis on:
\begin{itemize}\item 
  the construction of the sheaf theory operations in unbounded derived
  categories
\item elements of Verdier duality
\item and integration.
\end{itemize}
The main result is the construction of a functorial periodization 
associated to a $U(1)$-gerbe. 

As a application we verify the $T$-duality
isomorphism in periodic twisted cohomology and in periodic twisted orbispace
cohomology.
\end{abstract}

\begin{altabstract}
  La cohomology de de Rham tordue (periodique avec p\'eriode $2$) est une
  construction bien connue, et elle est importante comme codomaine d'un
  charact\`ere  de 
  Chern pour la K-theorie tordue. 

  La motivation principale de notre livre est  une interpretation
  topologique de la cohomology de de Rham tordue, une interpretation avec
  g\'eneralizations  \'a des espaces arbitraire,
  et aux co\'efficients quelconque. 

 \`A ce but, nous developpons une th\'eorie des faisceaux sur des stacks
 topologiques localement compacts. Nous appuyons
 \begin{itemize}
 \item la construction des operations de la th\'eorie des faisceaux dans les
   cat\'egories deriv\'ees non-born\'ees
 \item \'elements de la dualit\'e de Verdier
 \item et integration.
 \end{itemize}
 Le resultat principal est la construction d'une periodization fonctorielle
 associ\'e a une $U(1)$-gerbe.

Une application est la verification d'un isomorphisme de T-dualit\'e pour la
cohomologie periodique tordue et la cohomologie periodique tordue des
orbi-espaces. 
\end{altabstract}

\altkeywords{cohomologie tordue, Verdier dualit\'e, stack topologique,
  th\'eorie des faisceaux, T-dualit\'e, orbi-espaces, categorie deriv\'ee
  non-born\'ee, cohomologie de de Rham tordue}


\maketitle
\tableofcontents

\mainmatter

 \chapter{Introduction}

\section{Periodic twisted cohomology}

\subsection{}

The twisted de Rham cohomology $H_{dR}(M,\omega)$ of a manifold $M$ equipped with a closed   three form $\omega\in \Omega^3(M)$
is the two-periodic cohomology of the complex
\begin{equation}\label{system1002}
\Omega(M,\omega)_{per}\colon \dots\to \Omega^{ev}(M)\stackrel{d_\omega}{\to} \Omega^{odd}(M)\stackrel{d_\omega}{\to}\Omega^{ev}(M)\to \dots\ ,
\end{equation}
where $d_\omega:=d_{dR}+\omega$ is the sum of the de Rham differential and the operation
of taking the wedge product with the form $\omega$. The two-periodic twisted de Rham cohomology is interesting
as the target of the Chern character from twisted $K$-theory \cite{MR2172633}, \cite{MR1977885}, \cite{MR1911247}, or as a cohomology theory which admits a $T$-duality isomorphism  \cite{MR2080959}, \cite{MR2130624}.

\subsection{}

In \cite{bss} we developed a sheaf theory for smooth stacks. Let $f\colon G\to X$ be
a gerbe with band $U(1)$ over a smooth stack $X$, and consider a closed three-form   $\omega\in \Omega_X^3(X)$ which
represents the image of the Dixmier-Douady class of the gerbe $G\to X$ in de
Rham cohomology. The main result of \cite{bss} states
that there exists an isomorphism 
\begin{equation}\label{system1000}Rf_*f^*\uR_\bX\xleftarrow{\sim} \Omega_X[[z]]_\omega\end{equation}
in the bounded below derived category $D^+(\Sh_\Ab\bX)$ of sheaves of abelian groups on $X$. Here
$\uR_\bX$ denotes the constant sheaf with value $\R$ on $X$. Furthermore, 
$\Omega_X[[z]]_\omega$ is the sheaf of formal power series of smooth forms on $X$,
where $\deg(z)=2$, and its differential is given by
$d_\omega:=d_{dR}+\omega\frac{d}{dz}$. The isomorphism is not canonical, but depends on the
choice of a connection on the gerbe $G$ with characteristic  form $\omega$.

\subsection{}

The complex (\ref{system1002}) can be defined for a smooth stack $X$ equipped with a three-form $\omega\in \Omega^3_X(X)$. It is the complex of global sections of a sheaf of two-periodic complexes $\Omega_{X,\omega, per}$ on $X$.
The complex of sheaves $\Omega_X[[z]]_\omega$ is not two-periodic. The relation between $\Omega_X[[z]]_\omega$
and $\Omega_{X,\omega, per}$ has been discussed in \cite[1.3.23]{bss}. Consider the diagram
\begin{equation}\label{uifhuwefewf111}\cD\colon \Omega(X)[[z]]_\omega\stackrel{\frac{d}{dz}}{\leftarrow}\Omega(X)[[z]]_\omega\stackrel{\frac{d}{dz}}{\leftarrow}\Omega(X)[[z]]_\omega\stackrel{\frac{d}{dz}}{\leftarrow}\dots\ .\end{equation}
Then there exists an isomorphism
\begin{equation}\label{uifhuwefewf}
\Omega_{X,\omega, per} \cong \holim \cD\ .
\end{equation}

\subsection{}

As mentioned above, the isomorphism (\ref{system1000}) depends on the choice
of a connection on the gerbe $G$. Moreover, the diagram
$\cD$ depends on these choices via $\omega$.
In order to construct a natural two-periodic cohomology one must find a
natural replacement 
of the operation $\frac{d}{dz}$ which acts on the left-hand side
$Rf_*f^*\uR_\bX$ of (\ref{system1000}). It is the first goal of this paper to
carry this out properly.

\subsection{}

One can do this construction in the framework of smooth stacks developed in \cite{bss}.
But for the present paper we choose the setting of topological stacks.  Only in Subsection \ref{ufefiwufwefe}
we work in smooth stacks and discuss the connection with \cite{bss}.
In Section \ref{system3003} we develop some aspects of the theory of locally compact stacks and the sheaf theory in this context. For the purpose of this introduction we freely use notions and constructions from this theory.
We hope that the ideas are  understandable by analogy with the usual case of sheaf theory on locally compact spaces.

\subsection{}

Let $G\to X$ be a $U(1)$-banded gerbe over a locally compact stack.
The main object of the present paper is a periodization functor
 $$P_G:D^+(\Sh_\Ab\bX)\to D(\Sh_\Ab\bX)$$
which is functorial in $G\to X$, and where $D^+(\Sh_\Ab\bX)$ and
$D(\Sh_\Ab\bX)$ denote the bounded below and unbounded derived categories of sheaves of abelian groups on the site $\bX$ of the stack $X$.
A simple construction of the isomorphism class of $P_G(F)$ is given in Definition \ref{system188}.
The functorial version is much more complicated. Its construction is completed in Definition \ref{system18}.

\subsection{}

Let us sketch the construction of $P_G$. Recall that gerbes with band $U(1)$ over a locally compact stack $Y$ are
classified by $H^3(Y;\Z)$, and automorphisms of a given $U(1)$-gerbe are
classified by $H^2(Y;\Z)$ \cite{heinloth}.
We
consider the diagram
$$\xymatrix{T^2\times G\ar[d]_p\ar[dr]\ar[rr]^u&&T^2\times G\ar[d]^p\ar[dl]\\
G\ar[dr]_f&T^2\times X\ar[d]&G\ar[dl]^f\\&X&}\ ,$$
where the  automorphism $u$ of gerbes over $T^2\times X$  is classified by
$\ori_{T^2}\times 1\in H^2(T^2\times X;\Z)$, and where   $\ori_{T^2}$ denotes the orientation class of the  two-torus. We define a natural transformation
$$D\colon Rf_*f^*\to Rf_* f^*\colon D^+(\Sh_\Ab\bX)\to D^+(\Sh_\Ab\bX)$$
of degree $-2$ as the composition 
$$D\colon Rf_* f^*\stackrel{units}{\to}Rf^*Rp_*Ru_*u^*p^*f^*\stackrel{fpu=fp}{\to}Rf_*Rp_*p^*f^*\stackrel{\int_p}{\to} Rf_*f^*\ ,$$
where $\int_p\colon Rp_*p^*\to \id$ is the integration map of the oriented $T^2$-bundle $T^2\times G\to G$.

For $F\in  D^+(\Sh_\Ab\bX)$ we form the diagram
$$\cS_G(F)\colon Rf_* f^*(F)\stackrel{D}{\leftarrow} Rf_*f^*(F)[2] \stackrel{D}{\leftarrow}Rf_* f^*(F)[4]\stackrel{D}{\leftarrow} \dots$$
in $D(\Sh_\Ab\bX)$.
\begin{ddd}
We define the periodization
$P_G(F)\in D(\Sh_\Ab\bX)$
of $F$ by 
$$P_G(F):=\holim \cS_G(F)\in D(\Sh_\Ab\bX)\ .$$
\end{ddd}

Note that this introduction is meant as a sketch. In particular, one has to be
aware of the fact that the notion of $\holim$ in a triagulated category is
ambiguous and has to be used with great care, as will be explained below and
in the body of the paper. At present, the above definition only fixes the
isomorphism class of $P_G(F)$.
\subsection{}

The same construction can be applied in  the case of smooth stacks $X$.
It is an immediate consequence of Theorem \ref{theo:D_is_derivative} that there exists an isomorphism 
of the diagrams $S_G(\uR_\bX)$ and $\cD$ (see (\ref{uifhuwefewf111})). Equation  (\ref{uifhuwefewf})
implies the following result.
\begin{kor}
If $X$ is a smooth manifold, then there exists an isomorphism
$$P_G(\uR_\bX)\cong \Omega_{X,\omega,per} $$
in $D(\Sh_\Ab\bX)$.
In particular we have an isomorphism of two-periodic cohomology groups
$H^*_{dR}(X,\omega)\cong H^*(X;P_G(\uR_\bX))$.
\end{kor}

The existence of this isomorphism played the role of a design criterion for the construction of the periodization functor $P_G$. 

%
%

\subsection{}

The operation $D\colon Rf_*f^*(F)\to Rf_*f^*(F)$ is a well-defined morphism in the derived category.
In particular, we get a well-defined diagram $\cS_G(F)\in D(\Sh_\Ab\bX)^{\nat^{op}}$,
where we consider the ordered set $\nat$ as a category. This determines the isomorphism class of the object $P_G(F)\in D(\Sh_\Ab\bX)$.
We actually want to define a periodization \emph{functor}
$$P_G\colon D^+(\Sh_\Ab\bX)\to D(\Sh_\Ab\bX)\ ,$$
which also depends functorially on the gerbe $G\to X$. These functorial properties are
required in our applications to $T$-duality, or if one wants to formulate a statement about the naturality of a Chern character from $G$-twisted $K$-theory with values in the periodic twisted  cohomology $H^*(X;P_G(\uR_\bX))$.

In order to define $P_G(F)$ in a functorial way we must refine the diagram
$\cS_G(F)\in D(\Sh_\Ab\bX)^{\nat^{op}}$ to a diagram in
$D((\Sh_\Ab\bX)^{\nat^{op}})$. This is the technical heart of the present
paper. The details of this construction are contained in Section
\ref{system7} and will be completed in Definition \ref{system18}. Along the
way, we have to use the enhancement of the category of sheaves to bounded
below complexes of flasque sheaves. 

\subsection{}

The periodization functor
$P_G$ can be applied to arbitrary objects in $D^+(\Sh_\Ab\bX)$. In Proposition \ref{system30}
we calculate examples which indicate some interesting arithmetic features of this functor.

\section{$T$-duality}\label{sec:subs-topol-t}

\subsection{}

Topological $T$-duality is a concept which models the underlying topology of mirror symmetry in algebraic geometry or $T$-duality in string theory. We refer to \cite{math.GT/0501487} for a more detailed discussion of the literature. In the present paper we introduce the concept of $T$-duality for pairs $(E,G)$ of a $U(1)$-principal bundle $E\to B$  over a topological stack $B$ together with a topological gerbe $G\to E$ with band $U(1)$ using the notion of a $T$-duality diagram.

\subsection{}

Consider a diagram
\begin{equation}\label{system1004}
\xymatrix{&p^*G\ar[dl]^q\ar[dr]\ar[rr]^u&&\hat p^* \hat G\ar[dl]\ar[dr]^{\hat q}&\\G\ar[dr]^f&&E\times_B\hat E\ar[dl]^p\ar[dr]^{\hat p}&&\hat G\ar[dl]^{\hat f}\\&E\ar[dr]^\pi&&\hat E\ar[dl]^{\hat \pi}&\\&&B&&}\ ,
\end{equation}
where $\pi,\hat \pi$ are $U(1)$-principal bundles, and $f,\hat f$ are gerbes with band $U(1)$.
 In \ref{system1003} 
we describe the isomorphism class of the universal $T$-duality diagram over the classifying stack $\cB U(1)$. 
\begin{ddd}[Definition \ref{eruihfrvc}]
The diagram (\ref{system1004}) is a $T$-duality diagram, if it is locally   isomorphic to the universal $T$-duality diagram.
\end{ddd}
The pair $(\hat G,\hat E)$ is then called a $T$-dual of $(E,G)$.

\subsection{}

In Lemma \ref{lem:identify_T_duality} we will check that this generalizes the
concept of T-duality (for $U(1)$-bundles) from
the classical situation of
principal bundles in the category of spaces
\cite{MR2246781,math.GT/0501487} and the
slightly more general situation of such bundles in orbispaces \cite{MR2246781}
to arbitrary $U(1)$-actions. The situation of semi-free actions is discussed
(in a completely different way) in \cite{pande-2006}. It is an
interesting open problem to relate his approach to the approach used here.

\subsection{}

One of the main themes of topological $T$-duality is the $T$-duality transformation in twisted cohomology theories. 
In \cite{MR2246781} we observed that if the $T$-duality transformation is an isomorphism, then
the corresponding twisted cohomology theory must be two-periodic.

This applies e.g. to twisted $K$-theory. In fact, one can argue that twisted $K$-theory is the universal twisted cohomology theory for which the $T$-duality transformation is an isomorphism\footnote{We thank M. Hopkins for pointing out a proof of this fact.}.

%
%
%

\subsection{}\label{system1005}

Our construction of $P_G$ is designed such that the corresponding $T$-duality transformation is an isomorphism. To this end we define the periodic $G$-twisted cohomology of $E$
with coefficients in $\pi^*F$, $F\in D^+(\Sh_\Ab\bB)$, by
$$H^*_{per}(E,G;\pi^*F):=H^*(E;P_G(\pi^*F))\ .$$
In this case the $T$-duality transformation
$$T\colon H^*_{per}(E,G;\pi^*F)\to H^*_{per}(\hat E,\hat G;\hat\pi^*F)$$
is induced by the composition
\begin{eqnarray*}
R\pi_*P_G(\pi^*F)&\stackrel{unit}{\to} &R\pi_*Rp_*p^*P_G(\pi^*F)\\&\cong& R\pi_*Rp_*P_{p^*G}(p^*\pi^*F)\\&\stackrel{u^*}{\cong}& R\pi_*Rp_*P_{\hat p^*\hat G}(p^*\pi^*F)\\&\stackrel{\pi p=\hat \pi \hat p}{\to}&R\hat \pi_*R\hat p_*P_{\hat p^*\hat G}(\hat p^*\hat \pi^*F)\\&\stackrel{\cong}{\to}&R\hat \pi_*R\hat p_*\hat p^* P_{\hat G}(\hat \pi^*F)\\&\stackrel{\int_{\hat p}}{\to}&
R\hat \pi_* P_{\hat G}(\hat \pi^*(F))\ .
\end{eqnarray*}
Note that here we use the functoriality of the periodization in an essential way.
 
\begin{theorem}[Theorem \ref{system66}]
The $T$-duality transformation in twisted periodic cohomology is an isomorphism.
\end{theorem}

\subsection{}

If $G\to X$ is a gerbe over a nice non-singular space $X$, then
$H^*_{per}(X,G;\uR_{\bX})$ is the correct target of a Chern character from
twisted $K$-theory. 
If $X$ is a topological stack with non-trivial automorphisms of points, then
this no longer correct. At the moment we do understand the special case of
orbispaces. In \cite[Sec.~1.3]{math.KT/0609576} we give a detailed motivation
for the introduction of the twisted delocalized cohomology.

Let $G\to X$ be a topological gerbe with band $U(1)$ over an orbispace $X$.
In \cite[Definition 3.4]{math.KT/0609576} we show that it gives rise to a
sheaf $\cL\in \Sh_\Ab \bLX$, where $LX$ is the loop orbispace of  $X$.

The $G$-twisted delocalized periodic cohomology of $X$ (with complex
coefficients) is defined as (see \cite[Definition 3.5]{math.KT/0609576}) 
$$H^*_{deloc,per}(X,G):=H^*(LX;P_{G_L}(\cL))\ ,$$where
$G_L\to LX$ is defined by the pull-back
$$\xymatrix{G_L\ar[d]\ar[r]&G\ar[d]\\LX\ar[r]&X}\ .$$

Let us now consider a $T$-duality diagram (\ref{system1004}) over an orbispace $B$.
Then we define a $T$-duality transformation
$$T\colon H^*_{deloc,per}(E,G)\to H^*_{deloc,per}(\hat E,\hat G)$$
by a modification of the construction \ref{system1005}.

\begin{theorem}[Theorem \ref{system1006}]
The $T$-duality transformation in twisted delocalized periodic cohomology is an isomorphism.
\end{theorem}

So the situation with twisted delocalized periodic cohomology is better  than with orbispace $K$-theory.
At the moment we do not know a proof that the $T$-duality transformation in twisted orbifold $K$-theory
is an isomorphism (see the corresponding comments in \cite{MR2246781}).
Using the fact that the Chern character is an isomorphism, our result implies that
the $T$-duality transformation in  twisted orbifold and orbispace $K$-theory
is an isomorphism after complexification.

\section{Duality for sheaves on locally compact stacks}

\subsection{}

In Section \ref{system3003} of the present paper we develop some features of a  sheaf theory for locally compact stacks.
Our main results are the construction of the basic setup, of the functor $f^!$,  and the integration $\int_f$ for oriented fiber bundles. Section \ref{system3003}  not only  provides the technical  background for the applications of sheaf theory in the previous sections, but also contains some additional material of independent interest (in particular the results connected with  $f^!$).

\subsection{}

A presheaf $F$ of sets on a topological space $X$ associates to each open subset $U\subseteq X$ a set of sections $F(U)$, and to every inclusion $V\to U$ of open subsets a functorial restriction map
$F(U)\to F(V)$, $s\mapsto s_{|V}$. In short, a presheaf  it is contravariant a functor from the category $(X)$ open subsets of $X$ to sets. A presheaf is a sheaf of it has the following two properties:
\begin{enumerate}
\item If $s,t\in F(U)$ are two sections and there exists an open covering $(U_i)$ of $U$ such that $s_{|U_i}=t_{|U_i}$ for all $i$, then $s=t$.
\item If $(U_i)$ is an open covering of $U$ and $(s_i)$ is a collection of sections $s_i\in F(U_i)$
such that $s_{i|U_i\cap U_j}=s_{j|U_i\cap U_j}$ for all pairs $i,j$, then there exists a section $s\in F(U)$ such that
$s_{|U_i}=s_i$ for all $i$.
\end{enumerate}
The notion of a sheaf is thus determined by the Grothendieck topology on $(X)$ given by the collections of open coverings of open subsets. We will call $(X)$ the small site associated to $X$.

If $X$ is a topological stack, then  the open substacks form a two-category which  does not give the appropriate setting for sheaf theory on $X$. For example, if $G$ is a finite group, then the quotient stack $[*/G]$
is quite non-trivial but does not have  proper open substacks.
On the other hand its identity one-morphism has the two-automorphism group $G$, and in a non-trivial theory sheaves should reflect the two-automorphisms.

\subsection{}

For applications to twisted cohomology a setting for sheaf theory on smooth stacks has been introduced in \cite{bss}. In the present paper we develop a similar theory for topological stacks. There are various choices to be made in order to define the site of a stack in topological spaces. The   sheaf theories associated to these choices will have many features in common, but will differ in others. The main goal of the present paper is the construction of the periodization functor $P_G$ associated to a $U(1)$-banded gerbe $G\to X$. One of the main ingredients of the construction is an integration $\int_f$ for oriented fiber   bundles $f$ with a closed topological manifold as fiber.
In order to define the integration map we need a projection formula which
expresses a compatibility of the pull-back and push-forward operations  with
tensor products, see Lemma \ref{system81}. Already  for the projection formula
in ordinary sheaf theory one needs local compactness assumptions. For this
reason we decided to work generally with locally compact stacks and spaces
though much of the theory would go through under more general or different
assumptions.

\subsection{}

A stack in topological spaces is topological if it admits an atlas $A\to X$. From the atlas we can derive a groupoid $A\times_XA\rightrightarrows A$ which represents $X$ in an appropriate sense.
The stack is called locally compact if one can find an atlas $A\to X$ such that the resulting groupoid
is locally compact (i.e. $A$ and $A\times_XA$ are locally compact spaces).

The site $\bX$ associated to a locally compact stack is the category of locally compact spaces $(U\to X)$ over $X$ such that the  morphisms are  morphisms of spaces over $X$ (i.e. pairs of a morphism between the spaces and a two-morphism filling the obvious triangle.) We require that the structure morphism $U\to X$ has local sections.
 The topology on $\bX$
is again given by the collections of coverings by open subsets of the objects $(U\to X)$.
For many constructions and calculations the restriction functors from sheaves on $\bX$ to sheaves on $(U)$ play a distinguished role. They are used to build the connection between operations with sheaves on the stack $X$ and corresponding classical operations in sheaf theory on the spaces $U$.

\subsection{}

For the theory of stacks in topological spaces in general we refer to \cite{heinloth}, \cite{math.KT/0609576}, \cite{math.AG/0503247}.
Some special aspects of locally compact stacks are discussed in Subsection \ref{wuiefwefwefqwd} of the present paper.

In our treatment of sheaf theory on the site $\bX$ 
we  give a description of the closed monoidal structure on the categories of sheaves  and presheaves of abelian groups $\Sh_\Ab\bX$ and $\Pr_\Ab\bX$ on $\bX$.   The interplay between sheaves and presheaves will be important when we study the compatibility of the monoidal structures with the functors  $$f^*:\Sh_\Ab\bY \leftrightarrows \Sh_\Ab\bX: f_*$$
associated to a morphism of locally compact stacks $f:X\to Y$. 
In general these functors do not come from a morphisms of sites but are constructed in an ad-hoc manner. Because of this we must check under which conditions properties expected from the classical theory  carry over to the present case.

The  derived versions of these functors on the bounded below and unbounded derived categories $D^+(\Sh_\Ab\bX)$ and $D(\Sh_\Ab\bX)$ will play an important role in the present paper. In order to deal with the unbounded derived category we use an approach via model categories.

\subsection{}

Besides the development of the basic set up which we will not discuss further in the introduction
let us now explain the two main results which may be of independent interest. 
\begin{theorem}[Theorem \ref{main123}]\label{zueefewf}
If $f:X\to Y$ is a proper representable map between locally compact stacks such that $f_*$ has finite cohomological dimension, then the functor $Rf_*:D^+(\Sh_\Ab\bX)\to D^+(\Sh_\Ab\bY)$ has a right-adjoint,
i.e. we have an adjoint pair
\begin{equation}\label{hiduwqdwqdwqdq}
Rf_*:D^+(\Sh_\Ab\bX)\leftrightarrows D^+(\Sh_\Ab\bY):f^!\ .
\end{equation}

\end{theorem}
We  think that one could prove a more general theorem stating  the existence of a right adjoint
of a functor $Rf_!$ where $f_!$ is the push-forward with proper support along an arbitrary map between
locally compact stacks such that $f_!$ has finite cohomological dimension, though we have not checked all details.  

This theorem generalizes a well-known result  (\cite{MR1610971}, \cite[Ch.~3]{MR1299726} in ordinary sheaf theory. Its importance is due to the classical calculation 
\begin{equation}\label{wuiefwefwef}
f^!(F)\cong f^*(F)[n]
\end{equation} (compare \cite[Prop.3.3.2]{MR1299726}) for $F\in D^+(\Sh_\Ab(Y))$, if $f:X\to Y$ is an oriented  locally trivial bundle of closed connected topological $n$-dimensional manifolds on a locally compact space $Y$. If we would know such an  isomorphism in the present case (for sheaves on the sites $\bX, \bY$ and stacks $X,Y$), then  we could define the integration map as the composition
$$\int_f:Rf_*f^*(F)\stackrel{\sim}{\to} Rf_*f^!(F)[-n]\stackrel{counit}{\to} F[-n]\ ,$$
where the last map is the co-unit of the adjunction (\ref{hiduwqdwqdwqdq}).

Unfortunately, at the moment we are not able to calculate $f^!(F)$ in any interesting example.
However, we can construct the integration map in a direct manner avoiding the knowledge of (\ref{wuiefwefwef}).

Some elements of the theory developed here are formally similar to the work \cite{MR2312554}
on sheaves on the lisse \'etale site of an Artin stack. In this framework in \cite{laszlo-2005} a functor
$f^!$ was introduced  between derived categories of constructible sheaves.  On the one hand the methods seem to be completely different. On the other hand this functor has the expected behavior for
smooth maps, i.e. it satisfies a relation like (\ref{wuiefwefwef}).  At the moment we do not see  even a formal relation between  the construction of \cite{laszlo-2005} with the construction in the present paper which could be exploited for a calculation of $f^!(F)$. 

\subsection{}

The following Theorem is the result of Subsection \ref{system103}.
\begin{theorem}
If the map $f:X\to Y$ of locally compact stacks  is an oriented locally trivial fiber bundle with a closed connected topological $n$-dimensional manifold as fiber,
then there exists an integration map, a natural transformation of functors
$$\int_f:Rf_*f^*\to \id[-n]:D^+(\Sh_\Ab \bX)\to D^+(\Sh_\Ab\bX)$$
which has the expected compatibility with pull-back and compositions.
\end{theorem}
In Subsection \ref{system4001} we extend the push-forward and pull-back operations to the unbounded derived categories and construct the integration map in this setting.

\chapter{Gerbes and periodization}\label{system3001}

\section{Sheaves on the locally compact site of a stack}

\subsection{}\label{system1007}

Let $\Top$ denote the site of topological spaces. The topology is generated by  covering families
$\cov_\Top(A)$ of the objects $A\in \Top$, where $\cov_\Top(A)$ is the set  of
coverings by collections of open subsets.

A stack will be a stack on the site $\Top$. Spaces are considered as stacks
through the Yoneda embedding.

A map  $A\to X$ from a space $A$ to a stack $X$ which is surjective,
representable, and has local sections is called an atlas. 
We refer to \ref{staqwndwqodwqd} for definitions and more details about stacks in topological spaces.
\begin{ddd}
A topological stack
is a stack which admits an atlas. 
\end{ddd}

\begin{ddd}\label{qwuidiuwqdwqdwqd1fwefw}
A topological space is locally compact if it is Hausdorff and every point admits a compact neighborhood.
 A stack is called locally compact if it admits an atlas $A\to X$ such that $A$ 
and  $A\times_XA$ are locally compact.
\end{ddd}

If $X$ is a locally compact stack, then the site of $X$ is the subcategory $\Top_{lc}/X$ of locally compact spaces over $X$ such that the structure map $A\to X$ has local sections. The  topology is induced from $\Top$. We denote this site by $\bX$ or $\Site(X)$. See \ref{hshsqiwhswqs} for more details.

\subsection{}

As will be explained in \ref{sec:morph_of_sheaves}, a morphism of locally compact
stacks $f\colon X\to Y$ gives rise to an adjoint pair of functors
$$f^*\colon \Sh\bY\leftrightarrows\Sh\bX:f_*\ .$$
The functor
$f_*$ is left-exact on the categories of sheaves of abelian groups and admits
a right-derived
$$Rf_*\colon  D^+(\Sh_\Ab\bX)\to D^+(\Sh_\Ab\bY)$$
between the bounded below derived categories, compare  \ref{lem:lrexact}.

%
%

\subsection{}\label{system102}


Let $M$ be some space.
\begin{ddd}
A map between topological stacks $f:X\to Y$ is a locally trivial fiber bundle with fiber $M$ if
for every space $U\to X$ the pull-back $U\times_YX\to U$ is a locally trivial
fiber bundle of spaces with fiber $M$.
 \end{ddd}

Assume that $M$ is a closed connected and  orientable $n$-dimensional topological manifold.
\begin{ddd}\label{iquhuiqwdddqwfefefefefeffefeffefe}
Let $f\colon X\to Y$ be a map of locally compact stacks which
is a locally trivial fiber bundle  with fiber $M$. It is called orientable if
there exists an isomorphism $R^nf_*(\uZ_\bX)\cong \uZ_\bY$. An orientation of $f$ is a choice of such an isomorphism.
\end{ddd}


\subsection{}

Let $f\colon X\to Y$ be a locally trivial oriented fiber bundle with $n$-dimensional fiber $M$
over a locally compact stack $Y$. Under these assumption we can generalize the integration map
(see \cite[Sec.~3.3]{MR1299726}) 
\begin{theorem}[Definition \ref{rdphi}]\label{system2}
If $f\colon X\to Y$ be a locally trivial oriented fiber bundle over a locally compact stack with fiber a closed topological manifold of dimension $n$, then we have an integration map, i.e. a natural transformation of functors
$$\int_f\colon Rf_*\circ f^*\to \id\colon D^+(\Sh_\Ab\bY)\to D^+(\Sh_\Ab\bY)$$
of degree $-n$.
\end{theorem}

\subsection{}


We consider a map of locally compact stacks $f\colon X\to Y$ which is a locally trivial oriented fiber bundle  with fiber a closed topological manifold of dimension $n$. 
Furthermore let $U\to X$ be a morphisms of locally compact stacks which has local sections. Then we form the Cartesian\footnote{{In the present paper by a Cartesian diagram in the two-category of stacks we mean a $2$-Cartesian diagram. In particular, the square commutes up to a $2$-isomorphism which we often omit to write in order to simplify the notation. More generally, when we talk about a commutative diagram in stacks, then we mean a diagram of $1$-morphisms together with a collection of $2$-isomorphism filling all faces in a compatible way, and again we will usually not write the $2$-isomorphisms explicitly.}} diagram
$$\xymatrix{V\ar[r]^v\ar[d]^g&X\ar[d]^f\\U\ar[r]^u&Y}\ .$$
Note that $g:V\to U$ is again a locally trivial oriented fiber bundle  with fiber a closed topological manifold of dimension $n$. The orientation of $f$ (which gives the marked isomorphism below) induces an orientation of $g$ by
$$R^ng_*(\uZ_\bV)\cong R^ng_*v^*(\uZ_\bX)\stackrel{(\ref{wqdwqdqw})}{\cong}u^*R^nf_*(\uZ_\bX)\stackrel{!}{\cong} u^*(\uZ_\bY)\cong \uZ_{\bU}\ .$$

 
\begin{lem}\label{ghtq}
The following diagrams commute
\begin{equation}\label{fde4}\xymatrix{u^*\circ Rf_*\circ f^*\ar[r]^\cong\ar[d]^{u^*\int_f}&Rg_*\circ v^*\circ f^*\ar[d]^\cong\\
u^*&Rg_*\circ g^*\circ u^* \ar[l]^{\int_g}}
\quad 
 \xymatrix{Ru_*\circ Rg_*\circ g^*\ar[r]^\cong\ar[d]^{Ru_*\int_g}&Rf_*\circ Rv_*\circ g^*\ar[d]^\cong\\
   Ru_*&Rf_*\circ f^*\circ Ru_* \ar[l]^{\int_f Ru_*}}
\ .\end{equation}
\end{lem}
\begin{proof}
Commutativity of the first
diagram follows immediately from the stronger (because valid in the derived
category of unbounded complexes)  Lemma
\ref{system300}. Commutativity of the second diagram is proved in Lemma
\ref{lem:funct_int1}, but only for the bounded below derived category. 
\end{proof}

\section{Algebraic structures on the cohomology of a gerbe}

\subsection{}\label{ddef2}

Let $X$ be a locally compact  stack and
$f\colon G\to X$ be a topological gerbe with band $U(1)$.  
Then $G$ is a locally compact stack.
Indeed, we can choose an atlas $A\to X$ such that $A$ and $A\times_XA$ are locally compact, and there exists a section 
$$\xymatrix{&G\ar[d]\\A\ar@{.>}[ur]\ar[r]&X}\ .$$
Then $A\to G$ is an atlas and $A\times_GA\to A\times_XA$ is a locally trivial $U(1)$-bundle. In  particular, $A\times_GA$ is a locally compact space.

\subsection{}\label{uidqwdqwdqdqwdd}

By $T^2$ we denote the two-dimensional torus. We fix an orientation of $T^2$.
We consider the pull-back $\pr_2^*G\cong T^2\times G\to T^2\times X$. The isomorphism classes of automorphisms of this gerbe are classified by $H^2(T^2\times X;\Z)$. Let
$$\xymatrix{\pr^*_2G\ar[rr]^\phi\ar[dr]& &\pr_2^*G\ar[dl]\\&T^2\times X&}$$ 
be an automorphism classified by
$\ori_{T^2}\times 1_X\in H^2(T^2\times X;\Z)$.
We consider the diagram
\begin{equation} 
  \xymatrix{\pr_2^* G\ar[rd]\ar[rr]^\phi\ar[d]^{p}&&\pr_2^*G\ar[ld]\ar[d]^{p}\\
G\ar[dr]^f&T^2\times X\ar[d]&G\ar[dl]^f \\&X&}\ .
\label{eq:dtdul}
\end{equation}

{
Notice that $\phi$ is unique up to a non-canonical $2$-isomorphism. In the present paper we
prefer a more canonical choice. We will fix the morphism $\phi$ once and for all
in the special case that $X$ is a point and $G=\cB U(1)$, i.e. we fix a diagram
$$\xymatrix{T^2\times \cB U(1)\ar[d]\ar[rr]^{\phi_{univ}}\ar[dr]&&T^2\times \cB U(1)\ar[d]\ar[dl]\\\cB U(1)\ar[dr]&T^2\ar[d]&\cB U(1)\ar[dl]\\&{*}&}\ .$$
If $G\to X$ is a topological gerbe with band $U(1)$, then we obtain the induced diagram by taking products
$$\xymatrix{G\times T^2\times  \cB U(1)\ar[d]\ar[rr]^{\id_G\times \phi_{univ}}\ar[dr]&&G\times T^2\times \cB U(1)\ar[d]\ar[dl]\\G\times \cB U(1)\ar[dr]&X\times T^2\ar[d]&G\times \cB U(1)\ar[dl]\\&X&}\ .$$
We now replace the products $\cB U(1)\times G$ by the tensor product of gerbes
as explained in \cite[6.1.9]{math.AT/0701428} and identify $\cB U(1)\otimes G$ with $G$ using the canonical isomorphism in order to get
$$\xymatrix{\pr_2^*G\ar[d]^p\ar[rr]^{\phi}\ar[dr]&&\pr_2^* G \ar[d]^p\ar[dl]\\G \ar[dr]^f& T^2\times X\ar[d]&G\ar[dl]_f\\&X&}\ .$$
In this way we have constructed a $2$-functor from the $2$-category of $U(1)$-banded gerbes over $X$ to the $2$-category of diagrams of the form (\ref{eq:dtdul}).
By taking prefered models for the products we can, if we want, assume a strict equality
$f\circ p\circ \phi_G=f\circ p$.
}

\subsection{}

Observe that the map of locally compact stacks $p\colon \pr_2^*G\to G$ is a locally trivial oriented fiber bundle with fiber $T^2$.
Therefore we have the integration map (see \ref{system2})
$$\int_p\colon Rp_*\circ p^*\to\id\ .$$
\begin{ddd}\label{defofd}
We define a natural endo-transformation $D_G$ of the functor 
$$Rf_*\circ f^*\colon D^+(\Sh_\Ab\bX)\to
D^+(\Sh_\Ab\bX)$$ of degree $-2$ which associates to $F\in D^+(\Sh_{\Ab}\bX)$ the morphism  

\begin{multline*}
  Rf_*\circ f^*(F)\stackrel{units}{\longrightarrow}Rf_*\circ Rp_*\circ
  R\phi_*\circ \phi^* \circ p^*\circ f^*(F)\\
\xrightarrow{f\circ p\circ
    \phi= f\circ p} Rf_*\circ Rp_*\circ p^*\circ f^*(F) \stackrel{\int_p}{\to}
  Rf_*\circ f^*(F)\ .
\end{multline*}
 \end{ddd}

\subsection{}
It follows from Lemma \ref{ghtq} that $D_G$ is compatible with pull-back diagrams.
In fact, consider a Cartesian diagram 
$$\xymatrix{ G^\prime\ar[d]^{f^\prime}\ar[r]&G\ar[d]^f\\
X^\prime\ar[r]^g&X}\ .$$
{Using the canonical construction explained in \ref{uidqwdqwdqdqwdd} we  extend this to a morphism between diagrams of the form (\ref{eq:dtdul}).}
 Then we have  the  commutative diagram
$$\xymatrix{g^*\circ Rf_*\circ f^*\ar[r]^\sim\ar[d]^{g^*D_G}&Rf^\prime_*\circ (f^\prime)^*\circ g^*\ar[d]^{D_{G^\prime}\circ g^*}\\
g^*\circ Rf_*\circ f^*\ar[r]^\sim&Rf^\prime_*\circ (f^\prime)^*\circ g^*}\ .$$

\subsection{}

We compute the action of $D_G$ in the case of the trivial gerbe $f:G\to *$ and the sheaf $\uF\in \Sh_\Ab\Site(*)$ represented by a discrete abelian group $F$. Note that
$Rf_*\circ f^*(\uF)$ is an object of $D^+(\Sh_\Ab \Site(*))$. We get an object 
$Rf_*\circ f^*(\uF)(*)\in D^+(\Ab)$ by evaluation at the object $(*\to *)\in \Site(*)$.
\begin{lem}\label{system5}
There exists an isomorphism
$$H^*(Rf_*\circ f^*(\uF)(*))\cong F\otimes \Z[[z]]\ ,$$
where $\deg(z)=2$. On cohomology the transformation $D_G$ is given by  $D_G=\id\otimes \frac{d}{dz}$.
\end{lem}
\begin{proof}
 We choose a lift
$*\to G$. Forming iterated fiber products we get a simplicial space
$$\dots *\times_G*\times_G*\times_G*\to *\times_G*\times_G*\to *\times_G*\to *\ .$$
Note that $*\times_G*\cong U(1)$.
One checks that the simplicial space is equivalent to the simplicial space
$BU(1)^\cdot$, the classifying space of the group $U(1)$,
$$U(1)\times U(1)\times U(1)\to U(1)\times U(1)\to U(1)\to *\ .$$
Let $(U\to *)\in \Site(*)$. 
If $H\in \Sh_\Ab\bG$, then  we consider an injective resolution $0\to H\to I^\cdot$.
The evaluation $I^\cdot(U\times BU(1)^\cdot)$ gives a cosimplicial complex, and after normalization, a double complex. Its total complex represents $Rf_*(H)(U\to *)$ (see \cite[Lemma 2.41]{bss} for a proof of the corresponding statement in the smooth context).
We calculate the cohomology of  $Rf_*(H)(U\to *)$ using the associated spectral sequence.
Its second page has the form $$E_2^{p,q}\cong H^p(U\times BU(1)^q;H)\ .$$

We now specialize to the sheaf $H=f^*(\uF)\cong \uF_G$, where $F$ is a discrete abelian group, and $U=*$. In this case the spectral sequence is the usual spectral sequence which calculates the cohomology of the realization of the simplicial space $BU(1)^\cdot$ with coefficients in $F$.
Note that $H^*(B U(1);\Z)\cong  \Z[[z]]$ as rings with $\deg(z)=2$. Since it is torsion free as an abelian group
we get
$$H^*(R^*f_*\circ f^*(\uF)(*))\cong  F\otimes  H^*(B U(1);\uZ)\cong F\otimes  \Z[[z]]\ .$$
In a similar manner we calculate $Rf_*\circ Rp_*\circ p^*\circ f^*(\uF)(*)$. Its cohomology is
$H^*(T^2\times B U(1);F)$, hence we have
$$H^*(Rf_*\circ Rp_*\circ p^*\circ f^*(\uF) (*))\cong F\otimes H^*(T^2\times B U(1);\Z)\cong F\otimes \Lambda(u,v)\otimes  \Z[[z]]
\ ,$$
where $u,v\in H^1(T^2,\Z)$ are the canonical generators.


For every topological group $\Gamma$ we have a natural map
$\Gamma \to \Omega(B\Gamma)$. By adjointness we get a map
$c:U(1)\times \Gamma\to U(1)\wedge \Gamma\to B\Gamma$.
We will need a simplicial model $c^\cdot$ of this map.
We consider the standard simplicial model $\SSSS^\cdot$ of $U(1)$ with two non-degenerate simplices, one in degree $0$, and one in degree $1$. Then $\SSSS^\cdot\times \Gamma$ is a simplicial model
of $U(1)\times \Gamma$. It suffices to describe the map $c^\cdot$ on the non-degenerate part of $\SSSS^\cdot\times \Gamma$. The component $c^0$ maps
$\SSSS^0\times \Gamma$ to the base point $*$ of $B\Gamma^\cdot$.
 The component $c^1$ is the natural identification of the non-degenerate copy of $\Gamma\subset
\SSSS^1\times \Gamma$ with $\Gamma\cong B\Gamma^1$.

We now specialize to the case $\Gamma=U(1)$.
We get a map
$c:T^2\cong U(1)\times U(1)\to B U(1)$, or on the simplicial level, a map
$c^\cdot:\SSSS^\cdot\times U(1)\to B U(1)^\cdot$.
We have $H^*(BU(1);\Z)\cong \Z[[z]]$ with $z$ odd degree $2$, and one checks that $uv=c^*(z)\in H^2(T^2;\Z)$
(after choosing an appropriate basis $u,v\in H^1(T^2;\Z)$).

Note that $BU(1)^\cdot$ is a simplicial abelian group. 
The discussion above shows that the automorphism $\phi\colon G\to G$ in (\ref{eq:dtdul}) with $X=*$ and classified by $uv\in H^2(T^2;\Z)$ can be arranged so that it induces
an automorphism of bundles of $BU(1)^\cdot$-torsors
\begin{equation}\label{uqfbqwdqwqwd}
\xymatrix{\SSSS^\cdot\times U(1)\times BU(1)^\cdot\ar[dr]\ar[rr]_{\phi^{\cdot}}^{(t,x)\mapsto (t,c^\cdot(t)x)}&&\SSSS^\cdot\times U(1)\times BU(1)^\cdot\ar[dl]\\&\SSSS^\cdot\times U(1)&}\ .
\end{equation}

Under this isomorphism 
the action of
\begin{equation}
\phi^*\colon H^*(Rf_*\circ Rp_*\circ p^*\circ f^*(\uF)(*)) \to H^*(Rf_*\circ Rp_*\circ
p^*\circ f^*(\uF)(*))
\label{eq:action_of_phi}
\end{equation}
 is induced by  $z\mapsto z+u v$, $u\mapsto u$, $v\mapsto v$. 
In order to see this note that
$m^*(z)=z_1+z_2$, where  $m:BU(1)\times BU(1)\to BU(1)$ is the multiplication,
and $H^*(BU(1)\times BU(1);\Z)\cong \Z[[z_1,z_2]]$.
After realization the  map $\phi^\cdot$ leads to  the composition
$$T^2\times BU(1)\stackrel{(\id_{T^2}, c)\times \id}{\to}T^2\times  BU(1)\times BU(1)\stackrel{\id_{T^2}\times m}{\to} T^2\times BU(1)$$
which maps $$z\stackrel{(\id_{T^2}\times m)^*}{\mapsto} z_1+z_2\stackrel{((\id_{T^2},c)\times \id)^*}{\mapsto} uv+z\ .$$
In cohomology of the evaluations at the point the integration
map
$$\int_p\colon 
Rf_*\circ Rp_*\circ p^*\circ f^*(\uF)\to Rf_*\circ f^*(\uF)$$
induces the map 
$F\otimes \Lambda(u,v)\otimes  \Z[[z]] \to F\otimes  \Z[[z]]$ which takes
the coefficient at $uv$. This implies the assertions of Lemma \ref{system5}. 
\end{proof}

\section{Identification of the transformation $D_G$ in the smooth case}\label{ufefiwufwefe}

\subsection{}
In this subsection we work in the context of \cite{bss} of manifolds and smooth stacks.
It can be considered as a supplement to \cite{bss} concerning the  transformation $D_G$ introduced in Definition \ref{defofd} which can be defined in the smooth context in  a parallel manner.

If $X$ is a smooth stack, then $\Omega_X$ denotes the sheaf of de Rham complexes on $X$.
It associates to $(U\to X)\in \bX$ the de Rham complex $\Omega_X(U\to X):=\Omega(U)$ of the manifold $U$. Note that in this subsection $\bX$ denotes the site of a smooth stack introduced in \cite{bss}.

If $\omega\in \Omega_X^3(X)$ is a closed $3$-form, then we form the sheaf of twisted de Rham complexes $ \Omega_X[[z]]_\omega$. Its evaluation at
$(U\to X)\in \bX$ is the complex $\Omega_X[[z]]_\omega(U\to X):=\Omega(U)[[z]]\cong \Omega(U)\otimes_\Z\Z[[z]]$ with differential
$d_{dR}+\omega \frac{d}{dz}$.
In this formula the form $\omega$ acts by wedge multiplication with the pull-back of $\omega$ to $U$.
 
Let  $f\colon G\to X$ be a gerbe with band $U(1)$ over a \textit{smooth manifold} $X$. 
The choice of a gerbe connection determines a closed $3$-form $\omega\in \Omega_X^3(X)$ which  represents the
Dixmier-Douady class of the gerbe. By
\cite[Theorem 1.1]{bss} we have an
isomorphism
\begin{equation}\label{wukhwefwefw}
  Rf_* f^* \uR_\bX\stackrel{\sim}{\to} \Omega_X[[z]]_\omega 
\end{equation}
in the derived category $D^+(\Sh_\Ab\bX)$.


\subsection{}

\begin{theorem}\label{theo:D_is_derivative}
We have a commutative diagram
  \begin{equation*}
    \begin{CD}
      Rf_* f^*\uR_\bX @<{\iso}<(\ref{wukhwefwefw})< \Omega_X[[z]]_\omega\\
      @VV{D_G}V   @VV{\frac{d}{dz}}V\\
      Rf_*f^*\uR_\bX @<{\iso}<(\ref{wukhwefwefw})< \Omega_X[[z]]_\omega.
    \end{CD}
  \end{equation*}
\end{theorem}
\proof
The isomorphism (\ref{wukhwefwefw})  was constructed in \cite[Section 3]{bss} using a particular model of
$Rf_*f^*(\uR_\bX)$.  We first recall its construction.
Let $A\to G$ be  an atlas.
For $(U\to X)\in \bX$ 
we form the simplicial object
$(A_U^\cdot\to G)\in \bG^{\Delta^{op}}$ with $n$th piece $$
A_U^n:=\underbrace{A\times_G \dots \times_G A}_{n+1\text{ factors}}\times_X U\to G\ .$$
The boundaries and degenerations are given by the projections and diagonals as usual.

If $F\in C^+(\Pr_\Ab\bG)$ is a bounded below complex of  presheaves, then
we form the simplicial complex of presheaves
$(U\to X)\mapsto F(A_U^\cdot\to G)$.
We let $C_A(F)\in C^+(\Pr_\Ab\bX)$ denote the presheaf of associated total complexes.
Sometimes we will write $C_A^{m,n}(F) $ for the summand of bidegree $(m,n)$, where the first entry
$m$ denotes the cosimplicial degree.

If $F$ is a complex of flabby sheaves, then by \cite[Lemma 2.41]{bss}
we have a natural isomorphism $Rf_*(F)\cong C_A(F)$. Here we use in particular that the functor $C_A$ preserves sheaves.

Note that the resolution $\uR_G\to \Omega_G$ of the constant sheaf with value $\R$ by the sheaf of de Rham complexes is a flabby resolution (see \cite[Subsection 3.1]{bss}). Therefore we have a natural isomorphism $Rf_*(\uR_G)\cong C_A(\Omega_G)$.


We choose an atlas $A\to X$ given by the disjoint union of a collection of open subsets of $X$  such that there exists a lift in
$$\xymatrix{&G\ar[d]^f\\A\ar[r]\ar@{.>}[ur]&X}\ .$$
This lift is an atlas $A\to G$ of $G$.
 We furthermore
 choose a connection datum $(\alpha,\beta)\in \Omega^1(A\times_G A)\times
\Omega^2(A)$. The one-form  $\alpha$ is a connection of
the $U(1)$-principal bundle $A\times_G A\to A\times_X A$. It is related with the two-form $\beta$ by 
$d_{dR}\alpha=\delta \beta$. This equation implies that $\delta d_{dR}\beta=0$ so that 
$d_{dR}\beta$
assembles to a uniquely determined closed form $\omega\in \Omega^3_X(X)$ (compare \cite[Section 3.2]{bss}).
The  $3$-form $\omega$ represents the Dixmier-Douady class of the gerbe $G\to X$ and will be used for twisting the de Rham complex.

 The isomorphism (\ref{wukhwefwefw}) is given by an explicit quasi-isomorphism \begin{equation}\label{dgqwudwqdwq}
\Omega_X[[z]]_\omega\to C_A(\Omega_G)\ .
\end{equation}\
Note that $\Omega_X[[z]]_\omega$ and $C_A(\Omega_G)$ are sheaves of
associative $DG$-algebras central over the sheaf of $DG$-algebras $\Omega_X$, and that $z$ generates
 $\Omega_X[[z]]_\omega$.
The quasi-isomorphism (\ref{dgqwudwqdwq}) is the unique morphism of sheaves of associative $DG$-algebras, central over $\Omega_X$, with
\begin{equation*}
z\mapsto (\alpha,\beta)\in
C_A^{1,1}(\Omega_G)(X)\oplus C_A^{0,2}(\Omega_G)(X) \ .
\end{equation*}
For more details we refer to \cite[Subsection 3.2]{bss}

\subsection{}

For $i=1,\dots,n$ there are $U(1)$-principal bundle  structures  
$$p_i:\underbrace{A\times_G\cdots \times_G A}_{n+1\text{ factors}}\to \underbrace{A\times_G\cdots \times_G A}_{i \text{ factors}}\times_X\underbrace{
  A\times_G\cdots \times_G A}_{n-i+1 \text{ factors}}\ .$$
Furthermore, we have embeddings
$$j_i:\underbrace{A\times_G\cdots \times_G A}_{n\text{ factors}}
\to \underbrace{A\times_G\cdots \times_G A}_{i \text{ factors}}\times_X\underbrace{
  A\times_G\cdots \times_G A}_{n-i+1 \text{ factors}}$$
given by 
$$j_i:= \underbrace{\id_A\cdots \times \id_A}_{i-1 \text{ factors}}\times \Delta_A\times \underbrace{
  \id_A\cdots \times \id_A}_{n-i \text{ factors}}\ ,$$
where $\Delta:A\to A\times_XA$ is the diagonal.

If $(U\to X)\in \bX$, then 
the maps $p_i$ and $j_i$ induce similar maps on the product $\dots\times_XU$ of these manifolds over $X$ with $U$
which we denote by the same symbols.
For $i=1,\dots,n$ we define the map of degree $-1$
$$v_i\colon \Omega(A_U^n)\to \Omega(A_U^{n-1})$$
as the composition of the integration over the fiber of $p_i$ with the pull-back along $j_i$, i.e.
$v_i:=j_i^*\circ \int_{p_i}$.
Since the construction of $v_i$ is natural with respect to $U$
we can view $v_i$ as a morphism of sheaves $C_A^{n,m}(\Omega_G)\to C_A^{n-1,m-1}(\Omega_G)$.  
We define the family of morphisms 
$$D_n:=\sum_{i=1}^n (-1)^i v_i:
C_A^{n,*}(\Omega_G)\to C_A^{n-1,*-1}(\Omega_G)
$$
and let $D:C_A(\Omega_G)\to C_A(\Omega_G)$ be the endomorphism of sheaves of degree $-2$
given by $D_n$ in bidegree $(n,*)$.

%
%

\subsection{}

\begin{lem}
The map $D:C_A(\Omega_G)\to C_A(\Omega_G)$ is a derivation
of  $\Omega_X$-modules.
\end{lem}
\begin{proof}
Note that $v_j$ commutes with the de Rham differential. Moreover, if
$$q_k\colon \underbrace{A\times_G\dots\times_GA}_{n+1 \text{ factors}}\to \underbrace{A\times_G\dots\times_GA}_{n \text{ factors}}$$ is the projection which leaves out the $k$-th
factor $(k=0,\dots,n$), then we have the relations
\begin{equation*}
  \begin{split}
   v_jq_k^* &= q_{k-1}^* v_j, \qquad  j<k\\
    v_jq_k^* &= q_k^* v_{j-1}, \qquad j>k+1\\
    v_jq_k^* &= 0, \qquad j=k,k+1.
  \end{split}
\end{equation*}
Observe that in the last case $q_k$ factors over the bundle 
 which is used for the integration in the definition
of $v_k$ or $v_{k+1}$, and the composition of a pullback along a bundle projection followed by an 
integration along the same bundle projection vanishes.
These relations imply by a direct calculation that $D$ is a chain map 
  for the \v{C}ech-de Rham differential of
$C_A(\Omega_G)$. 

Moreover, it follows immediately from the definition of $D$  that it is
$\Omega_X$-linear (even $\Omega_A$-linear).

It is again a straightforward calculation to verify that $D$ is a derivation for that
associative product on $C_A(\Omega_G)$ (compare \cite[2.4.9]{bss} for the
product structure).
\end{proof}

\subsection{}

\begin{lemma}\label{dz_D}
  We have a commutative diagram
\begin{equation*}
  \begin{CD}
    \Omega_X[[z]]_\omega @>{(\ref{dgqwudwqdwq})}>> C_A(\Omega_G) \\
    @VV{\frac{d}{dz}}V @VV{D}V\\
    \Omega_X[[z]]_\omega @>{(\ref{dgqwudwqdwq})}>> C_A(\Omega_G) .
  \end{CD}
\end{equation*}
\end{lemma}
\begin{proof}
Since $\alpha$ is the connection one-form of  a $U(1)$-connection on the total space   of the $U(1)$-principal bundle $p_1:A\times_G A\to A\times_X A$  we have
$\int_{p_1}\alpha=1$. Consequently,
$D(\alpha,\beta)=1$.
This implies the assertion,  since $D$  and $\frac{d}{dz}$  are $\Omega_X$-linear
derivation, and $z$ generates $ \Omega_X[[z]]_\omega$.
\end{proof}
 In view of Lemma \ref{dz_D}, in order to finish the proof of Theorem \ref{theo:D_is_derivative} is suffices to show that the operation $D$ coincides with
the operation of 
$\int_p\circ \phi^*\circ p^*$ on $C_A(\Omega_G)$.

\subsection{}

Let $M^\cdot$ be a simplicial manifold and consider the bundle 
$U(1)\times M^\cdot\to M^\cdot$. 
We describe the integration map
$$\int:\Omega(U(1)\times M^\cdot)\to  \Omega(M^\cdot)$$ in the
simplicial picture, i.e. as a map
$$\int:\Omega(\SSSS^\cdot\times M^\cdot)\to  \Omega(M^\cdot)\ .$$
For $n\ge 1$ the manifolds
$\SSSS^{n}\times M^n$ consists of
$n$ copies $\sigma_1(M^n),\dots,\sigma_n(M^n)$ of $M^n$ which correspond to the points of $\SSSS^n$ which are degenerations of the non-degenerated point of $\SSSS^1$  (where the index measures which $1$-simplex in the boundary
  is non-degenerate), and an additional copy of $M^n$ corresponding the point of $\SSSS^n$ which is the  degeneration of the point in $\SSSS^0$.  For $k=1,\dots,n+1$ let  $j_k: M^n\to \SSSS^{n+1}\times M^{n+1}$ be the map $M^n\to \sigma_k(M^{n+1})\subset \SSSS^{n+1}\times M^{n+1}$, which corresponds the $k$th degeneration
$[n+1]\to [n]$.   We now define a chain map of total complexes
$$\int:\Omega(\SSSS^\cdot\times M^\cdot)\to \Omega(M^\cdot)$$
of degree $-1$ which is given by 
\begin{equation}\label{uefqwefqewfqfwqw}
\sum_{k=1}^{n+1} (-1)^kj_k^*:\Omega(\SSSS^{n+1}\times M^{n+1})\to \Omega(M^n)\ ,
\end{equation}
and is zero on $\Omega(\SSSS^0\times M^0)$.
This map realizes the integration in the simplicial picture.
\subsection{}

For $(U\to X)\in \bX$
the automorphism of gerbes $\phi:T^2\times G\to T^2\times G$ induces an automorphism
of simplicial sets
$$\phi^\cdot:\SSSS^\cdot\times U(1)\times A_U^\cdot\to \SSSS^\cdot\times U(1)\times A_U^\cdot$$
which we now describe explicitly by an extension of the special case (\ref{uqfbqwdqwqwd}) 
to general base spaces.

If  $t\in \SSSS^n\times U(1)$ belongs to $U(1)\cong\sigma_k(U(1))\subset \SSSS^n\times U(1)$, $k=1,\dots,n$, then $\phi^\cdot(t,a)=(t, m_k(t,a))$, where $m_k:U(1)\times A_U^n\to A_U^n$ is the action of $U(1)$ on the principal fibration $p_k$. 
If $t\in  \SSSS^n\times U(1)$ belongs to the degeneration of $\SSSS^0\times U(1)$, then $\phi^\cdot(t,a)=(t, a)$.
This formula provides a simplicial description of the action of
$$\phi^*:C_{\SSSS^\cdot\times U(1)\times A}(\Omega_{G})\to C_{\SSSS^\cdot\times U(1)\times A}(\Omega_{G})\ .$$

Combining the description of the integration map (\ref{uefqwefqewfqfwqw}) with this formula for the action of $\phi^*$
it is now straightforward to show the equality of maps
$$D=\int_p\circ \phi^*\circ p^*:C_A(\Omega_G)\to C_A(\Omega_G)\ .$$
\hB

 \section{Two-periodization --- up to isomorphism}\label{system12}

\subsection{}

Let $f\colon G\to X$   be a topological gerbe with band $U(1)$ over a locally compact stack $X$. In Definition \ref{defofd} we have constructed a natural endomorphism
$D_G\in \End(Rf_*\circ f^*)$ of degree $-2$.
To any object
$F\in D^+(\Sh_\Ab \bX)$ we  associate the inductive system
\begin{equation}\label{system10}
\cS_G(F)\colon Rf_*\circ f^*(F)\stackrel{D_G}{\leftarrow} Rf_*\circ f^*(F)[2]\stackrel{D_G}{\leftarrow}Rf_*\circ f^*(F)[4]\stackrel{D_G}{\leftarrow}\dots
\end{equation} indexed by $\{0,1,2\dots\}$.

Using the inclusion $D^+(\Sh_\Ab\bX)\to D(\Sh_\Ab\bX)$ of the bounded below into the unbounded derived category of sheaves of abelian groups on $X$ we can consider
$\cS_G(F)\in D(\Sh_\Ab \bX)^{\nat^{op}}$, where the ordered set of integers  $\nat$
is considered as a category.

\subsection{}

Using the triangulated structure of $D(\Sh_\Ab \bX)$ one can define for each object
$\cS\in D(\Sh_\Ab \bX)^{\nat^{op}}$   an object
$\holim\cS\in D(\Sh_\Ab \bX)$ which is unique up to non-canonical isomorphism
(see \cite{MR1812507}).
An explicit construction of this homotopy limit uses the extension of maps in 
$D(\Sh_\Ab \bX)$ to exact triangles by a mapping cylinder construction.
In particular, we obtain $\holim \cS_G(F)$  by the  extension to a triangle of the map $1-\hat D$ in
$$\holim \cS_G(F) \to \prod_{i\ge 0}Rf_*\circ f^*(F)[2i]\stackrel{1-\hat D}{\longrightarrow} \prod_{i\ge 0}Rf_*\circ f^*(F)[2i]\to \holim \cS_G(F)[1] \ ,$$
where
$$\hat D\colon \prod_{i\ge 0}Rf_*\circ f^*(F)[2i]\to \prod_{i\ge 0}Rf_*\circ f^*(F)[2i]$$
maps the sequence $(x_i)_{i\ge 0}$ to the sequence
$(D_Gx_{i+1})_{i\ge 0}$. 

\subsection{}

We can now define the periodization $P_G(F)\in D(\Sh_\Ab \bX)$ of an object
$F\in D^+(\Sh_\Ab \bX)$.

\begin{ddd}\label{system188}
 For $F\in D^+(\Sh_\Ab \bX)$ we define $P_G(F)\in D(\Sh_\Ab \bX)$  by
$$P_G(F):=\holim \cS_G(F)\ .$$
Note that $P_G(F)$ is well-defined up to non-canonical isomorphism.
\end{ddd}

\subsection{}\label{choicesiasg}

The operator 
$$\prod_{i\ge 0} D_G\colon \prod_{i\ge 0}Rf_*\circ f^*(F)[2i]\to (\prod_{i\ge 0}Rf_*\circ f^*(F)[2i])[-2]$$ commutes with $\hat D$ and  therefore induces a map $W\colon P_G(F)\to P_G(F)[-2]$ via an extension in the diagram
\begin{equation*}
  \begin{CD}
    P_G(F) @>W>> P_G(F)[-2]\\
    @VVV @VVV\\
   \prod_{i\ge 0}Rf_*\circ f^*(F)[2i] @>{\prod_{i\ge 0} D_G}>> \prod_{i\ge
     0}Rf_*\circ f^*(F)[2i])[-2]\\
   @VV{1-\hat D}V @VV{1-\hat D}V\\
    \prod_{i\ge 0}Rf_*\circ f^*(F)[2i] @>{\prod_{i\ge 0} D_G}>> \prod_{i\ge
      0}Rf_*\circ f^*(F)[2i])[-2]\\
    @VVV @VVV\\
    P_G(F)[1] @>{W}>> P_G(F)[1][-2]\ .
  \end{CD}
\end{equation*}
Note that such an extension exists by the axioms of a triangulated category, but it might not
be unique.  

The following proposition asserts that $P_G(F)$ is two-periodic.
\begin{prop}\label{system3}
The map
$W\colon P_G(F)\to P_G(F)[-2]$ is an isomorphism.
\end{prop}
\proof
For notational convenience, we consider the following general situation.
Let $D(A)$ be the unbounded derived category of a Grothendieck abelian
category. Note that $\Sh_\Ab(\bX)$ is such a category (see Section \ref{system11}).
We consider an object $X\in D(A)$ together with a morphism
$D\colon X\to X[-2]$. We can assume that $D$ is represented
by a map of complexes $D\colon X\to X[-2]$. We obtain the extension $1-\hat D$ to a triangle
\begin{equation}\label{system4}
Y\to \prod_{i\ge 0} X[2i]\stackrel{1-\hat D}{\to} \prod_{i\ge 0} X[2i]\to Y[1]
\end{equation}
where $Y:= \prod_{i\ge 0} X[2i]\oplus ( \prod_{i\ge 0} X[2i])[1]$ with the differential 
$$\delta:=\left(\begin{array}{cc}d&1-\hat D\\0&-d\end{array}\right)\ ,$$
where $d$ is the differential of $X$. The induced map $W\colon Y\to Y[-2]$ is given by
$$W:=\left(\begin{array}{cc}\prod_{i\ge 0}D&0\\0&\prod_{i\ge 0}D\end{array}\right)\ .$$
Let $$E\colon \prod_{i\ge 0}X[2i]\to (\prod_{i\ge 0}X[2i])[2]$$
be the shift $E(x_i)_{i\ge 0}:=(x_{i+1})_{i\ge 0}$.
Note that $E$ commutes with $1-\hat D$, too. Therefore
we obtain the extension
 $S\colon Y\to  Y[2]$ in the  diagram
$$\xymatrix{ Y\ar[r]\ar[d]^S& \prod_{i\ge 0}X[2i]\ar[d]^{E}\ar[r]^{1-\hat D}&\prod_{i\ge 0}X[2i]\ar[d]^{E}\ar[r]&Y[1]\ar[d]^S\\
Y[2] \ar[r]&( \prod_{i\ge 0}X[2i])[2]\ar[r]^{1-\hat D}&(\prod_{i\ge 0}X[2i])[2]\ar[r]& Y[1][2]}\ .$$
by the matrix 
 $$S:=\left(\begin{array}{cc}E&0\\0&E\end{array}\right)\ .$$
Proposition \ref{system3} is a consequence of the  following Lemma.
\begin{lem}\label{uwzgfwehjcfbsdac}
We have the equalities  $W\circ S=\id=S\circ W$.
\end{lem}
\begin{proof}
First observe that $\prod_{i\ge 0} D\circ E=\hat D=E\circ \prod_{i\ge 0} D$. Therefore
$W\circ S=S\circ W= \left(
  \begin{smallmatrix}
    \hat D & 0\\ 0 & \hat D
  \end{smallmatrix}\right)$.
In order to show that $W\circ S=\id$
we  show that the map
$$I:=\left(\begin{array}{cc}\hat D&0\\0&\hat D\end{array}\right)\ .$$
on $Y$ is homotopic to the identity and therefore is equal to the identity in
the derived category. This follows from   
$$1-I=\delta\circ J+J\circ \delta$$
with
$$J:=\left(\begin{array}{cc}0&0\\1&0\end{array}\right)
\ 
.$$ 
\end{proof}

\subsection{}

We continue with the notation introduced in the proof of Proposition \ref{system3}.
Applying a homological functor to the triangle (\ref{system4})  we get the long exact sequence
$$\dots\to H^*(Y)\to\prod_{i\ge 0}H^*(X[2i])\to \prod_{i\ge 0}H^*(X[2i])\to H^{*}(Y[1])\to \ .$$
If we analyze the middle map and compare it with the ordinary definition of limits in abelian categories we get the following result.
\begin{kor} \label{lim1seq}
We have an exact sequence:
$$0\to {\lim_i}^1 H^*(X[2i])[-1] \to H^*(Y)\to \lim_i H^*(X[2i]) \to 0\ .
$$
 \end{kor}

\subsection{}

Note that the construction
$$\holim\colon D(A)^{\nat^{op}}\to D(A)$$
is not a functor. The construction of the homotopy limit $\holim(S)$ for $S\in D(A)^{\nat^{op}}$  via mapping cylinders uses explicit representatives of the maps of the system $S$ and depends non-trivially on this choice.

A homotopy limit functor $\holim\colon D(A^{\nat^{op}} )\to D(A)$ can be defined as the right-derived functor 
of $\lim\colon A^{\nat^{op}} \to A$. Note that in the domain we take the derived category of the abelian category of $\nat^{op}$-diagrams in $A$ as opposed to $\nat^{op}$-diagrams in the derived category of $A$. In Section \ref{system7} we will use this idea and refine $P_G$ to a periodization functor
$$P_G\colon D^+(\Sh_\Ab\bX)\to D(\Sh_\Ab\bX)$$
which is a triangulated functor and natural in $G\to X$.
The main idea is the construction of a refinement of the diagram
(\ref{system10}) to a diagram in $D((\Sh_\Ab\bX)^{\nat^{op}})$, see \ref{system16} (the details are in fact more complicated).

\section{Calculations}

\subsection{}

In this subsection we calculate $P_G(\uF)$ in the special case, where
$G\to *$ is the (trivial) $U(1)$-gerbe over the point, and $\uF\in \Sh_\Ab\Site(*)$ is the sheaf
represented by a discrete abelian group $F$. We will calculate the abelian group $H^*(*;P_G(\uF))$.
This cohomology is two-periodic so that we only have to distinguish the even and the odd-degree case.
In the table below $\bbA^\Q_f$ denotes the group of finite adeles of $\Q$, which contains $\Q$ via the diagonal embedding.
\begin{prop}\label{system30}
We have the following table for the cohomology $H^*(*;P_G(\uF))$.
$$\begin{array}{|c|c|c|}\hline F&H^{ev}(*;P_G(\uF))&H^{odd}(*;P_G(\uF))\\\hline\hline
\Z&0&\bbA^\Q_f/\Q\\\hline
\Q&\Q&0\\\hline
\Z/n\Z&0&0\\\hline
\Q/\Z&\bbA_f^\Q&0\\\hline
  \end{array}\ .$$
\end{prop}

\subsection{}
To prove Proposition \ref{system30}, we use the exact sequence \ref{lim1seq} where
$$H^*(X)=H^*(*;Rf_*\circ f^*(\uF))\cong F\otimes \Z[ [z]]\cong F[[z]]$$
by Lemma \ref{system5} with $z$ of degree $2$. We must 
 discuss the cohomology of the complex
$$0\to \prod_{i\ge 0} F[[z]][2i]\stackrel{1-\hat D}{\to}  \prod_{i\ge 0}F[[z]][2i]\to 0\ ,$$
where $\hat D(x_i)_{i\ge 0}=(D_Gx_{i+1})_{i\ge 0}$ with $D_G=\frac{d}{dz}$.
This means that we have to study the solution theory for the system
\begin{equation}\label{system6}
x_i-\frac{d}{dz}x_{i+1}=a_i\ ,\quad  i\ge 0\ ,\quad x_i\in F[[z]]\ .
\end{equation}
 
\subsection{}\label{sec:calculate_period_cohom}
Let us start with the case $F=\Q$. Since we can divide by arbitrary integers
the operator $D_G$ is surjective and we can for any $(a_i)_{i\in\naturals}$ solve this system inductively. Therefore the cokernel $\lim_i^1 \Q[u]$
of $1-\hat D$ is trivial. A solution of the homogeneous system is uniquely determined by the choice of $x_0$ and the constant terms of the $x_i$, $i\ge 1$. Note that the constant term of $x_i$ is in degree $-2i$. It follows that
$$H^{ev}(*;P_G(\underline{\Q}))\cong   \Q\ ,\quad  H^{odd}(*;P_G(\underline{\Q}))\cong 0\ .$$
 \subsection{}
We now discuss torsion coefficients $F=\Z/n\Z$. Write $x_i=\sum x_{i,k} z^k$,
$a_i=\sum a_{i,k}z^k$ with $x_{i,k}, a_{i,k}\in \Z/n\Z$. Then we have to solve 
\begin{equation*}
  \sum_{k=0}^\infty x_{i,k}z^k -(k+1) x_{i+1,k+1} z^k=\sum_{k=0}^\infty
  a_{i,k}z^k\qquad\forall i\ge 0.
\end{equation*}
Equating coefficients  this system
decouples into finite systems 
\begin{eqnarray*}
x_{i,kn}-(kn+1)x_{i+1,kn+1}&=&a_{i,kn}\\
x_{i,kn+1}-(kn+2)x_{i+1,kn+2}&=&a_{i,kn+1}\\
&\vdots&\\
x_{i,kn+n-2}-(kn+n-1)x_{i+1,kn+n-1}&=&a_{i,kn+n+2}\\
x_{i,kn+n-1-r}+\underbrace{(kn+n)x_{i+1,kn+n}}_{=0}&=&a_{ikn+n-1} \  ,\end{eqnarray*} 
where $i,k\ge 0$.
We see that we can always solve this system uniquely by backwards induction.
We get $$H^{ev}(*;P_G(\underline{\Z/n\Z}))\cong 0\ ,\quad H^{odd}(*;P_G(\underline{\Z/n\Z}))\cong 0\ .$$
 
 \subsection{}
Let us now assume that $F=\Q/\Z$. Since this group is divisible we can solve
the system (\ref{system6}) for every $(a_i)_{i\in\naturals}$. It follows that
$$H^{odd}(*;P_G(\underline{\Q/\Z}))\cong 0\ .$$
We now discuss the solution of the homogeneous system in degree $0$.
We can choose $x_0$ arbitrary. If we have found $x_i$ for $i=0,\dots,n-1$, then we must solve
$x_{n-1}=nx_n$ in the next step. We see that $x_n$ is well-defined up to the image of $\Z/n\Z\cong n^{-1}\Z/\Z\subset \Q/\Z$.
We see that 
$H^{ev}(*;P_G(\underline{\Q/\Z}))$ admits a sequence of quotients
$$H^{ev}(*;P_G(\underline{\Q/\Z}))\to \dots \to Q^n\to Q^{n-1}\to \dots \to Q^0$$
where $Q^{n}\cong \Q/\Z$ and $Q^n\to Q^{n-1}$ is given by multiplication
with $n$ for all $n\in\naturals$.
The limit 
\begin{equation*}
\bbA^\Q_f\cong \lim_{\stackrel{\longleftarrow}{n\in \nat}}
(\Q/n!\Z)
\end{equation*}
is the ring $\bbA_f^\Q$ of finite adeles of $\Q$, and $\Q\subset \bbA^\Q_f$ is a subgroup.
We thus get
$$H^{ev}(*;P_G(\underline{\Q/\Z}))\cong \bbA^\Q_f\  .$$

 \subsection{}
Finally assume that $F=\Z$. We must again consider the system (\ref{system6}) of equations above. 
Let us discuss this system in degree $2r$. Then the relevant coefficients of
$x_i$ and $a_i$ are sequences of integers, and (writing out only these) 
$dx_{i+1}=(r+i+1)x_{i+1}$.
We see that the homogeneous equation has only the trivial solution since otherwise the integer
$x_0$ must be divisible by $n+i+1$ for all $i\ge 0$. Hence
$$H^{ev}(*;P_G(\uZ))\cong 0\ .$$ 
In order to calculate  $H^{odd}(*;P_G(\uZ))$ we consider the  exact sequence
$$0\rightarrow \Z\to \Q\to \Q/\Z\to 0\ .$$
It gives rise to an exact sequence of sheaves 
$$0\rightarrow \uZ\to \uQ\to \underline{\Q/\Z}\to 0\ .$$ and a long exact cohomology sequence.
In Section \ref{funct_per} we will construct a functorial version of $P_G$
which is a triangulated functor, and which coincides with the isomorphism class
constructed above. Using this functor, we get a triangle
$$P_G(\uZ)\to P_G(\uQ)\to P_G(\underline{\Q/\Z})\to P_G(\uZ)[1]$$ 
and therefore a long exact cohomology sequence 
$$H^{*}(*;P_G(\uZ))\to H^*(*;P_G(\uQ))\to H^*(*;P_G(\underline{\Q/\Z}))\to H^{*}(*;P_G(\uZ))[1]\ .$$
By the calculations for $\Q$ and $\Q/\Z$ we get exact sequences
$$0\to \Q\stackrel{c}{\to} \bbA_f^\Q \to H^{odd}(*;P_G(\uZ))\to 0\ ,$$
where $c$ is the canonical embedding. Therefore
$$  H^{odd}(*;P_G(\uF))\cong  \bbA^\Q_f/\Q \ .$$
\hB
 
\chapter{Functorial periodization}\label{system7}

\section{Flabby resolutions}

\subsection{}\label{nr1}

Let $\bX$ be a  site, e.g. the site of a locally compact stack. For $U\in \bX$
let $\tau:=(U_i\to U)_{i\in I}\in \cov_\bX(U)$ be a covering family.
Then we consider $V:=\bigsqcup_{i\in I} U_i\to U$. Forming iterated fiber products  we obtain a  simplicial object
$V^\cdot$ in $\bX$ with $$V^n=\underbrace{V\times_U\dots\times_UV}_{n+1 \:\:\mbox{\scriptsize factors}}\ .$$
If $F\in \Pr\bX$ is a presheaf on $X$, then we form the cosimplicial set
$C^\cdot(\tau,F):=F(V^\cdot)$.
 
\subsection{}

If $F$ is a presheaf of abelian groups, then we form the \v{C}ech complex
$\check C(\tau,F)$ which is the  chain complex associated to the cosimplicial abelian group  $C^\cdot(\tau,F)$.

If $F$ is a sheaf, then $H^0\check C(\tau,F)\cong F(U)$. We recall the following definition (see \cite[Definition 3.5.1]{MR1317816}).

\begin{ddd}[see  3.5.1, \cite{MR1317816}]\label{def:flabby}
A sheaf $F\in \Sh_\Ab\bX$ is called flabby if for all $U\in \bX$ and $\tau\in \cov_\bX(U)$ we have
$H^i\check C(\tau,F)\cong 0$ for all $i\ge 1$.
\end{ddd}
By \cite[Cor. 3.5.3]{MR1317816} a sheaf $F\in \Sh_\Ab\bX$ is flabby if and only if $R^ki(F)=0$ for all $k\ge 1$, where $i:\Sh_\Ab\bX\to \Pr_\Ab\bX$ is the inclusion of sheaves into presheaves.

As an immediate consequence of the definition a sheaf $F\in \Sh_\Ab\bX$ is flabby if and only if the restriction $F_U$  of $F$ to the site $(U)$ is flabby for all $(U\to X)\in \bX$ (see \ref{obstr1a} for the notation).

\subsection{}

Let now $X$ be a locally compact stack and $\bX$ be the site of $X$.
Occasionally, in the present paper we need the stronger notion of a flasque sheaf.
\begin{ddd}\label{flasquedef}
A sheaf $F\in \Sh_\Ab\bX$ is called flasque if for every $(U\to X)\in \bX$ and
open subset $V\subseteq U$ the restriction $F(U\to X)\to F(V\to X)$ is surjective.
\end{ddd}
In the literature, e.g. in \cite{MR1299726} or \cite{MR1481706}, 
this is used as the  definition of flabbiness. 
\begin{lem}\label{wgszguwqsws}
A flasque sheaf is flabby.
\end{lem}
\begin{proof}
For $U\in \bX$ let $\Gamma_U:\Sh_\Ab\bX\to \Ab$ be the section functor $F\mapsto \Gamma_U(F):=F(U)$.
For $V\subseteq U$ we have
$\Gamma_V(F_U)=\Gamma_V(F)$.
A sheaf  $F\in \Sh_\Ab\bX$ is flasque by definition if and only  if $F_U$ is flasque for all $U\in \bX$. But a flasque sheaf is $\Gamma_V$-acyclic
for every $V\subseteq U$ by \cite[Ch. 2, Thm. 5.4]{MR1481706} (note that in this reference our flasque is called flabby).  By \cite[Cor. 3.5.3]{MR1317816} it is flabby in the sense of \ref{def:flabby}.

This argument shows that $F_U$ is flabby for all $(U\to X)\in \bX$ and   implies that $F$ itself is flabby.
\end{proof}

We do not know if the converse of Lemma \ref{wgszguwqsws} is true.
Therefore we must be careful when using results from the literature.

\subsection{}

\begin{lem}\label{hdqoiwdwqw}
If $f:X\to Y$ is a representable map of locally compact stacks, then
a flabby sheaf is $f_*$-acyclic.
\end{lem}
\begin{proof}
Let $F\in \Sh_\Ab\bX$ be a flabby sheaf. We must show that $R^kf_*(F)=0$ for all $k\ge 1$.
We have a morphism of sites $f^\sharp:\bY\to \bX$, see \ref{obensharp}.
The functor ${}^pf_*:\Pr\bX\to \Pr\bY$ is given by ${}^pf_*F:=F\circ f^\sharp$.
It is in particular exact. Therefore we have $Rf_*\cong i^\sharp\circ {}^pf_*\circ Ri$.
Since a flabby sheaf is $i$-acyclic we conclude that $R^ki(F)=0$ for $k\ge 1$.
This implies $R^kf_*(F)=0$ for $k\ge 1$. 
\end{proof}


\subsection{}

  \begin{lem}\label{flabbypres}
If a morphism $f\colon X\to Y$ of locally compact stacks has local sections, then the functor $f^*\colon \Sh_\Ab\bY\to
\Sh_\Ab\bX$ preserves flabby sheaves.
\end{lem}
\begin{proof}
Let $F\in \Sh_\Ab\bY$ be flabby.
We consider an object  $(U\to X)\in \bX$ and a covering family $\tau \in \cov_\bX(U)$.
Then we must show that the higher cohomology groups of $\check{C}(\tau,f^*F)$ vanish.

We obtain  a covering family $f_\sharp\tau \in \cov_\bY(f_\sharp U)$, see \ref{sharp}.
Let $V^\cdot$ be the simplicial object associated to $\tau$ as in \ref{nr1}. Since
$f_\sharp$ preserves fiber products in the sense of \cite[1.2.2(ii)]{MR1317816} we see that
$f_\sharp V^\cdot$ is the simplicial object in $\bY$ associated to $f_\sharp \tau$.
The rule  $f^*F(U)\cong F(f_\sharp U)$ (see again \ref{sharp}) gives the isomorphism of
cosimplicial sets
$f^*F(V^\cdot)\cong F(f_\sharp V^\cdot)$ and hence
an isomorphism of complexes
$$\check C(\tau,f^*F)\cong \check C(f_\sharp\tau,F)\ .$$ 
Since $F$ is flabby the higher cohomology groups of the right-hand side vanish.
\end{proof}



\subsection{}
 
We now construct a canonical flabby resolution functor $$\Fl\colon \Sh_\Ab\bX\to C^+(\Sh_\Ab\bX)\ ,\quad  \id\to \Fl\ .$$ 
It associates to a $F$ a sort of Godement resolution which consists in fact of flasque sheaves.

 For a space $U$ let $(U)$ denote the site of open subsets of $U$
with the topology of open coverings. We will first construct flabby resolution
functors  
$$\Fl_U\colon \Sh_\Ab(U)\to C^+(\Sh_\Ab (U))\ ,\quad \id\to \Fl_U$$ for all  $(U\to X)\in \bX$
which are compatible with the morphisms $V\to U$ in $\bX$.
For $F\in \Sh_\Ab\bX$ we obtain a collection
of flabby resolutions $(F_U\to \Fl_U(F_U))_{U\in \bX}$, which by 
\ref{desc_sheaves_on_U} 
give rise to a resolution
$F\to \Fl(F)$.
In the following we discuss these steps in detail.

\subsection{}

Let $p_U\colon \hat U\to U$ be the identity map, where $\hat U$ is the set $U$ with the discrete topology. Let $F\in \Sh_\Ab(U)$. We set
$\Fl_U^0(F):=(p_U)_*\circ p_U^*(F)$ and let $F\to \Fl_U^0(F)$ be given by the unit
$\id\to (p_U)_*\circ p_U^*$.

\begin{lem}\label{lem:flab_exact}
  The sequence $0\to F\to (p_U)_*\circ p_U^*F$ is exact.
\end{lem}
\begin{proof}
Let $w\in U$. We must show that the induced map on stalks
$F_w\to ((p_U)_*\circ p_U^*F)_w$ is injective. This immediately follows from the description
$$((p_U)_*\circ p_U^*F)_w=\colim_{w\in W\subseteq U}
\prod_{v\in W} F_v\ .$$  
\end{proof}

\hB

\subsection{}

We now  construct $\Fl_U(F)$ inductively. Assume that we have
already  constructed $\Fl_U^0(F)\to \dots \to \Fl_U^k(F)$. Then we let
$$\Fl_U^{k+1}(F):=(p_U)_*\circ p_U^*(\coker(\Fl_U^{k-1}(F)\to \Fl_U^{k}(F))$$ and
$\Fl_U^k(F)\to \Fl_U^{k+1}(F)$ be again given by
$$\Fl_U^k(F)\to \coker(\Fl_U^{k-1}(F)\to \Fl_U^{k}))\stackrel{unit}{\to} \Fl^{k+1}_U(F)\ .$$
In this way we construct an exact complex
$$0\to F\to\Fl^0_U(F)\to \Fl^1_U(F)\to\dots\to \Fl^k_U(F)\to \dots\ .$$
All pieces of the construction are functorial. Hence, the association $F\mapsto \Fl_U(F)$ is functorial in $F$. The inclusion $F\to \Fl_U^0(F)$ gives the natural transformation $\id\to \Fl_U$.

\subsection{}

\begin{lem}\label{ufla32}
For any sheaf $F\in \Sh_\Ab(U)$ the sheaf $(p_U)_*\circ p_U^*(F)$ is flasque and flabby.
\end{lem}
\begin{proof}
For $W\subseteq U$ we have 
\begin{equation}\label{uwdhqiwduq}
(p_U)_*\circ p_U^*(F)(W)
\cong\prod_{w\in W} F_w\ .
\end{equation}
It is now obvious  that
$(p_U)_*\circ p_U^*(F)(U)\to (p_U)_*\circ p_U^*(F)(W)$
is surjective. 
A flasque sheaf is flabby by Lemma \ref{wgszguwqsws}. 
\end{proof}

%
%
%

\subsection{}\label{system14}

We now consider a sheaf $F\in \Sh_\Ab\bX$. For $(U\to \bX)$ let $F_U\in \Sh_\Ab(U)$ denote its restriction to $(U)$.
We apply the previous construction to all objects $(U\to X)\in \bX$ and the sheaves $F_U$. Then we get a collection of complexes of sheaves
$\Fl_U(F_U)$ for all $(U\to X)\in \bX$.
Let $f\colon V\to U$ be a morphism in $\bX$. We shall construct a functorial morphism
$f^*\Fl_U(F_U)\to \Fl_V(F_V)$. 

Let $G\in\Sh(U)$, $H\in\Sh(V)$, and $f^*G\to H$ be a morphism of sheaves. We
consider the diagram 
$$\xymatrix{\hat V\ar[r]^{\hat f}\ar[d]^{p_V}&\hat U\ar[d]^{p_U}\\V\ar[r]^f&U}\ .$$ 
It induces the transformation, natural in $G$ and $H$,
\begin{eqnarray*}
f^*\circ (p_U)_*\circ p_U^*(G)&\to&(p_V)_*\circ \hat f^*\circ p_U^*(G)\\
&\cong&(p_V)_*\circ p_V^*\circ f^*(G)\\
&\to&(p_V)_*\circ p_V^*(H)
\end{eqnarray*}

We now construct the map $f^*\Fl_U(F_U)\to \Fl_V(F_V)$ of complexes inductively.
Assume that we have already constructed the morphisms
$f^*(\Fl^i_U(F_U))\to \Fl_V^i(F_V)$ for all $i\le k$ compatible with the differential.
Using that $f^*$ is right exact (Lemma \ref{lem:lrexact}), we have an induced morphism
\begin{equation*}
f^*\coker(\Fl_U^{k-1}(F_U)\to \Fl_U^k(F_U))\to \coker(\Fl^{k-1}_V(F_V)\to
 \Fl_V^k(F_V)).
\end{equation*}
The construction above gives a morphism
$f^*\Fl_U^{k+1}(F_U)\to \Fl_V^{k+1}(F_V)$, again compatible with the differential of the complexes.

In this way we get the morphism
$f^*\Fl_U(F_U)\to \Fl_V(F_V)$.
By an inspection of the construction we check that for a second morphism
$g\colon W\to V$ in $\bX$ the morphisms
$g^*f^*\cF_U(F_U)\to g^*\cF_V(F_V)\to \cF_W (F_W)$ and
$(f\circ g)^*\cF_U(F_U)\to \cF_W(F_W)$
coincide.

The collections of resolutions $F_U\to \Fl_U(F_U)$, $(U\to X)\in \bX$, determines a resolution $F\to \Fl(F)$ in $C^+(\Sh_\Ab \bX)$.

\subsection{}

\begin{lem}\label{tutswassoll1}
The association $F\mapsto (F\to \Fl(F))$  is a functorial flabby resolution.
\end{lem}
\begin{proof}
The local constructions $F_U\mapsto \Fl_U(F_U)$ are functorial in $F_U$.
The connecting maps $f^*\Fl_U(F_U)\to \Fl_V(F_V)$ are compatible with this functoriality. It follows that the construction $F\to \Fl(F)$ is functorial in $F$. 

The restrictions $\Sh\bX\to \Sh(U)$ detect flabbiness and exact sequences (see  \ref{desc_sheaves_on_U}). Therefore the local statements \ref{lem:flab_exact} and \ref{ufla32}
imply that the sequence 
$0\to F\to \Fl(F)$ is a quasi-isomorphism, and that the  sheaves $\Fl^k(F)$ are flabby for all $k\ge 0$.
\end{proof}

\begin{ddd}\label{system100}
We call $F\to \Fl(F)$ the functorial flabby resolution of $F$.
\end{ddd}
Note that it actually produces resolutions by flasque sheaves.
%
%
%

\subsection{}

Let $f\colon \bX\to \bY$ be a map of locally compact stacks which has local sections.
Let $\Fl_\bX$ and $\Fl_\bY$ denote the flabby resolution functors for $\bX$
and $\bY$ according to Definition \ref{system100}.
\begin{lem}\label{ewkuwejahh}
We have a natural isomorphism of functors
$f^*\circ \Fl_\bY\cong \Fl_\bX\circ f^*$.
\end{lem}
\begin{proof}
  For $(U\to X)\in \bX$ we have by  \ref{lem:identify_star_sharp} a natural isomorphism
$f^*F_U\cong F_{f_\sharp U}$.
It gives natural isomorphisms  $\Fl_U((f^*F)_U)\cong \Fl_{f_\sharp
    U}(F_{f_\sharp U})$ and thus $\Fl _\bX(f^*F)_U\cong
  (f^*\Fl_\bY)_U$. 
Finally this collection of
  isomorphisms gives a natural isomorphism
  $$\Fl_\bX(f^*F)\cong f^*\Fl_\bY(F)\ . $$
\end{proof} 

\subsection{}

\begin{lem}\label{flat-preserv}
The functorial flabby resolution functor preserves flatness.
\end{lem}
\begin{proof}
Consider a space $U$, $p:\hat U\to U$ as above and a flat sheaf $F\in \Sh_\Ab(U)$.
Then $\coker(F\to p_*p^*(F))$ is flat as shown in the proof of \cite[Lemma 3.1.4]{MR1299726}.
This implies inductively that the sheaves $\Fl_U^k(F)$ are flat for all $k\ge 0$.
The result for the functorial flabby resolution functor on $\Sh_\Ab\bX$ now follows from the fact that the restriction functors
$\Sh_\Ab \bX\to \Sh_\Ab(U)$ detect flatness (see \ref{flatdetect}).
\end{proof}

\subsection{}

We can extend
 the flabby resolution functor \ref{system100}
to a quasi-isomorphism preserving functor
$$\Fl\colon C^+(\Sh_\Ab\bX)\to C^+(\Sh_\Ab\bX)$$ by applying $\Fl$ to a complex term-wise and forming the total complex of the resulting double  complex.

\section{A model for the push-forward}\label{iowefefwewqfqfefewf}

\subsection{}\label{afixas}

Let $f\colon G\to X$ be a morphism of  locally compact stacks which has local sections. 
 Following \cite[Sec.~2.4]{bss} we construct a nice model for the functor $Rf_*\circ f^*\colon D^+(\Sh_\Ab \bX)\to D^+(\Sh_\Ab \bX)$.  We choose an atlas $a\colon A\to G$. Then by Proposition \ref{lem:representability} the composition $f\circ a\colon  A\to G\to X$ is representable. 
Then we can define the functor
$${}^pC_A\colon C^+(\Pr_\Ab \bG)\to C^+(\Pr_\Ab \bX)$$
as in \cite[Sec.~2.4]{bss} (the subscript ${}^p$ indicates that it acts
between categories of presheaves). 

\subsection{}\label{recall-ca}

We recall the definition ${}^pC_A$.
For $(U\to X)$ consider the Cartesian diagram
$$\xymatrix{G_U\ar[d]\ar[r]&G\ar[d]^f\\U\ar[r]&X}\ .$$
Then for $F\in \Pr_\Ab\bG$ we have
\begin{equation}\label{eq:def_of_CA}
  {}^pC^k_A(F)(U\to X)=F((\underbrace{A\times_G\dots\times_GA}_{k+1 \:factors})\times_G G_U\to G)\ .
  \end{equation}
The differential 
${}^pC^k_A(F)(U\to X)\to {}^pC^{k+1}_A(F)(U\to X)$ is induced as usual as an
alternating sum 
by the projections $$(\underbrace{A\times_G\dots\times_GA}_{k+2 \:factors})\to (\underbrace{A\times_G\dots\times_GA}_{k+1 \:factors})\ .$$

We extend the functor ${}^pC_A$ to sheaves by the formula
$$C_A:=i^\sharp\circ {}^pC_A\circ i\ .$$

\subsection{}\label{reww1}

The functor $$i^\sharp\colon C^+(\Pr_\Ab\bX)\to C^+(\Sh_\Ab\bX)$$ is exact by
\ref{lechejwwcc}.
The functor ${}^p C_A$ is exact, see \cite[2.4.8]{bss}.
Since flabby sheaves are $i$-acyclic the functor $i\circ \Fl:C^+(\Sh_\Ab\bX)\to C^+(\Pr_\Ab\bX)$
preserves quasi-isomorphisms.

Therefore the composition $$i^\sharp\circ {}^pC_A\circ i\circ \Fl=C_A\circ \Fl:C^+(\Sh_\Ab\bG)\to C^+(\Sh_\Ab \bX)$$ preserves quasi-isomorphisms and descends
to the homotopy categories
\footnote{By abuse of notation we use the same symbol}
$$C_A\circ \Fl\colon hC^+(\Sh_\Ab\bG)\to hC^+(\Sh_\Ab \bX)\ .$$
After identification of the homotopy categories with the derived categories
we have  by \cite[2.41]{bss} that
$$C_A\circ \Fl\cong Rf_*\colon D^+(\Sh_\Ab \bG)\to D^+(\Sh_\Ab \bX)\ .$$

%
%

\subsection{}\label{reww2}

Since $f\colon G\to X$ has local sections the functor $f^*$ is exact.
It therefore descends to 
$$f^*\colon hC^+(\Sh_\Ab\bX)\to hC^+(\Sh_\Ab \bG)\ .$$
The composition
$$C_A\circ \Fl\circ f^*\colon hC^+(\Sh_\Ab\bX)\to hC^+(\Sh_\Ab \bX)$$
thus represents
$$Rf_*\circ f^*\colon D^+(\Sh_\Ab \bX)\to D^+(\Sh_\Ab \bX)\ .$$



\subsection{}\label{gghaaase}
We now study the dependence of $C_A$ on the choice of the atlas $A\to G$.
Let us consider a diagram 
\begin{equation}\label{afauniz}
\xymatrix{A^\prime\ar[rr]^\phi\ar[dr]^{a'}&&A\ar[dl]^a\\&G&}\ ,
\end{equation} 
where $a^\prime$ satisfies the same assumptions as $a$ (see \ref{afixas}).
The map $\phi$ induces maps
$$\xymatrix{(A^\prime\times_G\dots\times_GA^\prime)\times_GG_U \ar[rr]^{\phi^{k+1}\times \id_{G_U}} \ar[dr]&&(A\times_G\dots\times_GA)\times_GG_U \ar[dl]\\&G&}$$
and therefore
\begin{eqnarray*}
{}^pC^k_A(F)(U\to X)&=&F((\underbrace{A\times_G\dots\times_GA}_{k+1 \:factors})\times_G G_U\to G)\\
&\to&
F((\underbrace{A^\prime\times_G\dots\times_GA^\prime}_{k+1 \:factors})\times_G G_U\to G)\\
&=&
{}^pC^k_{A^\prime}(F)(U\to X)\ .
\end{eqnarray*}
This map is natural in $F$ and preserves the cosimplicial structures.
In other words, the diagram (\ref{afauniz}) induces a natural transformation
$${}^p C_{\phi}\colon  {}^pC_A\to {}^pC_{A^\prime}\ .$$
Composing with $i^\sharp$ and $i\circ \Fl$ we 
get a morphism of functors
$$C_\phi\colon C_A\circ \Fl\to C_{A^\prime}\circ \Fl\colon hC^+(\Sh_\Ab\bG)\to hC^+(\Sh_\Ab \bX)\ .$$
Both $C_A\circ \Fl$ and $C_{A^\prime}\circ \Fl$ represent $Rf_*$.
Using the explicit constructions and the proof of \cite[Lemma 2.36]{bss}
one checks that the diagram
$$\xymatrix{H^0(C_A\circ \Fl)(F)\ar[rr]^{\hspace{1cm}H^0(C_\phi)}\ar[dr]&&H^0(C_{A^\prime}\circ \Fl)(F)\ar[dl]\\&f_*(F)&}
$$
commutes.
Therefore, on the level of derived categories,  $C_\phi:C_A\circ \Fl
\to C_{A^\prime}\circ \Fl$ is the canonical isomorphism between two realizations of $Rf_*$.


\subsection{}\label{sec:let-qcolon-hto}

Let $q\colon H\to G$ be a representable morphism with local sections.
Consider the pullback diagram
$$\xymatrix{B\ar[r]^b\ar[d]^l&H\ar[d]^q\\A\ar[r]^a&G\ar[d]^f\\&X}\ .
$$
Then $b:B\to H$ is an atlas, and  we can form the functor
$C_B\colon C^+(\Pr_\Ab \bH)\to C^+(\Pr_\Ab\bX)$.

Observe that
$$B\times_H\dots \times_HB\cong (A\times_G\dots\times_GA)\times_G H\ .$$
For $(U\to X)$ consider the diagram
$$\xymatrix{H_U\ar[r]\ar[d]&H\ar[d]^q\\G_U\ar[d]\ar[r]&G\ar[d]^f\\U\ar[r]&X}\ .
$$
Observe further that
$$(B\times_H\dots \times_HB)\times_HH_U\cong (A\times_G\dots\times_GA)\times_G G_U\times_G  H\ .$$
For a presheaf $F\in \Pr \bH$ and $(V\to G)\in \bG$ we have
${}^pq_*(F)(V)=F(V\times_G H)$.
We now have the following identity
\begin{eqnarray*}
{}^pC_A^k\circ {}^pq_*(F)(U\to X)&\cong&{}^pq_*(F)(\underbrace{(A\times_G\dots\times_GA)}_{k+1 factors}\times_GG_U\to G)\\
&\cong &F(((\underbrace{A\times_G\dots\times_GA}_{k+1 factors})\times_GG_U)\times_GH\to H)\\
&\cong&
F((\underbrace{B\times_H\dots\times_HB}_{k+1 factors})\times_HH_U\to H)\\
&\cong&{}^pC^k_B(F)(U\to X)
\end{eqnarray*}
This isomorphism is functorial in $F$ and induces a natural isomorphism
$${}^pC_A\circ {}^pq_*\cong {}^p C_{q^*A}\ ,$$
where we write $q^*A:=B$.

The functor ${}^pq_*$ preserves sheaves  \cite[Lemma 2.13]{bss}. Therefore we get the identity 
$$i\circ i^\sharp \circ{}^p q_*\circ i={}^p q_*\circ i$$
 and thus an  isomorphism
\begin{equation}\label{zughj293}C_A\circ q_*\cong i^\sharp\circ {}^p C_A\circ i\circ i^\sharp\circ {}^pq_*\circ i\cong   i^\sharp\circ {}^p C_A\circ  {}^pq_*\circ i \cong i^\sharp\circ {}^pC_{q^*A}\circ i\cong C_{q^*A}\ .\end{equation}

\subsection{}

Consider a Cartesian diagram
$$\xymatrix{H\ar[d]^g\ar[r]^v&G\ar[d]\\Y\ar[r]^u&X}$$
where $u$ has local sections.  We extend the diagram to
$$\xymatrix{B\ar[d]\ar[r]&A\ar[d]\\H\ar[r]^v\ar[d]^g&G\ar[d]^f\\Y\ar[r]^u&X}\ .$$
The map $B\to H$ is again an atlas.  
\begin{lem}\label{pulcomghjdf}
We have a natural isomorphism of functors
$$ u^*\circ C_A\cong C_B\circ v^*\ .$$
\end{lem}
\begin{proof}
We first find a natural isomorphism
$${}^pu^* \circ {}^p C_A\cong {}^pC_B\circ {}^pv^*.$$
Let $(U\to Y)\in \bY$ and $F\in \Pr_\Ab\bG$. Then we have
$${}^pu^*\circ {}^pC_A(F)(U)\cong {}^pC_A(F)(u_\sharp U)\ .$$
 We have a diagram
$$\xymatrix{H_U\cong G_{u_\sharp U}\ar[r]\ar[d]&H\ar[r]^v\ar[d]^g&G\ar[d]\\U\ar[r]&Y\ar[r]^u&X}\ .$$
We calculate
\begin{eqnarray*}
(A\times_G\dots\times_GA)\times_G G_{u_\sharp U}&\cong&
(A\times_G\dots\times_GA)\times_G H\times_H G_{u_\sharp U}\\
&\cong&v_\sharp (B\times_H\dots\times_HB)\times_H H_{U}
\end{eqnarray*}
This implies that
\begin{eqnarray*}
{}^pu^*\circ C_A(F)(U)&\cong&{}^pC_A(F)(u_\sharp U)\\
&\cong&F((A\times_G\dots\times_GA)\times_G G_{u_\sharp U})\\
&\cong&F(v_\sharp ((B\times_H\dots\times_HB)\times_H H_{U}))\\
&\cong&({}^pv^*F)((B\times_H\dots\times_HB)\times_H H_{U})\\
&\cong&{}^pC_B\circ {}^pv^*(F)(U)
\end{eqnarray*}
Since $u$ and $v$ have local sections, by \ref{prexact} the functors
${}^pu^*$ and ${}^pv^*$ commute with $i\circ i^\sharp$, and this isomorphism
induces
$$u^*\circ C_A\cong C_B\circ v^*$$
(compare with the calculation (\ref{zughj293})). 
\end{proof}

\subsection{}

 The isomorphisms of Lemma \ref{pulcomghjdf} and Lemma \ref{ewkuwejahh}
 induce an isomorphism
\begin{equation}\label{uiidwqdwqdwqd}
u^*\circ C_A\circ \Fl\cong C_B\circ u^*\circ \Fl\cong  C_B\circ \Fl\circ v^*\ .
\end{equation}

On the other hand, by Lemma \ref{lem:pullpush} we have an isomorphism
$$u^*\circ Rf_*\cong Rg_*\circ v^*\ .$$
\begin{lem}\label{uiqehewqdqwdwqdqd}
The following diagram of natural isomorphisms of functors 
$$D^+(\Sh_\Ab\bG)\to D^+(\Sh_\Ab\bH)$$ commutes.
$$\xymatrix{u^*\circ C_A\circ \Fl\ar[d]^\cong\ar[r]^\cong&C_B\circ \Fl\circ v^*\ar[d]^\cong\\
u^*\circ Rf_*\ar[r]^\cong&Rg_*\circ v^*}
$$
\end{lem}
\begin{proof}
It is easy to check that this commutativity holds true on the level of zeroth cohomology sheaves.
Since all functors are the derived versions of their zeroth cohomology functors
the required commutativity follows.
\end{proof}

\begin{kor}\label{gdashdgasd}
The following diagram of natural isomorphisms commutes
$$\xymatrix{u^*\circ C_A\circ \Fl\circ f^*\ar[d]^\cong\ar[r]^\cong&C_B\circ \Fl\circ g^*\circ u^*\ar[d]^\cong\\
u^*\circ Rf_*\circ f^*\ar[r]^\cong&Rg_*\circ g^*\circ  u^*}
$$
\end{kor}

\section{Zig-zag diagrams and limits}

\subsection{}\label{system11}

We define the unbounded derived category $D(\cA)$ of an abelian category
as the homotopy category $hC(\cA)$ of  complexes (with no restrictions) in $\cA$.

\begin{ddd}\label{system104}
An abelian category $\cA$ with the following properties 
\begin{enumerate}
\item $\cA$ is cocomplete,
\item filtered colimits are exact,
\item $\cA$ has a generator, i.e.~there is an object $Z$ such that for every
  object $F$ with proper subobject $F'\subset F$, $\Hom(Z,F')\to \Hom(Z,F)$ is
  not surjective.
\end{enumerate}
is called a Grothendieck abelian category.
\end{ddd}

In this section, we will consider a Grothendieck category in which countable
products exist, e.g.~a complete Grothendieck category.
The category $\Sh_\Ab \bX$ of sheaves of abelian groups on a site $\bX$ is a complete
Grothendieck abelian category \cite[Chapter I, Thm. 3.2.1]{MR1317816}. 
\begin{lemma}\label{lem:diagramGroth}
  If $Z$ is a small category and $\cA$ is a Grothendieck abelian category in
  which countable products exists,
  then the diagram category $\cA^Z$ is again a Grothendieck abelian category
  in which countable products exist.
\end{lemma}
  This is proved in \cite[1.4.3]{MR1317816}.

\subsection{}\label{ztoou6}

We consider the category $C(\cA)$ of complexes in a Grothendieck abelian category  $\cA$.
It is known that $C(\cA)$ has a model category structure
(see \cite[Theorem 2.2]{MR1814077} where this fact is attributed to Joyal, \cite[Thm. 2.3.12]{MR1650134} for the example of the category of modules over a ring, and \cite{MR1780498} for a proof in general). This model structure is given by the following data:
\begin{enumerate}
\item The weak equivalences are the quasi-isomorphisms.
\item The cofibrations are the degree-wise injections.
\item The fibrations are defined by the right lifting property.
\end{enumerate}
By $hC(\cA)$ we denote the homotopy category of $C(\cA)$.
The category $hC(\cA)$ is triangulated with
%
the shift functor $T\colon hC(\cA)\to hC(\cA)$  given
by the shift of complexes $T(X)=X[1]$. The class of distinguished triangles is generated
by the mapping cone sequences on $C(\cA)$, 
$$\dots\to A\stackrel{f}{\to} B\to C(f)\to T(A)\dots\ .$$

The extension of a morphism in $[f]\in hC(\cA)$ with chosen representative $f\in C(\cA)$ 
to a triangle can thus naturally be  defined using the mapping cone  $C(f)$.

\subsection{}\label{ztoou}
Let $\cA$ be a Grothendieck abelian category, and consider
  a small category $Z$. Then we have an equivalence
$C(\cA)^Z\cong C(\cA^Z)$. Because $\cA^Z$ is a Grothendieck category by Lemma
  \ref{lem:diagramGroth},  
we can equip the category of $Z$-diagrams $C(\cA)^Z$
with the injective model category structure.
By translation of \ref{ztoou6} we get the following description.
\begin{enumerate}
\item The weak equivalences are the level-wise quasi-isomorphisms.
\item The cofibrations are the level-wise injections.
\item The fibrations are defined by the right lifting property.
\end{enumerate}

\subsection{}\label{hcddef}

We consider the category $U$ pictured by
$$
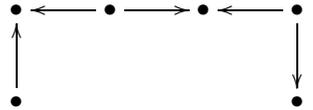
\xymatrix{\bullet&\bullet\ar[l]\ar[r]&\bullet&\bullet\ar[l]\ar[d]\\\bullet\ar[u]&&&\bullet}.$$

%

We let $\cD(\cA)\subset C^+(\cA)^U$ be the subcategory of objects of the form
\begin{equation}\label{square1}
\xymatrix{Y_0&Y_1\ar[l]^\sim\ar[r]&Y_2&Y_3\ar[l]^\sim\ar[d]\\X\ar[u]&&&X[-2]}\ 
\end{equation}
with bounded below complexes $Y_i,X$.
A morphism in the category $\cD(\cA)$ is given by maps $Y_i\to Y_i^\prime$,
$i=0,1,2,3$, and $X\to X^\prime$ which are compatible with the structure
maps. A quasi-isomorphism in this category is a morphism which is a
quasi-isomorphism level-wise.

\subsection{}\label{system15}

We let $Z$ be the category pictured by
$$
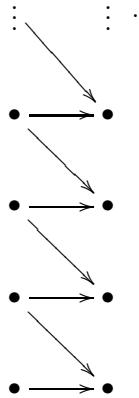
\xymatrix{\vdots\ar[dr]&\vdots\\\bullet \ar[r]\ar[dr]&\bullet\\
\bullet\ar[r]\ar[dr]&\bullet\\
\bullet\ar[r]\ar[dr]&\bullet\\
\bullet\ar[r]&\bullet}\ .
$$

Let $C(\cA)^Z$ be the category of $Z$-diagrams of complexes in $\cA$.  
We define a functor
$$R_1\colon \cD(\cA)\to  C(\cA)^Z$$
which maps the diagram (\ref{square1}) to the $Z$-diagram
$$\xymatrix{\vdots\ar[dr]&\vdots\\Y_3[4] \ar[r]^{}\ar[dr]^{}&Y_2[4]\\
Y_1[2]\ar[r]^{}\ar[dr]^{}&Y_0[2]\\
Y_3[2]\ar[r]^{}
&Y_2[2] 
}\ .
$$
The maps are induced by the shifted maps of the diagram  (\ref{square1}),
and the composition $Y_3[2k+2]\to X[2k]\to Y_0[2k]$.
The functor $R_1$  preserves quasi-isomorphisms, since those are defined
level-wise.

\subsection{}

 We now define a triangulated functor
$$\lim\colon h(C(\cA)^Z)\to hC(\cA)$$
by a direct construction 
on the level of complexes.  
Consider a $Z$-diagram $X\in C(\cA)^Z$
$$\xymatrix{C_3 \ar[r]^{c_3}\ar[dr]^{d_3}&B_3\\
C_2\ar[r]^{c_2}\ar[dr]^{d_2}&B_2\\
C_1\ar[r]^{c_1}\ar[dr]^{d_1}&B_1\\
C_0\ar[r]^{c_0}&B_0}\ .
$$
We define the morphism in $C(\cA)$
$$\phi_X\colon \prod_{i\ge 0} C_i \to \prod_{i\ge 0 }B_i$$
which maps
$(x_i)_{i\ge 0}$ to $(c_i(x_i)-d_{i+1}(x_{i+1}))_{i\ge 0}$.
Then we define $\lim(X)$ as a shifted cone of $\phi_X$:  $$\lim(X):=C(\phi_X)[-1]\in C(\cA)\ .$$  
Since quasi-isomorphisms in $C(\cA)^Z$ are defined level-wise, the
functorial  construction  $X\to \lim X$ preserves quasi-isomorphisms and thus
descends to a functor  
$$\lim\colon h(C(\cA)^Z)\to h C(\cA)\ .$$
Note that $\lim$ commutes with the shift and sum, so that it is  a triangulated functor.

\subsection{}

We now consider the composition
$\lim\circ R_1\colon \cD(\cA)\to hC(\cA)$.
The composition of the maps (or their inverses, respectively) in the diagram (\ref{square1}) gives rise to a morphism 
$D[-2]\colon X\to X[-2]$ in $hC(\cA)$.  We consider the  sequence
\begin{equation}\label{seq31o}
X^\bullet \colon X\stackrel{D}{\leftarrow} X[2]\stackrel{D[2]}{\leftarrow} X[4]\leftarrow ¸\dots\ .
\end{equation} 
in $hC(\cA)$. As already explained in \ref{system12}, for such a diagram in the triangulated category $hC(\cA)$ the homotopy limit $\holim(X^\bullet)\in hC(\cA)$ is a well-defined isomorphism class of objects. It is 
given by the mapping cone shifted by $-1$ of the morphism
$$\prod_{i\ge 0}X[2i]\to \prod_{i\ge 0}X[2i]$$
which maps $(x_i)_{i\ge 0}$ to $(x_i-D[2i]x_{i+1})_{i\ge 0}$
(see \cite[Sec.~1.6]{MR1812507}).

\begin{lem}\label{system17}
For a diagram $W\in \cD(\cA)$ of the form \eqref{square1}
we have a non-canonical isomorphism
$$\holim(X^\bullet)\cong \lim\circ R_1(W)\ .$$
\end{lem}
\begin{proof}
We use the dual statement of \cite[Lemma 1.7.1]{MR1812507}.
For $i\ge 1$ let  $C_{2i-1}=Y_3[2i]$, $C_{2i}:=Y_1[2i]$, $B_{2i-1}:=Y_2[2i]$
and $B_{2i}:=Y_{0}[2i]$. Note that we have morphisms
$v_i\colon C_{i}\to B_i$ in $C(\cA)$ which become isomorphisms in
$hC(\cA)$. Moreover, we have maps $w_{2i}\colon C_{2i}\to B_{2i-1}$ coming
from the map $Y_1\to Y_2$ of \eqref{square1}, and morphisms $w_{2i+1}\colon
C_{2i+1}\to B_{2i}$ coming from $Y_3[2]\to X\to Y_0$ of \eqref{square1}.
We consider the diagram in $hC(\cA)$, using the invertibility of $v_i$ in
$hC(\cA)$,
$$
\begin{CD}
  \prod_{i\ge 1} C_{i} @>{\prod v_i-\prod w_i}>{\phi_{R_1(W)}}> 
\prod_{i\ge 1}B_{i} 
\\
@VV{\id}V @VV{\prod_{i\ge 1} v_i^{-1}}V \\
    \prod_{i\ge 1} C_i @>>>
\prod_{i\ge 1} C_i\ ,
\end{CD}
$$
whose vertical maps are isomorphism. By definition, the mapping cone of the
upper horizontal map is $\lim\circ R_1(W)$. Because the vertical maps are
isomorphisms in $hC(\cA)$, this is isomorphic to the mapping cone of the lower horizontal map, which gives the homotopy limit of the sequence
$$ Y_3[2]\stackrel{}{\leftarrow}Y_1[2]\stackrel{}{\leftarrow}Y_3[4]\leftarrow Y_1[4]\leftarrow Y_3[6]\dots
\ .$$

We can expand this sequence to
\begin{multline}\label{longsetw}
X\stackrel{}{\leftarrow}Y_3[2]\stackrel{}{\leftarrow}
Y_2[2]\stackrel{}{\leftarrow}Y_1[2]\stackrel{}{\leftarrow}Y_0[2]\stackrel{}{\leftarrow}X[2]\\
\leftarrow Y_3[4]\stackrel{}{\leftarrow}Y_2[4]\leftarrow Y_1[4]\leftarrow Y_0[4]\leftarrow X[4]\leftarrow Y_3[6]\dots\ ,
\end{multline} 
and because the sequence (\ref{seq31o}) is just another contraction of
\eqref{longsetw}, by \cite[Lemma 1.7.1]{MR1812507} its homotopy limit
$\holim(X^\bullet)$ is then also isomorphic to 
$\lim\circ R_1(W)$. 
\end{proof}

%
%
%
%
%

\section{The functorial periodization}\label{funct_per}

\subsection{}\label{hdwidhwqdiwqdwdw}

Let $X$ be a locally compact stack. Define 
$C^+(\Sh^{flat}_\Ab\bX)\subseteq C^+(\Sh_\Ab\bX)$ to be the full subcategory of
bounded below complexes of flat sheaves. 
\begin{lem}\label{flat-inclu}
This inclusion induces an equivalence of homotopy categories
$$hC^+(\Sh^{flat}_\Ab\bX)\stackrel{\sim}{\to} hC^+(\Sh_\Ab\bX)\ .$$
\end{lem}
\begin{proof}
We first construct a functorial flat resolution functor
$$R:\Sh_\Ab\bX\to C^b(\Sh^{flat}_\Ab\bX)\ .$$
 Note that a torsion free sheaf is flat.
If $F\in \Sh_\Ab\bX$, then let $\hat F\in \Pr\bX$ denote the underlying presheaf of sets.
Let $\Z\hat F\in \Pr_\Ab\bX$ be the presheaf of free abelian groups generated by $\hat F$, and $\Z F:=i^\sharp \Z\hat F$ be its sheafification. Then we have a natural evaluation $\Z\hat F\to F$, which extends uniquely to $e:\Z F\to F$ since $F$ is a sheaf. We define
$R(F)$ to be the complex $\ker(e)\to \Z F$, where $\Z F$ is in degree zero.
The natural map $R(F)\to F$ is a quasi-isomorphism. Moreover,
$\Z F$ and its subsheaf $\ker(e)$ are torsion-free, hence flat.

We extend $R$ to a functor
$R:C^+(\Sh_\Ab\bX)\to C^+(\Sh^{flat}_\Ab\bX)$ by applying $R$ objectwise and taking the total complex of the resulting double complex.

The inclusion $C^+(\Sh^{flat}_\Ab\bX)\to C^+(\Sh_\Ab\bX)$
and $R: C^+(\Sh_\Ab\bX)\to C^+(\Sh^{flat}_\Ab\bX)$ induce mutually inverse
functors of the homotopy categories.
\end{proof}

\subsection{}

Let $f\colon G\to X$ be a topological gerbe with band $U(1)$ over a locally compact stack.
Recall the associated geometry introduced in \ref{ddef2}. {Using the functorial version we
 get the diagram
\begin{equation}\label{system40}\xymatrix{&T^2\times G\ar[dl]^{p}\ar[dr]^m &\\G\ar[dr]^f&&G\ar[dl]^f\\&X&}\end{equation}
which $2$-functorially depends on the gerbe $G\to X$.}
The map $p\colon T^2\times G\to G$ is the projection onto the second factor,
and $m:=p\circ \phi$.

\subsection{}\label{system21}

Observe that $p$ is a trivial oriented fiber bundle with fiber $T^2$.
Let $$0\to \uZ_{\Site(T^2\times G)}\to \Fl(\uZ_{\Site(T^2\times G)})$$ be the
functorial flat and flabby resolution of $\uZ_\bG$ 
constructed in \ref{system14}, see also \ref{flat-preserv} for flatness. By
$$K^\cdot\colon 0\to  K^0\to K^1\to K^2\to 0$$
we denote the truncation of $\Fl(\uZ_{\Site(T^2\times G)})$ after the
second term, i.e. with $$K^2:=\ker(\Fl^2(\uZ_{\Site(T^2\times G)})\to \Fl^3(\uZ_{\Site(T^2\times G)}))\ .$$
The complex $K^\cdot$ is still a flat and $p_*$-acyclic resolution of $\uZ_{\Site(T^2\times G)}$ (Lemma \ref{shortrest12}).  Let $$T\colon C^+(\Sh_\Ab\Site(T^2\times G))\to C^+(\Sh_\Ab\Site(T^2\times G))$$ be the functor given on objects by   
$$T_{K^\cdot}(F):=F\otimes K^\cdot\ .$$

\subsection{}\label{rcubd}
 
{We consider the commutative diagram \ref{system40}}.
Since $f\circ p\cong f\circ m$ {(recall that we actually can assume equality)}  we have by Lemma  \ref{feriueewwwwzzz} and Corollary \ref{uefhewiufuwefzzz} isomorphisms of functors
$m^*\circ f^*\cong p^*\circ f^*$ and $f_*\circ m_*\cong f_*\circ p_*$.
We fix an atlas $A\to G$ and define
$X\colon C^+(\Sh^{flat}\bX)\to C^+(\Sh\bX)$ by
$$X:=C_A\circ f^*\circ \Fl\ .$$ 
Since $f$ has local sections   we have
 $f^*\circ \Fl\cong \Fl\circ f^*$ by Lemma \ref{ewkuwejahh}.
It now follows from \ref{reww2} that $X\cong C_A\circ \Fl\circ f^*$ preserves quasi-isomorphisms.  
 It therefore descends to the homotopy categories and induces the functor $Rf_*\circ f^*$
$$D^+(\Sh_\Ab\bG)\stackrel{Lemma\  \ref{flat-inclu}}{\cong} hC^+(\Sh^{flat}_\Ab\bG)\stackrel{X}{\to} hC^+(\Sh_\Ab\bX)\cong D^+(\Sh_\Ab\bX)\ .$$

\subsection{}\label{rcubd1}

We further form $B:=m^*A\times_{T^2\times G} p^*A$.
It comes with a natural morphism $B\to m^*A$ over $T^2\times G$
which induces a transformation
$C_{m^*A}\to C_B$. 
Using the unit $\id\to m_*\circ m^*$, the inclusion
$\id\to T_{K^\cdot}$, and the isomorphisms $m^*\circ f^*\cong p^*\circ f^*$,
and using that by \ref{sec:let-qcolon-hto} $C_A\circ m_*\cong C_{m^*A}$, we
define a natural transformation
\begin{eqnarray*}
X&= &C_A\circ f^*\circ \Fl\\
&\to&C_A\circ m_*\circ m^*\circ  f^*\circ \Fl\\
&\to&C_A\circ m_*\circ T_{K^\cdot}\circ  m^* \circ f^*\circ \Fl\\
&\cong&C_{m^*A}\circ T_{K^\cdot} \circ m^* \circ f^*\circ \Fl\\
&\cong&C_{m^*A}\circ T_{K^\cdot} \circ p^* \circ f^*\circ \Fl\\
&\to&C_{m^*A}\circ \Fl\circ T_{K^\cdot} \circ p^* \circ f^*\circ \Fl\\
&\stackrel{}{\to}&C_{B}\circ \Fl\circ T_{K^\cdot} \circ p^* \circ f^*\circ \Fl\\
&=:&Y_0
\end{eqnarray*}

Using the other projection $B\to p^*A$ we define
\begin{eqnarray*}
Y_1&:=&C_{p^*A}\circ \Fl\circ T_{K^\cdot}\circ p^*\circ f^*\\
&\stackrel{\sim}{\to}&C_{B}\circ \Fl\circ T_{K^\cdot} \circ p^* \circ f^*\\
&\stackrel{\sim}{\to}&C_{B}\circ \Fl\circ T_{K^\cdot} \circ p^* \circ f^*\circ
\Fl\\
&=&Y_0 \ .
\end{eqnarray*}
Using the identity $C_{p^*A}\cong C_A\circ p_*$ we define
\begin{eqnarray*}
Y_1&=&C_{p^*A}\circ \Fl\circ T_{K^\cdot}\circ p^*\circ f^*\\
&\cong&C_{A}\circ p_* \circ \Fl\circ T_{K^\cdot}\circ p^*\circ f^*\\
&\to&C_{A}\circ \Fl\circ p_* \circ \Fl\circ T_{K^\cdot}\circ p^*\circ f^*\\
&=:&Y_2
\end{eqnarray*}
Note that $p_*\circ T_K$ is an exact functor by Lemma \ref{ggtre1} and
calculates 
$Rp_*$ by Corollary \ref{corol:calculate_Rf}. Since $p_*\circ \Fl\circ
T_K$ represents the same functor
the map $p_*\circ T_K\to p_*\circ \Fl\circ T_K$ induces a quasi-isomorphism which is preserved by $C_A\circ \Fl$.
The natural transformation
$T_{p_*K^\cdot}\xrightarrow{\sim} p_*\circ T_{K^\cdot}\circ p^*$
is an isomorphism, if applied to complexes of flat sheaves by \ref{system81}.
By Lemma \ref{prexact} the pull-back $f^*$ preserves flatness.
 
These two facts explain the quasi-isomorphisms in
 \begin{eqnarray*}
Y_3&:=&C_{A}\circ \Fl\circ T_{p_*K^\cdot}\circ f^*\\
   & \stackrel{\sim}{\to} & C_{A}\circ \Fl\circ p_*  \circ T_{K^\cdot}\circ p^*\circ f^*\\
&\stackrel{\sim}{\to}&C_{A}\circ \Fl\circ p_* \circ \Fl\circ T_{K^\cdot}\circ p^*\circ f^*\\
&=&Y_2\ .\end{eqnarray*}
Using the projection $T_{p_*K}\stackrel{[-2]}{\to}\id$ of (\ref{iofejoiwfwefwef}) we
define the natural transformation
\begin{eqnarray}\label{fret-wq}
Y_3&=&C_{A}\circ \Fl\circ T_{p_*K^\cdot}\circ f^*\\
&\to&C_{A}\circ \Fl\circ f^*[-2]\nonumber\\
&\cong&C_A\circ f^*\circ \Fl[-2]\nonumber\\
&=&X[-2]\nonumber\ .
\end{eqnarray}
Observe that all functors $Y_i$ preserve quasi-isomorphisms, using that $f^*$,
$p^*$, $C_A\circ \Fl$, $p_*\circ T_K$ (and by Lemma \ref{system81} therefore
also $T_{p_*K}$) do so.

\subsection{}\label{system16}

The construction \ref{rcubd}, \ref{rcubd1} gives a quasi-isomorphism preserving functor
$$R_0\colon C^+(\Sh^{flat}_\Ab\bX)\to \cD(\Sh_\Ab\bX)$$ (see \ref{hcddef} for the definition of the target).
By composition with the functor $R_1$ (see \ref{system15}) we get a functor
$$R:=R_1\circ R_0\colon C^+(\Sh^{flat}_\Ab\bX)\to C(\Sh_\Ab\bX)^Z\ .$$
It preserves quasi-isomorphisms and therefore descends to (again using Lemma \ref{flat-inclu})
$$R\colon D^+(\Sh_\Ab\bX)\to h(C(\Sh_\Ab\bX)^Z)\ .$$

\subsection{}\label{aabgfg}
\newcommand{\Stacks}{{\tt Stacks}}
The construction of the functor $R_0$ explicitly depends on the choice of an 
atlas $A\to G$.  These choices form a subcategory $\cZ\subset \Stacks/G$.
The choice of $A\to G$ enters the definition via the functor $C_A$.
For the moment let us indicate the dependence on $A$ in the notation and write
$R_0^A$ for the functor $R_0$ defined with the choice $A$. 

Observe, that
$A\to m^*A$, $A\to p^*A$ and $A\to m^*A\times_{T^2\times G}p^*A$ are functors
$\Stacks/G\to \Stacks/(T^2\times G)$.
The construction \ref{gghaaase} shows that for a given $F\in D^+(\Sh_\Ab\bX)$
the association $A\to R_0^A(F)$ extends to a functor 
$$R_0^{\dots}(F)\colon \cZ^{op}\to \cD(\Sh_\Ab\bX)\ .$$
The components $X\cong C_A\circ \Fl\circ f^*$ and $Y_i\cong C_*\circ \Fl\circ\dots$ (where $*\in \{A,p^*A,m^*A,m^*A\times_{T^2\times G}p^*A\}$)  all involve a flabby resolution functor in front of $C_*$. If $A\to A^\prime$ is a morphism in $\cZ$, then the transformation
$C_{A^\prime}\circ\Fl\to C_A\circ \Fl$ (or the similar transformations for
the other subscripts) produce quasi-isomorphisms by \ref{gghaaase}.


It follows that the functor $R_0^{\dots}(F)\colon \cZ^{op}\to \cD(\Sh_\Ab\bX)$
maps all morphisms to quasi-isomorphisms. We now consider the composition
$R^{\dots}(F):=R_1\circ R_0^{\dots}(F)\colon \cZ^{op}\to h(C(\Sh_\Ab\bX)^Z)$.

For two objects $A,B\in \cZ$ we consider the diagram
$$\xymatrix{&A\times B\ar[dl]^{s}\ar[dr]^t&\\A&&B}\ ,$$
where the fiber product is taken in $\Stacks/G$.
We consider the isomorphism
$$R(A,B):=R^t\circ (R^s)^{-1}\colon R^A(F)\to R^B(F)$$ in $h(C(\Sh_{\Ab}(X))^Z)$.
Using the commutativity of the squares in the diagram
$$\xymatrix{&&A\times B\times C\ar[dl]\ar[dr]\ar[d]\\&A\times B\ar[dl]\ar[dr]&A\times C\ar[dll]\ar[drr]&B\times C\ar[dr]\ar[dl]\\A&&B&&C}$$
we check that
$$R(A,B)\circ R(B,C)=R(A,C)\ .$$
This has the following consequence:.
\begin{lem}\label{lem:indofchoice}
The functor $R\colon D^+(\Sh_\Ab\bX)\to hC((\Sh_\Ab\bX)^Z)$ is independent of the choice of the atlas $A\to G$ up to canonical isomorphism. 
\end{lem}

Consider an  automorphism $\phi\colon A\to A$ in $\cZ$ and observe that it induces the identity on the level of cohomology, i.e. $H^*(R^\phi)=\id$.
It is an interesting question whether $R^\phi$ is the identity.

\subsection{}

\begin{ddd}\label{system18}
We define the periodization functor 
$$P_G:=\lim\circ R\colon D^+(\Sh_\Ab\bX)\to h(C((\Sh_\Ab\bX)^Z))\to
hC(\Sh_\Ab\bX)\ .$$  
\end{ddd}
By Lemma \ref{lem:indofchoice} it is well
defined up to canonical isomorphism.

\subsection{}\label{system19}

Let $F\in D^+(\Sh_\Ab\bX)$. By \ref{reww2} $X(F)=C_A\circ f^*\circ \Fl(F)$ represents $Rf_*\circ f^*(F)$.
The composition 
$D[-2]\colon X\to X[-2]$ of the maps (or their inverses, respectively) in
the diagram $R_0^A(F)\in \cD(\Sh_\Ab\bX)$ represents the map $D_G\colon Rf_*\circ f^*(F)\to Rf_*\circ f^*(F)[-2]$ defined in 
Definition \ref{defofd}. By Lemma \ref{system17} we see that $P_G(F)$ (according to \ref{system18})
is isomorphic to our former Definition \ref{system188} of the isomorphism class $P_G(F)$.

\section{Properties of the periodization functor}

\subsection{}

The domain and the target of $P_G$ are triangulated categories. Distinguished
triangles in both categories are all triangles which are isomorphic to
mapping cone sequences
$$\dots\to C(f)[-1]\to X\stackrel{f}{\to} Y\to C(f)\to\dots\ .$$
\begin{lem}
The functor $P_G\colon D^+(\Sh_\Ab\bX)\to hC(\Sh_\Ab\bX)$ is triangulated.
\end{lem}
\begin{proof}
We must show that it is additive, preserves the shift,  and maps distinguished triangles to distinguished triangles.
It follows from the explicit constructions that the functors $\lim$ and $R_1$ are additive and preserve the shift. The functorial flabby resolution $\Fl$ on sheaves is additive. On complexes of sheaves it is defined as the level-wise application of the flabby resolution functor composed with the total complex construction. Therefore it also commutes with the shift. All other functors  involved in the construction of $R_0$ (e.g. $C_A$, $q^*$, $T_{K^\cdot}$) are additive and commute with the shift, too.

Since the distinguished triangles in 
$D^+(\Sh_\Ab\bX)$, $h(C(\Sh_\Ab\bX)^Z)$, and $hC(\Sh_\Ab\bX)$ are defined as
triangles which are isomorphic to mapping cone sequences, and the latter only
depend on the additive structure and the shift, we see that
$\lim$ and $R$ preserve triangles.
\end{proof}

\subsection{}

\begin{lem}\label{zweipo}
For  $F\in D^+(\Sh_\Ab\bX)$ the object $P_G(F)\in hC(\Sh_\Ab\bX)$ is
two-periodic.
\end{lem}
\begin{proof}
The isomorphism $P_G(F)[2]\to P_G(F)$ is given by the isomorphism $W$ in 
  \ref{system3}.
\end{proof}
The two periodicity will be analyzed in more detail in Subsection \ref{sec:periodicity}.
%
%
%
%
%
%
%
%


\subsection{}\label{sec:naturality}

Let $u\colon Y\to X$ be a map of topological stacks which admits local sections.
Then we consider a Cartesian diagram
\begin{equation}\label{system211}
  \begin{CD} H @>{v}>> G\\
    @VVgV @VV{f}V\\
    Y @>{u}>> X .
\end{CD}
\end{equation}

\begin{lem}\label{system22}
The diagram (\ref{system211}) induces an isomorphism 
$ u^*\circ P_G\xrightarrow{\sim}  P_H\circ u^*$.
\end{lem}
\begin{proof}{
By taking the pull-back of  
(\ref{system40}) along $u$ we get  the extension}
of  the Cartesian diagram above to 
$$
\xymatrix{T^2\times H\ar@{=>}[d]^{n,q}\ar[r]^{w}&T^2\times
  G\ar@{=>}[d]^{m,p}\\H\ar[r]^v\ar[d]^g&G\ar[d]^f\\Y\ar[r]^u&X}\ .$$
{Note that there is no $2$-isomorphism between $n$ and $q$ or $m$ and $p$, respectively.}
Since $u$ has local sections the functor $u^*\colon\Sh_\Ab\bX\to \Sh_\Ab\bY$  is exact by Lemma
\ref{prexact}. It therefore  extends to functors 
 $u^*\colon \cD(\Sh_\Ab\bX)\to \cD(\Sh_\Ab\bY)$ and 
$u^*\colon  C(\Sh_\Ab \bX)^Z\to C(\Sh_\Ab \bY)^Z$ which both preserve quasi-isomorphisms.
We therefore also have corresponding functors on the derived categories which will all be denoted by $u^*$. 
In the following we are going to show that there are natural isomorphisms
\begin{enumerate}
\item $u^*\circ R_1\cong R_1\circ u^*$
\item $u^*\circ \lim \cong \lim \circ u^*$
\item $u^*\circ R_0\cong R_0\circ u^*$
\end{enumerate}
of functors on the level of homotopy categories.
 
In fact it  follows from an inspection of the construction of $R_1$ that already
$u^*\circ R_1\cong R_1\circ u^*$ on the level of functors
$\cD(\Sh_\Ab\bX)\to C(\Sh_\Ab \bY)^Z$, i.e. before descending to the homotopy category.
Assertion (1) follows.

Since $u^*\colon C(\Sh_\Ab\bX)^Z\to C(\Sh_\Ab\bY)^Z$ preserves products and mapping cones 
we again have  $u^*\circ \lim\cong \lim\circ u^*$ before going to the homotopy categories.
This implies (2).

In order to see (3), using $v$ we construct a canonical isomorphism 
$$u^*\circ R^A_0\cong R^{C}_0\circ u^*\colon C^+(\Sh^{flat}_\Ab\bX)\to \cD(\Sh_\Ab\bY)\ ,$$
where we indicate the dependence of the functor $R_0$ on the choices by a superscript as in \ref{aabgfg}. The atlas $C\to H$ is given by the 
diagram
$$\xymatrix{C\ar[d]\ar[r]&A\ar[d]\\H\ar[r]^v\ar[d]^g&G\ar[d]^f\\Y\ar[r]^u&X}\ ,$$
where the upper square is also Cartesian.

The isomorphism (3) is induced by a collection of isomorphisms indexed by the objects of the diagram $U$ (\ref{hcddef})
which induce a morphism of diagrams in $h\cD(\Sh_\Ab\bY)$.

First we have
\begin{equation}\label{eq:uXcommute}
  \begin{split}
    u^*\circ X&=u^*\circ C_A\circ f^*\circ \Fl\\
    &\cong C_C\circ v^*\circ f^*\circ \Fl\\
    &\cong C_C\circ g^*\circ u^*\circ \Fl\\
    &\cong C_C\circ g^*\circ \Fl\circ u^*\\
    &=  X\circ u^* \,
\end{split}
\end{equation}
where we use Lemma \ref{pulcomghjdf}, $v^*\circ f^*\cong g^*\circ u^*$ (see Lemma \ref{uefhewiufuwefzzz}) and the
fact that the flabby resolution functor commutes with the pull-back by $u$,
since $u$ has local sections (Lemma \ref{ewkuwejahh}). 

Let $D:=n^*C\times_{T^2\times H} q^*C$. We write
$K^\cdot_{T^2\times G}$ for the complex formerly denoted by $K^\cdot$.

Next we observe that there is a canonical isomorphism
$w^*K_{T^2\times G}^\cdot\cong K^\cdot_{T^2\times H}$. In fact 
$K_{T^2\times G}^\cdot$ and  $K^\cdot_{T^2\times H}$
are given by truncations of the complexes 
$\Fl(\uZ_{\Site(T^2\times G)})$ and $\Fl(\uZ_{\Site(T^2\times H)})$. The
isomorphism is induced by the fact that $w^*$ commutes with the flabby
resolution functor, and the isomorphism 
\begin{equation*}
w^*\uZ_{\Site(T^2\times G)}\cong \uZ_{\Site(T^2\times H)}.
\end{equation*}
This implies by Lemma \ref{tens-pres} that $w^*\circ  T_{K_{T^2\times G}^\cdot}\cong  T_{K_{T^2\times H}^\cdot}\circ w^*$. In order to increase readability of the formulas we will omit the double subscript from now on and write $T_{K^\cdot}$ for both functors.
Using this observation, Lemma \ref{pulcomghjdf},  and the other previously used isomorphisms, we get  
\begin{eqnarray*}u^*\circ Y_0&\cong&
u^*\circ C_{B}\circ \Fl\circ T_{K^\cdot} \circ p^* \circ f^*\circ \Fl\\
&\cong&C_D\circ w^*\circ \Fl\circ T_{K^\cdot} \circ p^* \circ f^*\circ \Fl\\
&\cong&C_D\circ \Fl\circ w^*\circ  T_{K^\cdot} \circ p^* \circ f^*\circ \Fl\\
&\cong&C_D\circ \Fl\circ T_{K^\cdot} \circ w^* \circ p^* \circ f^*\circ \Fl\\
&\cong&C_D\circ \Fl\circ T_{K^\cdot} \circ q^* \circ v^*\circ f^*\circ \Fl\\
&\cong&C_D\circ \Fl\circ T_{K^\cdot} \circ q^* \circ g^*\circ u^*\circ \Fl\\
&\cong&C_D\circ \Fl\circ T_{K^\cdot} \circ q^* \circ g^*\circ \Fl\circ u^*\\
&\cong&Y_0\circ u^*
\end{eqnarray*}
In a similar manner we get
\begin{eqnarray*}
u^*\circ Y_1&\cong&u^*\circ C_{p^*A}\circ \Fl\circ T_{K^\cdot}\circ p^*\circ f^*\\
&\cong&C_{q^*C}\circ w^*\circ \Fl\circ T_{K^\cdot}\circ p^*\circ f^*\\
&\vdots&\\
&\cong&Y_1\circ u^*\\
u^*\circ Y_2&\cong&Y_2\circ u^*\\
u^*\circ Y_3&\cong&Y_3\circ u^*
\end{eqnarray*}
For these isomorphisms, we use in particular Lemma \ref{lem:pullpush} to get
$v^*p_*\iso q_*w^*$, and moreover Lemma \ref{tens-pres} to get the chain of isomorphisms
\begin{equation*}
  v^*(F\tensor p_*K)\iso v^*F\tensor v^*p_*K \iso v^*F\tensor q_*w^*K\iso
  v^*F\tensor q_*K \iso T_{q_*K}(v^*F),
\end{equation*}
which gives the isomorphism $v^*\circ T_{p_*K}\iso T_{q_*K}\circ v^*$.

By a tedious check of the commutativity  of many little squares 
we see that
these maps indeed define an isomorphism of functors
$u^*\circ R^A_0\cong R^C_0\circ v^*$.
As an example of these checks, let us indicate some details of the argument
for the map
$Y_3\to X[-2]$. For $F\in D^+(\Sh_\Ab\bX)$ we have
the maps 
$\phi:Y_3(F)\to X[-2](F)$
and $\psi:Y_3(u^*F)\to X[-2](u^*F)$ given by (\ref{fret-wq}). We must show that
$$\xymatrix{u^*Y_3(F)\ar[r]^\cong\ar[d]^{u^*\phi}&Y_3(u^*F)\ar[d]^{\psi}\\u^*X[-2](F)\ar[r]^\cong&X[-2](u^*F)}$$
commutes.
 This indeed follows from the sequence of commutative diagrams
\begin{equation}\label{eq:biggcom}
  \begin{CD}
    u^*Y_3 @= u^*C_A\Fl T_{p_*K}f^* @>{T_{p_*K}\xrightarrow{[2]}\id}>> u^*C_A\Fl f^*[-2]@=      u^*X[-2]\\
&& @VV{\iso}V  @VV{\iso}V \\
&& C_Bv^*\Fl T_{p_*K}f^*       @>{T_{p_*K}\xrightarrow{[2]}\id}>>  C_B v^* \Fl f^*[-2]\\ 
&& @VV{\iso}V @VV{\iso}V \\
 &&   C_B \Fl v^* T_{p_*K}f^*     @>{T_{p_*K}\xrightarrow{[2]}\id}>> C_B \Fl v^*f^*[-2]\\
&& @VV{\iso}V @VV{\iso}V \\
Y_3u^* @= C_B \Fl T_{q_*K}g^*u^* @>{T_{q_*K}\xrightarrow{[2]}\id}>> C_B \Fl
g^*u^*[-2] @=X [-2]u^* 
  \end{CD}
\end{equation}
where for the last we use that $w$ preserves the orientation of the fiber $T^2$.
\end{proof}
The following statement directly follows from the constructions.
\begin{lem}\label{citebaleczhbwd}
The isomorphism of Lemma \ref{system22}
behaves functorially under compositions of diagrams of the form (\ref{system211}).\end{lem}

\subsection{}\label{sec:nat_of_eval}

Let $F\in D^+(\Sh_\Ab\bX)$. Recall that
$P_G(F)$ is the homotopy limit of
a $Z$-diagram consisting of sheaves $Y_0[2i]$, $Y_1[2i]$, $Y_2[2i]$, $Y_3[2i]$.
  For all $i\ge 0$
we construct an evaluation transformation
$$e_i\colon P_G(F)\to Rf_*\circ f^*(F)[2i]$$  
as the composition of the 
canonical map from the limit to $Y_3[2i+2]$  with the structure map
to $X[2i]$ and  the identification $X[2i](F)\cong  Rf_*\circ f^*[2i](F)$.
 To be precise we consider $Rf_*f^*(F)\in D(\Sh_\Ab\bX)$ via the inclusion
$D^+(\Sh_\Ab\bX)\to D(\Sh_\Ab\bX)$.
In the situation of \ref{sec:naturality} an inspection of the
proof of Lemma \ref{system22} together with Corollary \ref{gdashdgasd} shows that we have a
commutative diagram in $D(\Sh_\Ab\bX)$
\begin{equation}\label{eq:ei_natural}
  \begin{CD}
    u^* P_G(F) @>{\iso}>{v^*}> P_H(u^*F)\\
    @VV{u^*e_i}V @VV{e_i}V\\
    u^* Rf_*f^*(F)[2i] @>{\iso}>{v^*}> Rg_*g^*(u^*F)[2i] \quad .
  \end{CD}
\end{equation}

Note, however, that the morphism in the bottom line is only defined on
$D^+(\Sh_{\Ab}\bX)$ (or equivalently on its image in $D(\Sh_{\Ab}\bX)$), and
we do not know whether we can extend it to the full unbounded derived
category. Fortunately, we do not have to do this for the purposes of the
present paper.



\subsection{}

Consider the special case of the diagram (\ref{system211})
where $Y=X$, $u=\id_X$, $H=G$, and $v$ is an automorphism of the gerbe $G$.
Lemma \ref{system22} provides an  automorphism
$v^*\colon P_G\to P_G$ of periodization functors.

\subsection{}

Let us illustrate this automorphism by an example.
We consider the trivial $U(1)$-gerbe $G \to S^2$ over $S^2$ and let
$\phi\in \Aut(G/S^2)$ be classified by $1\in H^2(S^2;\Z)\cong
\Z$. 
It induces an automorphism of the cohomology
$H^*(S^2;P_G(\uF_{S^2}))$, where $\uF_{S^2}$ is the sheaf represented by a discrete abelian group $F$.
We have a Cartesian diagram
$$\xymatrix{G\ar[r]\ar[d]^g&\cB U(1)\ar[d]\\S^2\ar[r]^f&{*}}\ .$$
Since $f^*\uF_*\cong \uF_{S^2}$ we have
\begin{eqnarray*}
H^*(S^2;P_G(\uF_{S^2}))&\cong& H^*(S^2;P_G(f^*\uF_{*}))\\
&\stackrel{Lemma \:\ref{system22}}{\cong}&H^*(S^2;f^*P_{\cB U(1)}(\uF_{*}))\\
&\stackrel{Lemma \:\ref{projefoa}}{\cong}&H^*(S^2;\uZ)\otimes H^*(*;P_{\cB U(1)}(\uF_{*}))\\
&\cong&\Z[w]/(w^2)\otimes H^*(*;P_{\cB U(1)}(\uF_{*}))\ ,
\end{eqnarray*}
where $H^*(*;P_{\cB U(1)}(\uF_{*}))$ has been calculated in examples in Proposition \ref{system30}.
If $F=\Q$ or $\Q/\Z$, then $H^{ev}(*;P_{\cB U(1)}(\uF_{*}))\cong \Q$ or $\dots \cong \bbA^\Q_f$, respectively. If $F=\Z$, then 
$H^{odd}(*;P_{\cB U(1)}(\uZ_{*}))\cong \bbA^\Q_f/\Q$.
\begin{lem}
In all these cases the action of $\phi^*$ is given by 
$$\phi^*(1\otimes \lambda+ w\otimes \mu)=1\otimes \lambda+w\otimes
(\lambda+\mu)\ ,$$ 
where $\lambda,\mu\in \Q$, $\bbA^\Q_f$, or $\bbA^\Q_f/\Q$, respectively.
\end{lem}
\begin{proof}
We will use the  description of
$H^*(S^2,P_G(\uF_{S^2}))$ given in  Corollary \ref{lim1seq}.
In Lemma \ref{system5} have already calculated the automorphism on
$H^*(S^2, Rg_*g^*\uF_{S^2})\cong  F[w][[z]]/(w^2)$ induced by the diagram
$$\xymatrix{G\ar[rr]^\phi\ar[dr]_g&&G\ar[dl]^g\\&S^2}\ .$$  
It is given by $z\mapsto z+w$, $w\mapsto w$.
The operation induced by $D_G$ is $\frac{d}{dz}$, and
%
  the periodized cohomology  is given as the
kernel (in the cases $F=\Q$ and $F=\Q/\Z$) or cokernel (in the case $F=\Z$) of $\prod_{i\ge 0} \id[2i] -\prod_{i\ge 0} D_G[2i]$ on $\prod_{i\ge 0} F[w][[z]]/(w^2)[2i]$.
Recall from \ref{sec:calculate_period_cohom}
that the 
class 
$a\in H^0(S^2,P_G(\uQ_{S^2}))\cong \Q[w]/(w^2)$ is
represented by $(a,az,az^2/2,\dots, az^k/k!\dots )$, which is mapped by $\phi^*$ to
$(a,a(w+z),a(w+z)^2/2,\dots)$.
We must  read off a representative of this class in the form above.
If $a=w$ then $w(w+z)^k/k!=wz^k/k!$ and therefore
$\phi^*w=w$. On the other hand, if $a=1$, then 
$a(w+z)^k/k!=z^k/k!+w z^{k-1}/(k-1)!$, so that $\phi^*(1)=1+w$.

Exactly the same argument applies if $F=\Q/\Z$. Finally, the cohomology with
coefficients $F=\Z$ is the cokernel (up to shift of degree) of the map induced
by the inclusion $\Q\into \bbA_f^\Q$, which implies the assertion also for
$F=\Z$. 
\end{proof}

\section{Periodicity}\label{sec:periodicity}

\subsection{}

We consider a topological $U(1)$-gerbe $f\colon G\to X$ over a locally compact stack.
Let $F\in D^+(\Sh_\Ab\bX)$.
In Lemma \ref{zweipo} we have argued that $P_G(F)\in D(\Sh_\Ab\bX)$ is two-periodic.
The periodicity is implemented by a certain isomorphism
$W:P_G(F)[2]\to P_G(F)$ which may depend on additional choices, see also the discussion in \ref{choicesiasg}. In the present subsection we show that there is a canonical two-periodicity isomorphism.

\subsection{}\label{ewjdgednaa}

The gerbe $G\to X$ gives rise {in a $2$-functorial way} to
the diagram  (see \ref{ddef2} for details)
\begin{equation}\label{system53}
\xymatrix{\tilde G\ar[dr]^s\ar[d]^r\ar[rr]^\phi&&\tilde G\ar[dl]^s\ar[d]^r\\G\ar[dr]_f&X\times T^2\ar[d]^p&G\ar[dl]^f\\&X&
}\ .
\end{equation}
This diagram induces the desired periodization isomorphism
as the following composition of natural transformations
\begin{multline}\label{prechgh}
W\colon P_G(F)\stackrel{\text{unit}}{\rightarrow} Rp_*p^*P_G(F)\stackrel{\text{Lemma
    \ref{system22}}}{\rightarrow}  Rp_*P_{\tilde
    G}(p^*F)\\
\xrightarrow{\phi^*}Rp_*P_{\tilde G}(p^*F) \cong 
Rp_*p^*P_{G}(F)\xrightarrow{\int_p} P_G(F)[-2]\ .\end{multline}

\begin{prop}\label{dhewud82d}
The transformation (\ref{prechgh})
$$W\colon P_G(F)\to P_G(F)[-2]$$ 
is a canonical choice for the isomorphism in  
Proposition \ref{system3}.
\end{prop}

\subsection{}
To start the proof of Proposition \ref{dhewud82d},
recall the definition
$$D_G\colon Rf_*f^*(F)\to Rf_*f^*(F)[-2]$$ as the composition
$$Rf_*f^*(F)\xrightarrow{\text{unit}} Rf_*Rr_*R\phi_*\phi^* r^*f^*(F)\stackrel{!}{\cong} Rf_*Rr_*r^*f^*(F)\xrightarrow{\int_r} Rf_*f^*(F)[-2]\ ,$$
where at the marked isomorphism $"!"$ we use the natural isomorphisms
\ref{keykey} and \ref{uefhewiufuwefzzz} associated to the identity  $f\circ r= f\circ r\circ\phi$

Recall from  \ref{sec:nat_of_eval}
the definition
of the natural evaluation transformation $e_i\colon P_G(F)\to Rf_*f^*(F)[2i]$ for all $i\ge 0$.

\begin{lem}\label{uzasdda}
The following diagram commutes:
$$\xymatrix{P_G(F)\ar[d]^{e_{i+1}}\ar[r]^W&P_G(F)\ar[d]^{e_i}\\
Rf_*f^*(F)[2i+2]\ar[r]^{D_G}&Rf_*f^*(F)[2i]}\ .$$
\end{lem}
\begin{proof}
We split this square in parts.  First we observe that in $D(\Sh_\Ab \bX)$
$$
\begin{CD}
  P_G(F)
  @>{\text{unit}}>> Rp_*p^*P_G(F) @>{Rp_*r^*}>\cong> Rp_*P_{\tilde
  G}(p^*F)\\
   @VV{e_{i+1}}V @VV{Rp_*p^*e_{i+1}}V @VV{Rp_*e_{i+1}}V \\
    Rf_*f^*(F)[2i+2] @>{\text{unit}}>> Rp_*p^*
    Rf_*f^*(F)[2i+2] @>{Rp_*r^*}>\cong> Rp_*
    Rs_*s^*p^*(F)[2i+2]\\
    @VV=V && @VV{\cong}V \\
    Rf_*f^*(F)[2i+2] @>{Rf_*f^*\text{unit}}>> Rf_*f^*Rp_*p^*(F)  @>{\iso}>>
  Rf_*Rr_*r^*f^*(F)[2i+2] 
\end{CD}
$$
commutes (use Lemma   \ref{lem:pullpush} for the upper left and the lower and \ref{sec:nat_of_eval} for the upper right rectangle).

In the next step we observe that
$${\small
\begin{CD}
 Rp_*P_{\tilde G}(p^*F) @>{\id}>>  Rp_*P_{\tilde G}(p^*F)
      @>{Rp_*\phi^*}>{\cong}>   
 Rp_*P_{\tilde
      G}(p^*F)\\
 @VV{Rp_*e_{i+1}}V && @VV{Rp_*e_{i+1}}V \\
 Rp_*Rs_*s^*p^*(F)[2i+2] @>{\text{unit}}>> Rp_*Rs_*R\phi_*\phi^*s^*p^*(F)[2i+2]
      @>\cong>> Rp_*Rs_*s^*p^*(F)[2i+2] \\
 @VV{\cong}V @VV{\cong}V @VV{\cong}V\\
Rf_*Rr_*r^*f^*(F)[2i+2] @>{\text{unit}}>> Rf_*Rr_*R\phi_*\phi^*r^*f^*(F)[2i+2]
      @>\cong>> Rf_*Rr_*r^*f^*(F)[2i+2] 
\end{CD}
}$$
commutes, where we use for the upper rectangle again \ref{sec:nat_of_eval},  and $p\circ s\circ \phi=p\circ s$, $p\circ s=f\circ r$,  $f\circ r\circ \phi=f\circ r$ and
Lemma \ref{lem:pullpush} for the remaining squares. 

In the last step we observe the commutativity of
\begin{equation*}
  \begin{CD}
    Rp_*P_{\tilde G}(p^*F) @>{(r^*)^{-1}}>{\iso}> Rp_*p^*P_{G}(F) @>{\int_p}>>
    P_G(F)[-2]\stackrel{T^{-2}}{\cong} P_G(F)\\ 
    @VV{Rp_*e_{i+1}}V @VV{Rp_*p^*e_{i+1}}V @VV{T^{-2}e_{i+1}\cong e_i}V\\
    Rp_*Rs_*s^*p^*(F)[2i+2] @>{(r^*)^{-1}}>{\iso}> Rp_*p^*Rf_*f^*(F)[2i+2] @>{\int_p}>>
    Rf_*f^*(F)[2i]\\
    @VV{\iso}V && @VV{=}V\\
    Rf_*Rr_*r^*f^*(F)[2i+2] & @>{Rf_*(\int_r)}>> & Rf_*f^*(F)[2i].
  \end{CD}
\end{equation*}
Again, for the commutativity of the upper left rectangle we use
\eqref{eq:ei_natural} of \ref{sec:nat_of_eval}. For the upper right corner we
%
%
use the fact that $\int_p$ is a natural
transformation between the functors $Rp_*p^*$ and $\id$ on $D(\Sh_\Ab\bX)$.  For the lower rectangle we use Lemma \ref{system300}.
\end{proof}

\subsection{}

We now finish the proof of Proposition \ref{dhewud82d}.
We have an exact triangle
$$\dots\to P_G(F)\stackrel{\prod_{i\ge 0}e_i}{\to} \prod_{i\ge 0} Rf_*f^*(F)[2i]\stackrel{\alpha}{\to }\prod_{i\ge 0} Rf_*f^*(F)[2i]\stackrel{[1]}{\to}\dots$$
where (using the language of elements) the map $\alpha$ is given by
$$\alpha(x_i)_{i\ge 0}=(x_i-D_Gx_{i+1})_{i\ge 0} .$$ 
By Lemma \ref{uzasdda} we have a morphism of exact triangles
$$\xymatrix{ P_G(F)\ar[d]^W\ar[r]^{\hspace{-0.9cm}\prod_{i\ge 0}e_i} & \prod_{i\ge 0} Rf_*f^*(F)[2i]\ar[d]^\beta\ar[r]^\alpha&\prod_{i\ge 0} Rf_*f^*(F)[2i]\ar[d]^\beta\\
P_G(F)[-2]\ar[r]^{\hspace{-0.9cm}\prod_{i\ge 0}e_i}& \prod_{i\ge 0} Rf_*f^*(F)[2i-2]\ar[r]^\alpha&\prod_{i\ge 0} Rf_*f^*(F)[2i-2]
}\ ,$$
where the map $\beta$ is given by
$\beta(x_i)_{i\ge 0}:=(D_Gx_i)_{i\ge 0}$. 
In Lemma \ref{uwzgfwehjcfbsdac} we have shown that $W$ is an isomorphism.
\hB

\chapter{$T$-duality}\label{system4000}

\section{The universal $T$-duality diagram}

\subsection{}

Topological $T$-duality intends to model the underlying topology of string theoretic $T$-duality on the level of targets and quantum field theory.  In the special case of targets modeled by a gerbe on top of a $T^n$-principal bundle over a space,  topological $T$-duality is by now a well-defined mathematical concept, see \cite{math.AT/0701428}, \cite{math.GT/0501487} and the literature cited therein. In the case of $T$-principal bundles it was extended to orbifolds in \cite{MR2246781}.
In the present paper we propose a definition of $T$-duality in the case of $T$-bundles over arbitrary stacks.
This framework includes arbitrary $T$-actions on spaces. The  special case
of an almost free action (i.e.~every orbit is either free or a fixed point)
has been treated with completely different methods in \cite{pande-2006}.

\subsection{}

The notion of a $T$-duality diagram has first been introduced in \cite{math.GT/0501487}. 
In the present paper we first produce a universal $T$-duality diagram over the stack $\cB U(1)=[*/U(1)]$.
Then we proceed to define a $T$-duality diagram over a general stack as one which is locally isomorphic to the universal one.

\subsection{}\label{system1003}

The universal $T$-duality diagram is a diagram of stacks
\begin{equation}\label{univtdza1}
\xymatrix{&p_{univ}^*G_{univ}\ar[dl]\ar[dr]\ar[rr]^{u_{univ}}&&\hat p_{univ}^* \hat G_{univ}\ar[dl]\ar[dr]&\\G_{univ}\ar[dr]^{f_{univ}}&&F_{univ}\ar[dl]^{p_{univ}}\ar[dr]^{\hat p_{univ}}&&\hat G_{univ}\ar[dl]^{\hat f_{univ}}\\&E_{univ}\ar[dr]^{\pi_{univ}}&&\hat E_{univ}\ar[dl]^{\hat \pi_{univ}}&\\&&B_{univ}&&}\ .
\end{equation}
In the following we explain the stacks and the maps.
\begin{itemize}
\item $B_{univ}:=\cB U(1)$
\item $E_{univ}:=*$ and $\pi_{univ}$ is the map which classifies the trivial $U(1)$-bundle over the point $*$.
\item $G_{univ}:=\cB U(1)$, and $f_{univ}$ is the unique map.
\item $\hat E_{univ}:=\cB U(1)\times U(1)$, and $\hat \pi_{univ}$ is the projection onto the first factor.
\item $\hat f_{univ}\colon \hat G_{univ}\to \hat E_{univ}$ is a gerbe with band $U(1)$ classified by $z\otimes v\in H^2(\cB U(1);\Z)\otimes H^1(U(1);\Z)\cong H^3(\cB U(1)\times U(1);\Z)$, where
$z\in H^2(\cB U(1);\Z)$ and $v\in H^1(U(1);\Z)$ are the standard generators.
\item $F_{univ}:=E_{univ}\times_{B_{univ}}\hat E_{univ}\cong U(1)$, and $p_{univ},\hat p_{univ}$ are the canonical projections.
\item Since $H^2(F_{univ};\Z)\cong 0\cong H^3(F_{univ};\Z)$, the pull-back
$\hat p^*_{univ}\hat G_{univ}$ can be identified with the trivial gerbe
$p_{univ}^*G_{univ}\cong U(1)\times\cB U(1)$ by a unique isomorphism class of maps
represented by $u_{univ}$. 
\end{itemize}
Let us fix once and for all a universal $T$-duality diagram (i.e. a choice of $u_{univ}$ in its isomorphism class {and $2$-isomorphisms filling the faces}).

\subsection{}

Let $B$ be a topological stack and consider a diagram 
\begin{equation}\label{tzfrzrzr}
\xymatrix{&p^*G\ar[dl]\ar[dr]\ar[rr]^u&&\hat p^* \hat G\ar[dl]\ar[dr]&\\G\ar[dr]^f&&F\ar[dl]^p\ar[dr]^{\hat p}&&\hat G\ar[dl]^{\hat f}\\&E\ar[dr]^\pi&&\hat E\ar[dl]^{\hat \pi}&\\&&B&&}
\end{equation}
of topological stacks where the squares are Cartesian, $f\colon G\to E$ and $\hat f\colon \hat G\to \hat E$ are topological  $U(1)$-gerbes,  and $u$ is an isomorphism of gerbes over $F$.

An isomorphism between two such diagrams over $B$ is first of all a large commutative diagram in stacks, but we furthermore require that the horizontal morphisms are morphisms of $U(1)$-banded gerbes in all places 
where this condition makes sense.

\begin{ddd}\label{eruihfrvc}
The diagram (\ref{tzfrzrzr}) is called a $T$-duality diagram if for every object $(U\to B)\in \bB$ there exists a covering
$(U_i\to U)_{i\in I}\in \cov_\bB(U)$ such that for all $i\in I$ the pull-back of
the diagram (\ref{tzfrzrzr}) along the map $U_i\to U \to B$ is isomorphic to the pull-back of the  universal $T$-duality diagram (\ref{univtdza1}) along a map $U_i\to B_{univ}$.
\end{ddd}

\subsection{}

In the following we describe the concept of $T$-duality. 
Let $B$ be a topological stack. A pair $(E,G)$ over $B$ consists of a $T$-principal bundle
$\pi\colon E\to B$ and a $U(1)$-gerbe $f\colon G\to E$. 
\begin{ddd}\label{def:admits_dual}
We say that a pair $(E,G)$ admits a $T$-dual, if it appears as a part of a $T$-duality diagram
\ref{tzfrzrzr}. In this case the pair $(\hat E,\hat G)$ is called a $T$-dual of $(E,G)$.
\end{ddd}

This is our proposal for the mathematical concept of $T$-duality for pairs of $T$-principal bundles and gerbes.
Using the $T^n$-bundle variant of the universal $T$-duality diagram one can easily generalize this definition to the higher-dimensional case. But note that, in contrast to the case of
one-dimensional fibers, a unique isomorphism $u_{univ}$ does not exist for
$T^n$ if one uses the exact parallel setup. 
This explains why suitable
modifications are necessary in \cite{math.GT/0501487}. In particular, the universal base space is not simply
the $n$-fold product of copies of  $B_{univ}$ used in the one-dimensional case.

\subsection{}

In the following we show that the concept of topological $T$-duality as defined above really coincides with the former definitions.
\begin{lem}\label{lem:identify_T_duality}
 Definitions \ref{eruihfrvc} and \ref{def:admits_dual} reduce to the notion of $T$-duality as used in
  \cite{math.GT/0501487}, \cite{MR2130624}, if $B$ is a locally acyclic space.
\end{lem}
\begin{proof}
By Definition \ref{eruihfrvc} a $T$-duality triple over a space $B$ is given by the following data:
\begin{enumerate}
\item  locally trivial $U(1)$-principal bundles $E,\hat E$ over $B$,
\item $U(1)$-banded gerbes $G$, $\hat G$ over $E$ or $\hat E$,
  respectively,
\item  an isomorphism $u$ between the pullbacks of
  $G$ and $\hat G$ to the correspondence space $E\times_B \hat E$.
\end{enumerate} 
Every point $b\in B$ admits an acyclic neighborhood $b\in U\subseteq B$.
The bundles $E$ and $\hat E$ are trivial over $U$, i.e. we have $E_{|U}\cong U\times U(1)\cong \hat E_{|U}$.
Since  $H^3(U\times U(1);\Z)\cong 0$, the restrictions of the gerbes $G_{|E_{|U}}$ and $\hat G_{|\hat E_{|U}}$ are trivial, too. The Definition  \ref{eruihfrvc} requires that the isomorphism of trivial gerbes $u_{|E_{|U}\times_U\hat E_{|U}}$ is 
 classified by the
  generator of $H^2(E_{|U}\times_U\hat E_{|U};\Z)$ (note that $E_{|U}\times_U\hat E_{|U}\cong U\times U(1)\times U(1)$).
This reformulation of the definition of a $T$-duality triple over a locally acyclic space $B$
is exactly the definition of a $T$-duality triple in  \cite{math.GT/0501487}.

In the approach of  \cite{MR2130624}  to $T$-duality
we start with a pair $(E,G)$. We characterize $T$-dual pairs by topological conditions.
We then analyze the classifying space of pairs and observe that the universal pair has a 
unique $T$-dual pair which gives rise to the $T$-duality transformation.

It turns out that the classifying space of pairs in \cite{MR2130624} is equivalent to the classifying space of $T$-duality triples in \cite{math.GT/0501487}, and that the universal pair and its dual
are  parts of the  universal $T$-duality triple. This shows that the approaches of  \cite{MR2130624}
and \cite{math.GT/0501487} are equivalent.
\end{proof}

\section{$T$-duality and periodization diagrams}

\subsection{}

Recall that the construction of the periodization functor
$P_G$ was based on the diagrams introduced in \ref{ddef2}. In the present subsection we
relate these diagrams to $T$-duality. 

\subsection{}

The double of the universal $T$-duality diagram (\ref{univtdza1}) is (by
definition) the 
big universal periodization diagram
\begin{equation}\label{univtdza11}
{\small \xymatrix{&\pr_0^*p_{univ}^*G_{univ}\ar[dl]^{\tilde
      \pr_0}\ar[dr]\ar[r]^{\pr_0^*u_{univ}}& \pr_{\hat E}^*\hat
    G_{univ}\ar[r]^{\pr_1^*u_{univ}^{-1}}\ar[d]&\pr_1^*p_{univ}^* G_{univ}
    \ar[dl] \ar[dr]_{\tilde
      \pr_1}&\\p_{univ}^*G_{univ}\ar[dd]^{f_{univ}^*p_{univ}}\ar[dr]^{p_{univ}^*{f_{univ}}}&&F_{univ}\times_{\hat
      E_{univ}}F_{univ}\ar[dl]^{\pr_0}\ar[dr]^{\pr_1}&&p_{univ}^*
    G_{univ}\ar[dd]_{f_{univ}^*p_{univ}}\ar[dl]_{p_{univ}^*{f_{univ}}}\\&F_{univ}\ar[dr]^{p_{univ}}&&F_{univ}\ar[dl]^{p_{univ}}&\\G_{univ}\ar[rr]^{f_{univ}}&&E_{univ}&&G_{univ}\ar[ll]^{f_{univ}}}
}
\end{equation}
Note that all  squares are Cartesian, with the exception 
of the central square
$$\xymatrix{&F_{univ}\times_{\hat E_{univ}}  F_{univ}\ar[dl]\ar[dr]&\\F_{univ}\ar[dr]&&F_{univ}\ar[dl]\\&E_{univ}&}$$
which does not commute.
 The same remark applies to similar diagrams we introduce later.

\subsection{}

We form the diagram\footnote{{This diagram does not commute. It is a short-hand for a square of the form (\ref{system40}) with a $2$-isomorphism
between $f_{univ}\circ q_{univ}$ and $f_{univ}\circ m_{univ}$. We will adopt a similar convention for other diagrams written in this short-hand form below.}}

\begin{equation}\label{ajksdhbxaiudkj}
\xymatrix{\pr_0^*p^*_{univ}G_{univ}\ar@/^1pc/[r]^{q_{univ}}\ar@/_1pc/[r]_{m_{univ}}&G_{univ}\ar[r]^{f_{univ}}&E_{univ}}\ ,
\end{equation}
where $$m_{univ}:= f_{univ}^*p_{univ}\circ \tilde \pr_1\circ \pr_1^* u_{univ}^{-1}\circ \pr_0^*u_{univ}\ ,\quad 
q_{univ}:= f_{univ}^*p_{univ}\circ \tilde \pr_0\ .$$
\begin{ddd}
The diagram (\ref{ajksdhbxaiudkj}) is called the small universal periodization diagram.
\end{ddd}

\subsection{}

Let $f\colon G\to X$ be a topological gerbe with band $U(1)$ over a stack $X$. Then we consider the pull-back of the small universal periodization diagram to $X$ via the projection $r\colon X\to E_{univ}\cong *$. We form the tensor product with the gerbe $G$ (see \cite[6.1.9]{math.AT/0701428} for some details on such tensor products) and obtain the diagram 
\begin{equation}\label{jhsdjkhwud}
\xymatrix{\tilde H\ar@/^1pc/[r]^{q}\ar@/_1pc/[r]_{m}&H \ar[r]^f& X}\ ,
\end{equation}
where 
$$\tilde H:= \pr_X^*G\otimes \pr^*_{F_{univ}\times_{\hat E_{univ}}F_{univ}}\pr^*_0 p_{univ}^*G_{univ}\ , \quad H:=G\otimes r^*G_{univ}\ ,$$
$$\pr_X\colon X\times F_{univ}\times_{\hat E_{univ}}F_{univ}\to X\ ,$$ 
$$\pr_{F_{univ}\times_{\hat E_{univ}}F_{univ}}\colon X\times F_{univ}\times_{
\hat E_{univ}}F_{univ}\to F_{univ}\times_{\hat E_{univ}} F_{univ}$$ are the
projections, and $m$, $q$ are induced by the corresponding universal maps
$m_{univ}$ or $q_{univ}$, respectively.
\begin{ddd}
The diagram (\ref{jhsdjkhwud}) is called the small periodization diagram of $G\to X$.
\end{ddd}
In fact we have defined a {$2$-functor}  from $\gerbes/X$ to a $2$-category of such small
periodization diagrams.  
Using the fact that $G_{univ} = \cB U(1)$ we have a canonical identification
$H\cong G$. Furthermore, $F_{univ}\times_{
\hat E_{univ}}F_{univ}\cong T^2$, and we can identify
$\tilde H \to X\times F_{univ}\times_{\hat E_{univ}}F_{univ}$ with $G\times T^2 \to X\times T^2$.
\begin{lem}
With these identifications  the small periodization diagram
(\ref{jhsdjkhwud}) is isomorphic to  the diagram (\ref{system40})
used in the definition of $P_G$.
\end{lem}
\begin{proof}
  This follows directly from the definitions of these maps.
\end{proof}

\subsection{}\label{system52}

The $T$-duality diagram (\ref{tzfrzrzr}) gives rise to the big double $T$-duality diagram
 \begin{equation}\label{univtdzfe1}
\xymatrix{&\pr_0^*p^*G\ar[dl]^{\tilde \pr_0}\ar[dr]\ar[r]^{\pr_0^*u}& \pr_{\hat E}^*\hat G\ar[r]^{\pr_1^*u^{-1}}\ar[d]&\pr_1^*p^* G \ar[dl] \ar[dr]_{\tilde \pr_1}&\\p^*G\ar[dd]^{f^*p}\ar[dr]^{p^*f}&&F\times_{\hat E}F\ar[dl]^{\pr_0}\ar[dr]^{\pr_1}&&p^* G\ar[dd]^{f^*p}\ar[dl]_{p^*f}\\&F\ar[dr]^p&&F\ar[dl]^{p}&\\G\ar[rr]^f&&E&&G\ar[ll]^f}
\end{equation}
Note that the middle square does not commute.
We have 
$$F\times_{\hat E}F\cong (E\times_B\hat E)\times_{\hat E}(\hat E\times_B E)\cong 
E\times_B\hat E\times_B E\stackrel{\sim}{\leftarrow} E\times_B\hat E\times U(1)\ ,$$
where the last arrow is given by $(e,\hat e,eu)\leftarrow (e,\hat e,u)$.
Under this identification $\pr_0(e,\hat e,u)=(e,\hat e)$ and $\pr_1(e,\hat e,u)=(eu,\hat e)$.
We can correct this non-commutativity as follows.
Let $c:F\times_{\hat E} F\to F\times_{\hat E} F$ be the isomorphism, which under the above
identification is given by
$c(e,\hat e,u):= (eu^{-1},\hat e,u)$.
Note that $\pr_1\circ c=\pr_0$.
Furthermore note that $\pr_{\hat E}=\pr_{\hat E}\circ c:F\times_{\hat E} F\to \hat E$.
Therefore we get a canonical morphism $\hat c$ satisfying
$\overline{\pr_{\hat E}}=\overline{\pr_{\hat E}}\circ \hat c$
in the diagram
$$\xymatrix{\pr^*_{\hat E}\hat G\ar[d]\ar[r]^{\hat c}&\pr^*_{\hat E}\hat G\ar[d]\ar[r]^{\overline{\pr_{\hat E}}}&\hat G\ar[d]\\F\times_{\hat E} F\ar[r]^c&F\times_{\hat E} F\ar[r]^{\pr_{\hat E}}&\hat E}\ .$$
If we plug this in the big double $T$-duality diagram, then we get the big commutative $T$-duality diagram diagram
\begin{equation}\label{univtdzfe11wdqdqw}
\xymatrix{&\pr_0^*p^*G\ar[dl]^{\tilde \pr_0}\ar[dr]\ar[r]^{\pr_0^*u}& \pr_{\hat E}^*\hat G  \ar[d]                         \ar[r]^{\hat c}         &\pr^*_{\hat E} \hat G                          \ar[r]^{\pr_1^*u^{-1}}\ar[d]&\pr_1^*p^* G \ar[dl] \ar[dr]_{\tilde \pr_1}&\\p^*G\ar[dd]^{f^*p}\ar[dr]^{p^*f}&&F\times_{\hat E}F\ar[dl]^{\pr_0}                      \ar[r]^c& F\times_{\hat E}F              \ar[dr]^{\pr_1}                   &&p^* G\ar[dd]^{f^*p}\ar[dl]_{p^*f}\\&F\ar[dr]^p&&&F\ar[dl]^{p}&\\G\ar[rr]^f&&E  \ar@{=}[r]  &E                                                                                    &&G\ar[ll]^f}
\end{equation}

>From this we derive the diagram 
\begin{equation}\label{sdsdh}
\xymatrix{\pr_0^*p^*G\ar@/^1pc/[r]^{q_T}\ar@/_1pc/[r]_{m_T}&G\ar[r]^f&E}\ ,
\end{equation}
where $$q_T:=f^*p\circ \tilde \pr_0\ ,\quad m_T:= f^*p\circ \tilde \pr_1\circ \pr_1^*u^{-1}\circ \hat c\circ \pr_0^* u\ .$$
\begin{ddd}
The diagram (\ref{sdsdh}) is called the small double $T$-duality diagram associated to  (\ref{tzfrzrzr}).
\end{ddd}
\subsection{}
The following fact is an immediate consequence of the definitions.
\begin{prop}\label{eukdnwe}
The small double $T$-duality diagram (\ref{sdsdh})  is locally isomorphic to the small periodization diagram  (\ref{jhsdjkhwud})
 of $G\to E$.
\end{prop}

\section{Twisted cohomology and the $T$-duality transformation}\subsection{}
\newcommand{\fp}{\mathfrak{p}}\newcommand{\fq}{\mathfrak{q}}\newcommand{\fu}{\mathfrak{u}}
\newcommand{\fI}{\mathfrak{I}}
Let $E$ be a topological stack. In order to write out operations on twisted cohomology effectively
we introduce some notation for operations on $D^+(\Sh_\Ab E)$ or $D(\Sh_\Ab E)$.
If $p:F\to E$ is a map of topological stacks, then we let
$\fp^*:\id\to Rp_*p^*$ denote the unit.
If $p$ is an oriented fiber bundle, then
we let $\fp_!:Rp_*p^*\to \id$ denote the integration map.
If $\pi:E\to B$ is a second map, then we write $\pi_* \fp^*$,
$\pi_*\fp_!$ or simply also $\fp^*$ and $\fp_!$ for the induced transformations
$R\pi_*\pi^*\to R\pi_*Rp_*p^*\pi^*$ and 
$R\pi_*Rp_*p^*\pi^*\to R\pi_*\pi^*$.

If $$\xymatrix{G\ar[rr]^v\ar[d]&&H\ar[d]\\E\ar[rr]^u\ar[dr]^\pi&&F\ar[dl]^{\hat \pi}\\&B&}$$
is a diagram with $U(1)$-gerbes $H\to F$ and $G\to E$ such that the square is Cartesian,
then we write
$P(v)$ for the transformation
$u^*\circ P_H\to P_G\circ u^*$,
and we use the same symbol for the induced transformation
$R \pi_* u^* P_H \hat \pi^*\to R\pi_*P_G u^* \hat \pi^*$.

In a commutative diagram
$$\xymatrix{&F\ar[dl]^p\ar[dr]^{\hat p}&\\
E\ar[dr]^\pi && \hat E\ar[dl]^{\hat \pi} \\
&B&}$$
we will use the symbol $\fI$ or, if necessary,
$\fI_{\pi\circ p=\hat \pi\circ \hat p}$
in order to denote the transformation
$$R\pi_*Rp^*p^*\pi^*\stackrel{\sim}{\to} R\hat \pi_*R\hat p_*\hat p^*\hat \pi^*\ .$$

\subsection{}

We consider a topological gerbe $f\colon G\to E$ with band $U(1)$ over a locally compact stack.
In  \cite{bss}  we define the $G$-twisted cohomology of $E$ with coefficients in $F\in D^+(\Sh_\Ab\bE)$ by $$H^*(E,G;F):=H^*(E;Rf_*f^*(F))\ .$$ 

\subsection{}

Assume now that $f\colon G\to E$ is a part of a $T$-duality diagram
\begin{equation}\label{tzfrzrzr111}
\xymatrix{&p^*G\ar[dl]^q\ar[dr]\ar[rr]^u&&\hat p^* \hat G\ar[dl]\ar[dr]^{\hat q}&\\G\ar[dr]^f&&F\ar[dl]^p\ar[dr]^{\hat p}&&\hat G\ar[dl]^{\hat f}\\&E\ar[dr]^\pi&&\hat E\ar[dl]^{\hat \pi}&\\&&B&&}\ .
\end{equation}
Then we define the transformation
\begin{equation}\label{system41}
J:= \hat \fq_!\circ \fI\circ (\fu^{-1})^*\circ \fq^*:R\pi_*Rf_*f^*\pi^*\to R\hat \pi_*R\hat f_*\hat f^*\hat \pi^*\ .
\end{equation}
Note that here $\fI=\fI_{\pi f q u^{-1}=\hat \pi  f  \hat q }$.

Consider a sheaf $F\in D^+(\Sh_\Ab\bB)$. Note that, by definition,
$H^*(E,G;\pi^*F)=H^*(B;R\pi_*Rf_*f^*\pi^*F)$.
\begin{ddd}
For $F\in D^+(\Sh_\Ab\bE)$ the $T$-duality 
  transformation is defined as the map
 $$T\colon H^*(E,G;\pi^*F)\to H^{*-1}(\hat E,\hat G;\hat \pi^*F)$$
 induced by the natural transformation  (\ref{system41}).
\end{ddd}
 \subsection{}\label{system50}

Let us calculate the effect of the $T$-duality transformation in a simple example.
There is a unique isomorphism class of $T$-duality diagrams over the point $B=*$.
In this case $E=U(1)$ and $G=U(1)\times \cB U(1)$.
We consider a discrete abelian group $F$. Then we have
$$H^*(E,G;\pi^*\uF_{B})\cong \Z[[z]] [v]/(v^2)\otimes F\ ,
 H^*(\hat E,\hat
G;\hat\pi^*\uF_{B})\cong \Z[[z]] [\hat v]/(\hat v^2)\otimes F\ ,$$ where 
$\deg(v)=1=\deg(\hat v)$ and $\deg(z)=2$.

To explicitly calculate the effect of $T$ in this case, observe that the
  cohomology of $Rf_*Rq_*q^*f^* \underline{F}$ is $\Z[[z]]\tensor
  \Lambda(v,\hat v)\tensor F$ with $v$ and $\hat v$ the generators of the two
 $S^1$-factors $E$ and $\hat E$ in $F$. 
The automorphism $u$ induces
  in cohomology, i.e.~on $\Z[[z]]\tensor
  \Lambda(v,\hat v)\tensor F$, the algebra homomorphism given by $z\mapsto z+v\hat v$,
  $v\mapsto v$, $\hat v\mapsto \hat v$. It follows that
  \begin{equation*}
    \begin{split}
      T(z^n\tensor f) &= \int_{F/\hat E}(z^n\tensor f+nz^{n-1}v\hat v\tensor
      f)=nz^{n-1}\hat v\tensor f\\
      T(z^nv\tensor f)&= \int_{F/\hat E} z^nv\tensor f = z^n\tensor f.
  \end{split}
\end{equation*}
We see that the $T$-duality transformation is not an isomorphism.


\subsection{}

Our main motivation for introducing the periodization functor is the
construction of twisted sheaf cohomology
which admits a $T$-duality isomorphism.  Let $G\to E$ be a topological gerbe with band $U(1)$ over a locally compact stack $E$.
\begin{ddd}
We define the periodic $G$-twisted cohomology of $E$ with coefficients in $F\in D^+(\Sh_\Ab\bE)$ by
$$H^*_{per}(E,G;F):=H^*(E;P_G(F))\ .$$
\end{ddd}

Note that here we use the sheaf theory operations for the unbounded derived category, see Subsection \ref{system4001} for details.
\subsection{}

Assume again that $f\colon G\to E$ is part of a $T$-duality diagram (\ref{tzfrzrzr111}).
We define a natural transformation
\begin{equation}\label{system42}
J\colon R\pi_*\circ P_G\circ \pi^*\to R\hat \pi_*\circ  P_{\hat G}\circ \hat \pi^*
\end{equation}
by 
$$J:=\hat \fp_!\circ\fI\circ  P(u)^{-1} \circ \fp^*\ .$$

%
It again involves sheaf theory operations in the unbounded derived category.

Consider a sheaf $F\in D^+(\Sh_\Ab\bB)$.  
Note that by definition $H^*_{per}(E,G;\pi^*F) = H^*(B,R\pi_*P_G(\pi^*(F)))$.
\begin{ddd}\label{system51}
For $F\in D^+(\Sh_\Ab\bE)$ the $T$-duality 
 transformation in periodic twisted cohomology 
 $$T\colon H^*_{per}(E,G;\pi^*F)\to H_{per}^{*-1}(\hat E,\hat G;\hat \pi^*F)$$
is the map induced by the natural transformation  (\ref{system42}).
\end{ddd}

\subsection{}

As an illustration let us calculate the action of the $T$-duality transformation
in the example started in \ref{system50}. 
The sequence $\cS_G(\uF)$ for $F=\Z,\Q,\Q/\Z$ either has trivial $\lim$ or trivial
$\lim^1$. Therefore in this special case the morphism $T$ calculated in
\ref{system50} defines uniquely an
endomorphism of $H^*_{per}(E,G;\pi^*\underline{F}_B)$ (we identify $E\cong \hat E$). For example if  $F=\Q$, then
we read off directly from \ref{system50} that (with $H^{0}_{per}(E,G;\pi^*\underline{\Q}) \cong \Q[v]/v^2$) the $T$-duality morphism is
\begin{equation*}
  T\colon \Q[v]/v^2\to  \Q[v]/v^2\ ,\quad  T(v)=1\ ,\quad  T(1)=v\ .
\end{equation*}
In particular, we see in this example that
now we get an isomorphism. 


\subsection{}

In the remainder of the present subsection we show the following theorem.

\begin{theorem}\label{system66}
The $T$-duality transformation in twisted periodic cohomology \ref{system51}
is an isomorphism.
\end{theorem}
\begin{proof}
The opposite of the $T$-duality diagram (\ref{tzfrzrzr111}) is obtained by
reflecting it in the middle vertical, and by replacing $u$ by its inverse.
We let
$T^\prime\colon H^*_{per}(\hat E,\hat G;\hat \pi^*F)\to H_{per}^{*-1}( E, G;\pi^*F)$
be the associated $T$-duality transformation.

Both, the  $T$-duality diagram and its opposite
can be recognized as subdiagrams of the 
 (slightly extended) big commutative $T$-duality diagram 
\begin{equation}\label{univtdzfe11zuzuqwd}
\xymatrix{&\pr_0^*p^*G\ar[ddl]^{\tilde \pr_0}\ar[ddr]\ar[r]^{\pr_0^*u}& \pr_{\hat E}^*\hat G \ar[dr] \ar[dd]                         \ar[rr]^{\hat c}    &     &\pr^*_{\hat E} \hat G             \ar[dl]             \ar[r]^{\pr_1^*u^{-1}}\ar[dd]&\pr_1^*p^* G \ar[ddl] \ar[ddr]_{\tilde \pr_1}&\\
&\hat p^*\hat G\ar[dd]\ar[rr]&&\hat G\ar[dd]&&\hat p^*\hat G\ar[ll]\ar[dd]\ar[dr]^{u^{-1}}&
\\p^*G\ar[ru]^u\ar[dd]^{f^*p}\ar[dr]^{p^*f}&&F\times_{\hat E}F\ar[dl]^{\pr_0}\ar[dr] ^s                     \ar[rr]^{\hspace{-1cm}c}&& F\times_{\hat E}F              \ar[dr]^{\pr_1}         \ar[dl]            &&p^* G\ar[dd]^{f^*p}\ar[dl]_{p^*f}\\&F\ar[dr]^p\ar[rr]^{\hat p}&&\hat E&&F\ar[ll]^{\hat p}\ar[dl]^{p}&\\G\ar[rr]^f&&E  \ar@{=}[rr]  &&E                                                                                    &&G\ar[ll]^f}
\end{equation}
We now calculate the composition $T^\prime\circ T$.
The compatibility of the integration with pull-back in the Cartesian diagram
$$\xymatrix{F\ar[d]^{\hat p}&F\times_{\hat E}F\ar[d]^{\pr_1}\ar[l]_{\pr_0 }\\\hat E&F\ar[l]^{\hat p}}$$
 is employed in the equality marked by $!$ below. The equality $\hat p\circ \pr_0\circ c^{-1}=\hat p\circ \pr_0$ is used in the equality $!!$. Finally we use $\pr_0\circ c=\pr_1$ at $!!!$.
We have
\begin{eqnarray*}
J^\prime\circ J&=&\fp_!\circ\fI\circ  P(u) \circ \hat \fp^*\circ \hat \fp_!\circ\fI\circ  P(u)^{-1} \circ \fp^*\\
&\stackrel{!}{=}&\fp_!\circ\fI\circ  P(u) \circ \mathfrak{pr_1}_!\circ \fI\circ \mathfrak{pr_0}^* \circ \fI\circ  P(u)^{-1} \circ \fp^*\\&\stackrel{!!}{=}&\fp_!\circ\fI\circ  P(u) \circ \mathfrak{pr_1}_!\circ \fI\circ P(\hat c^{-1})\circ (\mathfrak{ c}^{-1})^*\circ \mathfrak{pr_0}^* \circ \fI\circ  P(u)^{-1} \circ \fp^*\\
&\stackrel{!!!}{=}&\fp_!\circ \mathfrak{pr_1}_!\circ P(\pr_1^* u)\circ  P(\hat c^{-1})\circ P(\pr_0^*u)^{-1}\circ  \mathfrak{pr_1}^* \circ \fp^*\\
&=&\fp_!\circ \mathfrak{pr_1}_!\circ P( \pr_1^* u\circ  \hat c^{-1}\circ(\pr_0^*u)^{-1})\circ  \mathfrak{pr_1}^*\circ \fp^*
\end{eqnarray*}

%
%
This is exactly the transformation associated to the associated small double $T$-duality diagram 
(\ref{sdsdh}) (actually its mirror). Since this is locally isomorphic to the small periodization diagram we see that locally $J^\prime\circ J$ coincides with $\pi_* W$, where $W$ is as in Proposition  \ref{dhewud82d}.  By Proposition \ref{dhewud82d} this transformation is an isomorphism on periodic sheaves
of the form $R\pi_*P_G(\pi^*F)$.
Therefore $T\circ T^\prime$ is an isomorphism.
We can interchange the roles of $T$ and $T^\prime$, hence $T \circ T^\prime$ is an isomorphism, too.
This implies the result.
\end{proof}

\chapter{Orbispaces}\label{system5000}

\section{Twisted periodic delocalized cohomology of orbispaces}

\subsection{}

Let us recall some notions related to orbispaces (compare \cite{math.KT/0609576}). Orbispaces as particular kind of topological stacks 
have previously been introduced in \cite[Sec.~2.1]{MR2246781} and \cite[Sec.~19.3]{math.AG/0503247}). In the present paper  we use the set-up of   \cite{MR2246781} 
but add the additional condition that an orbifold atlas should be separated. This condition is needed in order to show that the loop
stack of an orbifold is again an orbifold.

\begin{enumerate}
\item A topological groupoid $A\colon A^1\Rightarrow A^0$ is called separated if the identity $\eins_A\colon A^0\to A^1$ of the groupoid is a closed map.
\item A topological groupoid $A^1\Rightarrow A^0$ is called proper 
if $(s,r)\colon A^1\to A^0\times A^0$ is a proper map. 
\item A topological groupoid is called {\'e}tale if
the source and range maps $s,r\colon A^1\to A^0$ are  {\'e}tale.
\item A proper \'etale  topological groupoid $A^1\Rightarrow  A^0$ is called very proper if 
there exists a continuous function $\chi \colon  A^0\to [0,1]$ such that
\begin{enumerate}
\item $r\colon \supp(s^*\chi)\to A^0$ is proper
\item  $\sum_{y\in A^x} \chi(s(y))=1$ for all $x\in A^0$.
\end{enumerate}
\item A topological stack is called (very) proper (or {\'e}tale, separated,
  respectively),  if it admits an atlas $A\to X$ such that the topological
  groupoid 
$A\times_XA\Rightarrow A$ is (very) proper  (or {\'e}tale, separated, respectively).
 
\item An orbispace atlas of a topological stack $X$ is an atlas $A\to X$ such that $A\times_XA\Rightarrow A$ is a very proper  {\'e}tale and separated  groupoid. 
\item
An orbispace $X$ is a topological stack which admits an orbispace atlas.
\item If $X,Y$ are orbispaces, then a morphism of orbispaces $X\to Y$
is a representable morphism of stacks.
\item A locally compact orbispace is an orbispace $X$ which admits an orbispace atlas $A\to X$ such that $A$ is locally compact. 
 \end{enumerate}

\subsection{}

If $X$ is a stack, then its inertia stack (sometimes called loop stack) $LX$
is defined as the two-categorical equalizer 
of the diagram $$\xymatrix{X\ar@{=>}[r]^{\id_X}_{\id_X}&X}\ .$$ 
In \cite[Sec 2.2]{math.KT/0609576} we have introduced   an
explicit model of $LX$ and studied its properties. {
The loop stack $LX$ depends $2$-functorially on $X$. Indeed,  since $\underline{\Hom}_{\tt
  Cat}$ is a strict $2$-functor, the loop functor is
a strict functor between $2$-categories.} {As already
mentioned before, later we will suppress the 2-morphisms in 2-commutative diagrams
in 2-categories for better legibility.} If
$X$ is a topological stack (orbispace), then 
$LX$ is a topological stack (orbispace), too (see \cite[Lemma
2.25]{math.KT/0609576}, \cite[Lemma 2.33]{math.KT/0609576}).

\begin{lem}
If $X$ is a locally compact orbispace, then $LX$ is a locally compact orbispace, too.
\end{lem}
\begin{proof}
Let $A\to X$ be a locally compact orbispace atlas of $X$.
Then we have the proper, separated and \'etale topological  groupoid
$A\times_XA\Rightarrow A$. Since the source map of this groupoid is  \'etale, the space of
morphisms $A\times_XA$ of this groupoid is locally compact, too.

In the proof of Lemma \cite[Lemma 2.25]{math.KT/0609576} we constructed an orbispace atlas 
$W\to LX$ of $LX$, where $W$ was given by the pull-back of spaces
$$\xymatrix{W\ar[r]\ar[d]^w&A\times_XA\ar[d]^{(\pr_1,\pr_2)}\\
A\ar[r]^{\diag}&A\times A}\ .$$
This implies that $W$ is locally compact.
\end{proof}

\subsection{}\label{system60}

Let $G\to X$ be a topological gerbe with band $U(1)$ over a locally compact  orbispace. 
The truly interesting $G$-twisted cohomology of $X$ (with complex coefficients)
is not the cohomology $H^*_{per}(X,G;\uC)$ (see \ref{system51}), 
but a more complicated delocalized version $H^*_{deloc,per}(X,G)$,
which we will define below (see \cite[Sec.~1.3]{math.KT/0609576} for an explanation).

As shown in \cite[Sec.~2.5]{math.KT/0609576} the gerbe gives rise to a principal bundle $\tilde G^\delta\to LX$
with structure group $U(1)^\delta$ in a functorial way,
 where $U(1)^\delta$ denotes the group $U(1)$ with
the discrete topology. By $\cL\in \Sh_\Ab\bLX$
we denote the sheaf of locally constant sections of the associated vector bundle $\tilde G^\delta\times_{U(1)^\delta}\C\to LX$.

We define the gerbe $G_L\to LX$ as the pull-back
$$\xymatrix{G_L\ar[d]^{f_L}\ar[r]&G\ar[d]^f\\LX\ar[r]&X}\ .$$

\begin{ddd}\label{eiuwqd}
We define 
$$\cL_G:=P_{G_L}(\cL)\in D(\Sh_\Ab \bLX)\ .$$
The $G$-twisted delocalized periodic cohomology of $X$ is defined as
$$H^*_{deloc,per}(X,G):=H^*(LX;\cL_G)\ .$$
\end{ddd}

\section[$T$-duality in twisted periodic delocalized cohomology]{The $T$-duality transformation in twisted periodic delocalized cohomology}

\subsection{}

 We consider a $T$-duality diagram 
\begin{equation}\label{tzfrzrddzr123}
\xymatrix{&p^*G\ar[dl]\ar[dr]\ar[rr]^u&&\hat p^* \hat G\ar[dl]\ar[dr]&\\G\ar[dr]^f&&F\ar[dl]^p\ar[dr]^{\hat p}&&\hat G\ar[dl]^{\hat f}\\&E\ar[dr]^\pi&&\hat E\ar[dl]^{\hat \pi}&\\&&B&&}
\end{equation}
(see Definition  \ref{eruihfrvc}), where $B$ is a  locally compact  orbispace.

We apply the loops functor $L\colon orbispaces\to orbispaces$ to the subdiagram
$$\xymatrix{&F\ar[dr]^{\hat p}\ar[dl]_{p}&\\E\ar[dr]^\pi&&\hat E\ar[dl]_{\hat \pi}\\&B&}$$ and get
$$\xymatrix{&LF\ar[dr]^{L\hat p}\ar[dl]_{Lp}&\\LE\ar[dr]^{L\pi}&&L\hat E\ar[dl]_{L\hat \pi}\\&LB&}\ .$$ 
In the first diagram the maps $p,\hat p,\pi,\hat \pi$ are all $U(1)$-principal bundles.
The maps $Lp,L\hat p,L\pi,L\hat \pi$ are not necessarily surjective. Thus in general the derived diagram of loop stacks is not part of a $T$-duality diagram. But it is so locally in a certain sense which we will explain in the following.

\subsection{}

We can extend the second diagram by the local systems (see \ref{system60})
\begin{equation}\label{tzfrzrddzr444}
\xymatrix{&Lp^*\cL\ar[dl]\ar[dr]\ar[rr]^u&&L\hat p^* \hat \cL\ar[dl]\ar[dr]&\\\cL\ar[dr]&&LF\ar[dl]^{Lp}\ar[dr]^{L\hat p}&&\hat \cL\ar[dl]\\&LE\ar[dr]^{L\pi}&&L\hat E\ar[dl]^{L\hat \pi}&\\&&LB&&}
\end{equation}
and the pull-backs of gerbes
\begin{equation}\label{tzfrzrddzr}
\xymatrix{&Lp^*G_L\ar[dl]\ar[dr]\ar[rr]^u&&L\hat p^* \hat G_L\ar[dl]\ar[dr]&\\G_L\ar[dr]^{f_L}&&LF\ar[dl]^{Lp}\ar[dr]^{L\hat p}&&\hat G_L\ar[dl]^{\hat f_L}\\&LE\ar[dr]^{L\pi}&&L\hat E\ar[dl]^{L\hat \pi}&\\&&LB&&}
\end{equation}
In particular, we have an isomorphism
\begin{equation}\label{ezfgbdmnwdas}
u:Lp^*\cL_G\stackrel{\sim}{\to} L\hat p^* \hat\cL_{\hat G}\ .
\end{equation}

\subsection{}\label{system62}

Note that $\hat p\colon F\to \hat E$ is a $U(1)$-principal bundle.
In \cite[Lemma 2.34]{math.KT/0609576} we have constructed
a map $h\colon L\hat E\to U(1)^\delta$ which measures the action of the automorphisms of the points of $\hat E$
on the fibers of $\hat p$. We get a decomposition into a disjoint union of open substacks
$$L\hat E\cong \bigsqcup_{u\in U(1)} L\hat E_u\ ,$$
where $L\hat E_u:=h^{-1}(u)$. Here and in the following we use the simplified
notation $h^{-1}(u)$ for the pullback of $h\colon L\hat E\to U(1)^{\delta}$
along the inclusion $i_u\colon *\to U(1)$ with $ i_u(*):=u$.  By \cite[Lemma 2.36]{math.KT/0609576}, the map
$L\hat p\colon LF\to L\hat E$ factors over the inclusion $J:L\hat E_1\to L\hat
E$, and the corresponding map $L\hat p_1\colon LF\to L\hat E_1$ is a
$U(1)$-principal bundle. The integration 
$$\mathfrak{L\hat p_1}_!\colon R(L\hat p_1)_*\circ L\hat p_1^*\to \id$$
is well-defined. The open inclusion $J$ induces
a natural transformation $\mathfrak{J}_!\colon RJ_*\circ J^*\to \id$.
We can thus define
$$\mathfrak{L\hat p}_!:=\mathfrak{J}_! \circ \mathfrak{L\hat p_1}_! \colon  RL\hat p_*\circ L\hat p^*  \to \id\ .$$

\subsection{}

\begin{ddd}\label{system61}
The local $T$-duality transformation   associated to the diagram (\ref{tzfrzrddzr123}) is 
given by the composition
\begin{eqnarray*}
T_{loc}:=\mathfrak{L\hat p}_!\circ u\circ \mathfrak{L p}^*\colon RL\pi_*\cL_G\to RL\hat \pi_*\hat \cL_{\hat G}
\ ,\end{eqnarray*}
where $u$ is induced by (\ref{ezfgbdmnwdas}).
\end{ddd}

Note that $H^*_{deloc,per}(E,G)\cong H^*(LB;RL\pi_*\cL_G)$.
Hence we can make the following definition.
\begin{ddd}
The $T$-duality transformation in twisted periodic delocalized cohomology associated to the
$T$-duality diagram (\ref{tzfrzrddzr123}) is the transformation
$$T\colon H^*_{deloc,per}(E,G)\to H^*_{deloc,per}(\hat E,\hat G)$$ induced by the local
$T$-duality transformation $T_{loc}$ defined in \ref{system61}.
\end{ddd}

\section{The geometry of $T$-duality diagrams over orbispaces}

%
%
%
%

\subsection{}
We consider a $T$-duality diagram (\ref{tzfrzrddzr123}) over a locally compact orbispace. 
As explained in \cite[Sec.~2.5]{math.KT/0609576} (see also \ref{system60})
 the gerbe $G\to E$ naturally gives rise to a $U(1)^\delta$-principal bundle $\tilde G^\delta\to LE$.
Let  $g\colon LB_1\to U(1)^\delta$ be the function which describes the holonomy of the bundle $\tilde G^\delta\to LE$ along the fibers of $LE\to LB_{1}$ (see \cite[2.6.3]{math.KT/0609576}).
In the following we recall from \cite{math.KT/0609576} a cohomological description of the functions
$g$ and $h$ (introduced in \ref{system62}).
 
Let $c_1\in H^2(B;\Z)$ denote the first Chern class of the $U(1)$-principal bundle $\pi\colon E\to B$, and let
 $d\in H^3(E;\Z)$ denote the Dixmier-Douady class of the gerbe $f\colon G\to E$.
By integration over the fiber it gives rise to a class $\int_\pi d\in H^2(B;\Z)$.
 In \cite[2.4.11]{math.KT/0609576} we have shown that a class $\chi\in H^2(B;\Z)$ gives rise to a function $\bar\chi\colon LB\to U(1)^\delta$ in a natural way.  

\begin{prop}[Lemma 2.38 and Prop. 2.49  \cite{math.KT/0609576} ]\label{auqewmdq}
We have the equalities
\begin{enumerate}
\item
$$\overline{c_1}=h\colon LB\to U(1)^\delta\ .$$
\item
$$\overline{\int_\pi d}_{|LB_1}=g\colon  LB_1\to U(1)^\delta\ .$$ 
\end{enumerate}
\end{prop}

\subsection{}

We now have functions $h,\hat h\colon LB\to U(1)^\delta$ associated to the $U(1)$-principal bundles
$\pi\colon E\to B$ and $ \hat \pi\colon \hat E\to B$.  We define
$$LB_{(1,*)}:=h^{-1}(1)\ ,\quad LB_{(*,1)}:=\hat h^{-1}(1)\ .$$ 
We furthermore have functions (see \ref{tzfrzrddzr123})
$$g\colon LB_{(1,*)}\to U(1)^\delta\ ,\quad \hat g\colon LB_{(*,1)}\to U(1)^\delta$$
measuring the holonomy of $\tilde G^\delta\to LE$ and $\tilde {\hat G}^\delta\to L\hat E$ along the fibers.

\begin{prop}\label{ueidkjqwdq}
We have the equalities
$$\hat g=h^{-1}_{|LB_{(*,1)}}\ ,\quad g=\hat h^{-1}_{|LB_{(1,*)}}\ .$$
\end{prop}
\begin{proof}
Let $$d\in H^3(E;\Z)\ , \quad \hat d\in H^3(\hat E;\Z)$$
be the Dixmier-Douady classes of the gerbes $G_L\to E$ and $\hat G_L\to \hat E$.
Furthermore let  $$c_1,\hat c_1\in H^2(B;\Z)$$ denote the first Chern classes of the
 $U(1)$-principal bundles $\pi\colon E\to B$ and $\hat \pi\colon \hat E \to B$.
 The theory of $T$-duality for orbispaces \cite{MR2246781} gives the equalities
$$c_1=-\hat \pi_!(\hat d)\ ,\quad \hat c_1 = - \pi_!(d)\ .$$
Hence the assertion follows from Proposition \ref{auqewmdq}. 
\end{proof}

\section[$T$-duality is an isomorphism]{The $T$-duality transformation in  twisted periodic delocalized cohomology is an isomorphism}

\subsection{}
 Let us consider a $U(1)$-principal bundle $\pi\colon E\to B$ in locally compact orbispaces with first Chern class $c_1\in H^2(B;\Z)$ and a topological $U(1)$-banded gerbe $f\colon G\to E$ with Dixmier-Douady class $d\in H^3(E;\Z)$. In Definition \ref{eiuwqd} we have introduced the object $\cL_G\in D(\Sh_\Ab\bLE)$. Furthermore we have $U(1)^\delta$-valued functions
$h=\overline{c_1}$ and $g=\overline{\pi_!(d)}$ on $LB$. Let $LB_1:=h^{-1}(1)$
 and note that $L\pi\colon LE\to LB$ factors over the $U(1)$-principal bundle
$L\pi\colon LE\to LB_1$. We fix $u\in U(1)^\delta\setminus \{1\}$ and consider the component
$LB_{(1,u)}:=h^{-1}(1)\cap g^{-1}(u)$.
\begin{lem}\label{hjwegbd}
We have $R\pi_*(\cL_G)_{|LB_{(1,u)}}\cong 0$.
\end{lem}
\begin{proof}
Let $(T\to LB_{(1,u)})\in \bLB_{(1,u)}$. 
After refining $T$  by a covering we can assume that there is a diagram 
$$\xymatrix{\cB U(1)\ar[d]^y&U(1)\times \cB
  U(1)\ar[l]^z\ar[d]^x&s^*G_L\ar[d]\ar[l]\ar[r]&{G_L}_{(1,u)}\ar[d]\\{*}&U(1)\ar[l]^q\ar[d]^q&T\times
  U(1)\ar[l]^v\ar[r]^s\ar[d]^p&LE_{(1,u)}\ar[d]^\pi\\&{*}&T\ar[l]^w\ar[r]^t
  &LB_{(1,u)}}$$  
of Cartesian squares.
We get
\begin{eqnarray*}
t^*R\pi_*(\cL_G)&\cong&Rp_* s^*(\cL_G)\\
&= &Rp_* s^*(P_{G_L}(\cL))\\
&\cong&Rp_*P_{s^*G_L}(s^*\cL)\ .
\end{eqnarray*}
Let $\cH\in \Sh_\Ab(\Site(U(1)))$ be the locally constant sheaf over $U(1)$
with fiber $\C$ and holonomy $u\in U(1)\setminus \{1\}$. Then we have $s^*\cL\cong v^* \cH$.  
We calculate further
\begin{eqnarray*}
Rp_*P_{s^*G_L}(s^*\cL)&\cong&Rp_* P_{s^*G_L}(v^*\cH)\\
&\cong&Rp_* v^*P_{U(1)\times \cB U(1)}(\cH)\\&\cong&
w^*Rq_*P_{U(1)\times \cB U(1)}(\cH)\ .
\end{eqnarray*}
It remains to show that
$$Rq_*P_{U(1)\times \cB U(1) }(\cH)\cong 0\ .$$ Recall from \ref{system19} that the object
$P_{U(1)\times \cB U(1)}(\cH)\in D(\Sh_\Ab \Site(U(1)))$ is given (up to non-canonical isomorphism) by the $\holim$ of a diagram
$$0\leftarrow Rx_* x^*(\cH)\stackrel{D}{\leftarrow} Rx_* x^*(\cH) [2]\stackrel{D}{\leftarrow} Rx_* x^*(\cH)[4]\stackrel{D}{\leftarrow} Rx_* x^*(\cH)[6]\dots\ .$$
The functor $Rq_*$ commutes with this $\holim$\footnote{$Rq_*$ is a right-adjoint and commutes with products and mapping cones}. Therefore
$Rq_*P_{U(1)\times \cB U(1)}(\cH)$ is given by the $\holim$ of the diagram
\begin{multline*}
  0\leftarrow Rq_* Rx_* x^*(\cH)\stackrel{Rq_*(D)}{\leftarrow} Rq_* Rx_*
  x^*(\cH)[2]\\
 \stackrel{Rq_*(D)}{\leftarrow} Rq_*Rx_*
  x^*(\cH)[4]\stackrel{Rq_*(D)}{\leftarrow} Rq_*Rx_* x^*(\cH)[6]\dots\ .
\end{multline*}
The following calculation uses the projection formula twice, first by Lemma \ref{gr213} for the non-representable map 
$x$ and a tensor product with a one-dimensional local system of complex vector
spaces $\cH$, secondly using Lemma 
\ref{projefoa} for the proper representable map $q$ and the tensor product
with the bounded below object $Ry_*(i^\sharp\Z_{\Site([*/U(1)])})\in
D^+(\Sh_\Ab\Site(U(1)))$ 
\begin{eqnarray*}
Rq_*Rx_*x^*(\cH)&\cong&Rq_*Rx_*(\uZ_{\Site(U(1)\times \cB U(1)}\otimes x^*(\cH))\\
&\cong&Rq_* (Rx_*(\uZ_{\Site(U(1)\times \cB U(1))})\otimes \cH)\\
&\cong&Rq_*(Rx_*(z^*\uZ_{\Site(\cB U(1))})\otimes \cH)\\
&\cong&Rq_*(q^*(Ry_*\uZ_{\Site(\cB U(1))})\otimes \cH)\\
&\cong&Ry_*\uZ_{\Site(\cB U(1))}\otimes Rq_*(\cH)\ .
\end{eqnarray*}
Since the holonomy of $\cH$ along $U(1)$ is non-trivial, and the cohomology of
$S^1$ with coefficients in a non-trivial flat line bundle is trivial,  we have
$$Rq_*(\cH)\cong 0\ .$$
\end{proof}

\subsection{}

We now consider a $T$-duality diagram  (\ref{tzfrzrddzr123})
 where $B$ is a locally compact  orbispace.

\begin{theorem}\label{system1006}
The local  $T$-duality transformation (Definition \ref{system61}) 
$$T_{loc}\colon RL\pi_*(\cL_G)\to RL\hat \pi_*(\hat \cL_{\hat G})[-2]$$ is an isomorphism in $D(\Sh_\Ab\bLB)$.
In particular, the $T$-duality transformation
$$T\colon H^*_{deloc,per}(E,G)\to H^*_{deloc,per}(\hat E,\hat G)$$
is an isomorphism.
 \end{theorem}
\begin{proof}
We have functions
$h,\hat h\colon LB\to U(1)$ which define substacks
$LB_{(1,*)}:=h^{-1}(1)$ and $LB_{(*,1)}:=\hat h^{-1}(1)$.
By Proposition \ref{ueidkjqwdq} we have $g=\hat h^{-1}_{|LB_{(1,*)}}:LB_{(1,*)}\to U(1)^\delta$.
By Lemma \ref{hjwegbd} the object
$RL\pi_*(\cL_G)\in D(\Sh_\Ab \bLB)$ is supported on 
$$g^{-1}(1)=LB_{(1,*)}\cap
LB_{(*,1)}=: LB_{(1,1)}\ .$$
Note that
$\hat g=h^{-1}_{|LB_{(*,1)}}$, so that $RL\hat \pi_*\hat \cL_{\hat G}$ is supported on $LB_{(1,1)}$, too. Let $i\colon LB_{(1,1)}\to LB$ denote the inclusion. The following diagram is the pull-back of  (\ref{tzfrzrddzr123}) via the map $LB_{(1,1)}\to LB\to B$  
\begin{equation}\label{wedddbwed}
{\small 
\xymatrix{&p_L^*(G_L)_{|LE_{|LB_{(1,1)}}}\ar[dl]\ar[dr]\ar[rr]^{u_L}&&\hat
  p_L^* (\hat G_L)_{|L\hat
    E_{|LB_{(1,1)}}}\ar[dl]\ar[dr]&\\(G_L)_{|LE_{|LB_{(1,1)}}}\ar[dr]^{f_L}&&LF_{|LB_{(1,1)}}\ar[dl]^{Lp}\ar[dr]^{\hat
    Lp}&&(\hat G_L)_{|L\hat E_{|LB_{(1,1)}}}\ar[dl]^{\hat
    f_L}\\&LE_{|LB_{(1,1)}}\ar[dr]^{L\pi_1}&&L \hat
  E_{|LB_{(1,1)}}\ar[dl]^{L\hat \pi_1}&\\&&LB_{(1,1)}&&} 
}
\end{equation}
We consider $$\cL_1:=\cL_{|LE_{|LB_{(1,1)}}}\ ,\quad \hat \cL_1:=\hat \cL_{|L\hat E_{|LB_{(1,1)}}}\ .$$
Because we restrict to the subset $LB_{(1,1)}$ of trivial holonomy we have isomorphisms
 $$\cL_1\cong L\pi_1^*\uC_{\bLB_{(1,1)}}\,\quad \hat \cL_1\cong L\hat  \pi_1^* \uC_{\bLB_{(1,1)}}\ .$$
The local $T$-duality transformation $T_{loc}$ is now locally equal to
the transformation $J$ defined in \ref{system42} applied to the $T$-duality diagram
(\ref{wedddbwed}) and the sheaf $\uC_{\bLB_{(1,1)}}$. As in the proof of Theorem \ref{system66} one shows, using the commutative double $T$-duality diagram, that $T_{loc}$ is an isomorphism.

The global second assertion can be deduced directly from Theorem \ref{system66}. 
By
the observation on the support of $RL\pi_*(\cL_G)\in D(\Sh_\Ab \bLB)$ made above
we get
$$H^*_{deloc,per}(E,G)\cong H^*_{per}(LB_{(1,1)};RL( \pi_1)_*P_{( G_L)_{|L E_{|LB_{(1,1)}}}}(L\pi_1^*\uC_{\bLB_{(1,1)}}))\ ,$$
and similarly
$$H^*_{deloc,per}(\hat E,\hat G)\cong H^*_{per}(LB_{(1,1)};RL(\hat \pi_1)_*P_{(\hat G_L)_{|L\hat E_{|LB_{(1,1)}}}}(L\pi_1^*\uC_{\bLB_{(1,1)}}))\ .$$
With these identifications the  $T$-duality transformation in twisted periodic delocalized cohomology is then equal to the $T$-duality transformation in twisted periodic cohomology for the diagram (\ref{wedddbwed}) and the sheaf $\uC_{\bLB_{(1,1)}}\in D^+(\Sh_\Ab \bLB_{1,1})$.
 \end{proof}

\chapter{Verdier duality for locally compact stacks}\label{system3003}

\section{Elements of the theory of stacks on $\Top$ and sheaf theory}\label{wuiefwefwefqwd}

\subsection{}

In the present paper we consider stacks on the site $\Top$. 
A prestack is a lax presheaf $X$ of groupoids on $\Top$. The prefix "lax" indicates that for a pair of composable morphisms $u\colon U\to V$, $v\colon V\to W$ we have a natural transformation of functors $\phi_{u,v}\colon X(u)\circ X(v)\to X(v\circ u)$ which is not necessarily the identity, and  which satisfies a compatibility condition for triples.
A prestack is a stack  if it  satisfies the standard descent conditions on the level of objects and morphisms.
 A sheaf of sets can be considered as a stack in the canonical way. Via the Yoneda embedding
 $\Top\to \Sh\Top$ (note that the topology of $\Top$ is sub-canonical, i.e. representable presheaves are sheaves) we consider topological spaces as stacks in the natural way.

\subsection{}\label{staqwndwqodwqd}

In the following we collect some definitions and facts of the theory of stacks in topological spaces.
Stacks are objects of a two-category, and fibre products and more general limits in stacks are understood in the two-categorial sense. Note that two-categorial limits in stacks exists (see \cite{math.KT/0609576} for more information), and that the inclusion of spaces into stacks preserves those limits.
A useful reference for stacks in topological spaces and manifolds is the
survey \cite{heinloth}.

\begin{enumerate}
\item
A morphism of stacks $G\rightarrow H$ is called representable, if
for each space $U$ and map $U\rightarrow H$ the fibre product $U\times_HG$ is equivalent to a space.
\item A representable map $G\to H$ between stacks is called proper if for every map $K\to H$ from a compact space the fibre product $K\times_HG$ is a compact space.
\item
A map $f:A\to B$ of topological spaces has local sections if for each point $b\in B$ in the image of $f$ there exists a neighbourhood $b\in U\subseteq B$ and a map $s:U\to A$ such that $f\circ u=\id_{U}$.
\item
A representable morphism $G\rightarrow H$  has local sections  if for every map $U\to H$ from a space the induced map $U\times_HG\to U$ of spaces has local sections.
\item A representable map $G\to H$ is surjective if for every map $U\to H$ from a space
the induced map $U\times_HG\to U$ is a surjective map of spaces.
\item
A map $A\to X$ from a space $A$ to a stack $X$ is called an atlas of $X$, if
it is surjective, representable
and  admits local sections. A stack which admits an atlas is called a topological stack.
\item A morphism (not necessarily representable) between topological stacks $G\to H$ is surjective  (or has local sections, respectively) if for an atlas $A\to G$ the composition
$A\to G\to H$ is surjective (or has local sections, respectively) (note that this composition is representable by Proposition \ref{lem:representability} below).  
\item
 A composition of maps with local sections has local sections.
The corresponding  assertion is true for the following properties of maps:
\begin{enumerate}
\item representable
\item representable and proper
\item surjective.
\end{enumerate}
\item
Consider a two-cartesian diagram of stacks
 \begin{equation*}
    \begin{CD}
      H @>>v> G\\
      @VVgV @VVfV\\
      Y @>>u> X
    \end{CD}
  \end{equation*}
If $u$ has local sections, then so has $v$. If $f$ is representable, then so is $g$.
  \end{enumerate}

\subsection{}

The inclusion of spaces into sheaves and of sheaves into stacks preserves small limits, where limits in stacks are understood in the two-categorical sense. This implies that a map of spaces $X\to Y$ is representable.
In fact we have the following more general result.

\begin{prop}\label{lem:representability}
  Let $G$ be a topological stack and $X$ a space. Then every morphism $f\colon
  X\to G$ is representable.
\end{prop}
The proof will be given in \ref{jhsgbwqjswqs}
and needs some preparations.

\subsection{}

We will need the notion of an open substack.
  \begin{definition}\label{def:open_emb_of_stacks}
    Let $G$ be a stack in topological spaces. A morphism $H\to G$ of stacks is an
    embedding of an open substack, if it is representable and for each  map $T\to
    G$ from a space $T$ the induced map of spaces $T\times_GH\to T$ is an open embedding
    of topological spaces. 
  \end{definition}
Note that, via Yoneda, an open embedding of spaces is an open embedding of stacks.

  \begin{definition}\label{didwedqwdqwdqwdd}
   A morphism $U\to G$ of topological stacks is locally an open embedding if $U\cong \bigsqcup_{i\in I} U_i$ for a collection $(U_i)_{i\in I}$ of topological stacks and $U_i\to G$ is an
    embedding of an open substack for every $i\in I$.
\end{definition}
\newcommand{\coeq}{{\mathrm{coeq}}}
Let us first characterize spaces as stacks which can be covered by a collection of spaces.
\begin{lem}\label{kjhdkqwdqwwqd}
Let $X$ be a stack in topological spaces for which there exists a morphism $U\to X$ from a space which is surjective and locally an open embedding. Then $X$ is equivalent to a space.
 \end{lem}
\begin{proof}
Let $U\cong \sqcup_{i} U_i$ be such that $U_i\to X$ is an open embedding for all $i$.
Then we define the space $B$ as the coequalizer in spaces
\begin{equation}\label{asdwdisdwqdqdqw}
B:=\coeq(\bigsqcup_{i,j} U_i\times_X U_j\rightrightarrows \bigsqcup_i U_i)\ .
\end{equation}
Since $U_i\to X$ is an open embedding we see that $\pr_{U_i}\colon U_i\times_X U_j\to U_i$
is an open embedding. We can now refer to  \cite[Prop. 16.1]{math.AG/0503247} and deduce that
the equalizer in spaces $B$ is also the two-categorical equalizer in stacks of the diagram (\ref{asdwdisdwqdqdqw}), which is of course equivalent to $X$. Note that the difficulty at this point is that the embedding of the category of spaces (viewed as a two-category)  into the two-category of stacks does not preserve general small colimits, as opposed to the case of limits.

For completeness we will give an argument.
First note that $\pr_{U_i}\colon U_i\times_X U_i\stackrel{\sim}{\to} U_i$ is a homeomorphism. It thus follows from the groupoid structure of the coequalizer diagram that
$U_i\to B$ is injective for all $i$.  Since $\bigsqcup_i U_i\to B$
is a topological quotient map it is open. Therefore  $\bigsqcup_iU_i\to B$ is a open covering.
We further conclude that the natural map
$U_i\times_X U_j\to U_i\times_BU_j$ is in fact a homeomorphism.

The claim is that $X$ is equivalent to $B$. We first construct a morphism $X\to B$.
Let $(T\to X)\in X(T)$. Then $(T_i:=T\times_XU_i)_i$ is an open covering of $T$. 
Using the identification $T_i\times_TT_j\cong T\times_X(U_i\times_X U_j)$ we get a diagram
$$\xymatrix{\bigsqcup_{i,j}T_i\times_TT_j\ar@/_-0.2cm/[d]\ar@/^-0.2cm/[d]\ar[r]&U_i\times_XU_j\ar@/_-0.2cm/[d]\ar@/^-0.2cm/[d]\\\bigsqcup_i T_i\ar[r]\ar[d]&\bigsqcup_i U_i\ar[d]\\T\ar@{.>}[r]&B}\ ,$$
where the horizontal maps are induced by the projections $T\times_XU_i\to U_i$, and
the left vertical is the representation of $T$ as a coequalizer.
Therefore we obtain a unique factorization $(T\to B)\in B(T)$.
The construction is functorial in $T$ and therefore induces a morphism $X\to B$. 

In order to see that it has an inverse let $(T\to B)\in B(T)$ be given. Then we define the open covering
$(T_i:=T\times_BU_i)_i$ of $T$.
The compositions $$\phi_i:T_i\cong T\times_B U_i\stackrel{\pr_{U_i}}{\to} U_i\to X$$ can be considered as a collection of objects $(\phi_i\in X(T_i))_i$.
The induced map
\begin{multline*}
  T_i\cap T_j\cong T_i\times_TT_j\cong (T\times_BU_i)\times_T
  (T\times_BU_j)\cong T\times_B(U_i\times_B U_j)\\
\stackrel{\pr_{U_i\times_B
      U_j}}{\to} U_i\times_B U_j\cong U_i\times_X U_j\to X\times_XX
\end{multline*}
can be considered as a collection of isomorphisms $\phi_{ij}:(\phi_i)_{|T_i\cap T_j}\stackrel{\sim}{\to} (\phi_j)_{|T_i\cap T_j}$
 which satisfy the cocycle condition on triple intersections. Since $X$ is a stack we
 can therefore glue the local maps and get a map $(T\to X)\in X(T)$ which is unique
up to unique isomorphism. This construction is again functorial in $T$ and provides 
the map $B\to X$.

It is easy to see that both maps $X\to B$ and $B\to X$ constructed above are mutually inverse.
\end{proof}

\subsection{}\label{jhsgbwqjswqs}

%
 
We now show Proposition \ref{lem:representability}
\begin{proof}   
Consider a map $T\to G$ from a space $T$.
  We have to prove that the fiber product $T\times_G X$ is
  equivalent to a space. Using the assumption that $G$ is topological we choose 
  an atlas $A\to G$ of $G$.
  Because $A\to G$ has local sections, we can find an open covering $\bigsqcup_{i\in I}U_i=: U\to X$
  such that $U\times_G A\to U$ 
  has a section $s\colon U\to U\times_G A$.  We first want to show that $T\times_G U$ is a space.
  Since the structure map $A\to G$ of an atlas is representable we know that $U\times_GA$ and $T\times_GA$ are spaces.
Therefore,
  $T\times_GU\times_G A\cong  (T\times_G A)\times_A (U\times_G A)$ is a space, too. The section
  $s$ pulls back to a section $\hat s:T\times_GU\to
  T\times_G U\times_GA$ which implements $T\times_G U$ as a subspace of
  the space $T\times_G U\times_G A$.
$$\xymatrix{&T\times_GU\times_GA\ar[rr]\ar[dl]&&U\times_GA\ar[dl]\ar[ddd]\\
T\times_GU\ar@/_-1cm/[ur]^{\hat s}\ar[rr]\ar[d]&&U\ar[d]\ar@/_-1cm/[ur]^s&\\
T\times_GX\ar[rr]\ar[dd]&&X\ar[dd]&\\
&&&A\ar[dl]\\
T\ar[rr]&&G&}\ .$$
Since the map $U\to X$ is surjective and locally an open embedding its pull-back
$T\times_GU\to T\times_GX$ is surjective and locally an open embedding, too.
Therefore by  Lemma \ref{kjhdkqwdqwwqd} the stack $T\times_GX$
is equivalent to a space. 
\end{proof}

\subsection{}\label{hshsqiwhswqs}

Recall that
a topological stack is called locally compact if it admits a locally compact atlas $A\to G$ such that $A\times_GA$ is a locally compact space.
Furthermore recall that the site $\bX=\Site(X)$ associated to a locally compact stack $X$ is the full subcategory of locally compact spaces $U\to X$ over $X$ such that the structure map has local sections. A morphism in this site $\bX$ is a diagram
\begin{equation}\label{gsdhasddas}
\xymatrix{U\ar[dr]\ar[rr]&\ar@{:>}[d]&V\ar[dl]\\&X&}
\end{equation}
consisting of a morphism of spaces over $X$ and a two-morphism.
The topology on $\bX$ is given by the covering families of the objects $(U\to X)$ induced by open covering of $U$.

Much of the general theory would work without the assumption of local compactness.
But local compactness is important in connection with the projection formula Lemma \ref{system81}
which is a crucial ingredient of the theory of integration. Since the latter is our main goal
of the present section we generally adopt the restriction to locally compact stacks.

\subsection{}

The sheaf theory for topological stacks can be built in a parallel manner to the sheaf theory for smooth stacks developed in \cite{bss}.
The transition goes via the following replacements of words:
\begin{enumerate}
\item For the definition of stacks the site of smooth manifolds $\Mf^\infty$ is replaced by the site of topological spaces $\Top$. In the definition of the site of a locally compact stack manifolds are replaced by locally compact spaces.
\item The concept of a \textit{smooth stack} is replaced by the concept of a  \textit{locally compact  stack}.
\item The notion of a  \textit{smooth map} is replaced by the notion of a \textit{map which admits local sections}.
\end{enumerate}

In the present paper we freely use results in the general sheaf theory for topological stacks 
from \cite[Sec.~2]{bss} in the case of stacks in topological spaces
which are proved there for manifolds. It should be noted that with the
conventions just made, all statements and proofs carry over verbatim

\subsection{}\label{echejwwcc}\label{lechejwwcc}

Let $X$ be a locally compact stack. By
 $\Pr\bX$ and $\Sh\bX$ we denote the categories of presheaves and sheaves on $\bX$.
They are related by 
 a pair of adjoint functors
  $$i^\sharp\colon \Pr\bX\leftrightarrows\Sh\bX:i\ .$$
The sheafification functor $i^\sharp$ is exact.

\subsection{}\label{lem:lrexact} \label{sec:morph_of_sheaves}

Let $f:X\to Y$ be a morphism of locally compact stacks.
In induces a functor
${}^pf_*:\Pr\bX\to \Pr\bY$ by
$${}^pf_*F(V\to Y):=\lim F(U\to X)\ ,$$ where
the limit is taken over the category of diagrams
\begin{equation}\label{dergetzwe}
\xymatrix{U\ar[d]\ar[r]&X\ar[d]^f\\V\ar[r]\ar@{:>}[ur]&Y}
\end{equation}
with $(U\to X)\in \bX$. For details we refer to  \cite[Sections 2.1, 2.2]{bss}.
This functor fits into an adjoint pair
$${}^pf^*:\Pr\bY\leftrightarrows\Pr\bX:{}^pf_*\ .$$
The functor 
${}^pf^*$ is given by
$${}^pf^*G(U\to X)=\colim\: G(V\to Y)\ ,$$
where 
the colimit is again taken over the category of diagrams
with $(V\to Y)\in \bY$.


We extend these functors to sheaves by
$$f_*:=i^\sharp\circ  {}^pf_*\circ i\ ,\quad f^*:=i^\sharp\circ {}^pf^*\circ i$$
and obtain an adjoint pair
$$f^*:\Sh\bY\leftrightarrows\Sh\bX:f_*\ .$$
Note that
${}^pf_*$ preserves sheaves (see \cite[Lemma 2.13]{bss}).
The right-adjoint functor $f_*:\Sh_\Ab\bX\to \Sh_\Ab\bY$ is left exact and therefore admits a right-derived functor
$$Rf_*:D^+(\Sh_\Ab\bX)\to D^+(\Sh_\Ab\bY)$$
between the bounded below derived categories.

\subsection{}

If $g:Y\to Z$ is a second morphism of locally compact stacks, then  we have natural isomorphisms
of functors
$$(g\circ f)_*\cong g_*\circ f_*\ ,\quad f^*\circ g^*\cong (g\circ f)^*$$
(see \ref{uefhewiufuwefzzz}). 
Furthermore, we have
$$Rg_*\circ Rf_*\cong R(g\circ f)_*$$ on the level of bounded below derived categories by Lemma
\ref{keykey}.
The relation $ f^*\circ g^*\cong (g\circ f)^*$ descends to the derived categories if the pull-back functors are exact, e.g. if $f$ and $g$ have local sections (see \ref{sharp}).
 These facts generalize corresponding results shown in \cite{bss}.

\subsection{}\label{sharp}\label{system200}\label{lem:map_on_pr_vers_map_on_sh} \label{prexact}\label{lem:identify_star_sharp}
Let $f:G\to H$ be a morphism between topological stacks which has local sections.
It induces a morphism between sites
$f_\sharp:\bG\to \bH$ by composition. On objects it is given by
$f_\sharp(U\to G):=(U\to G\to H)$
(we will often use the short hand $U$ for $(U\to G)$ and write $f_\sharp U$).
In fact, since $U\to \bG$ and $f$ have local sections, the composition $U\to H$
has local sections. Furthermore, the map $U\to H$ from a space to a topological
stack is representable by Lemma \ref{lem:representability}. 
One checks that $f_\sharp$ maps covering families to covering families and preserves the fiber products
as in \cite[1.2.2]{MR1317816}.

If $f:G\to H$ has local sections, then 
the functor
$f^*:\Sh \bH\to \Sh \bG$ is the pull-back  $f^*=(f_\sharp)^*$ associated to a morphism of sites.  Explicitly it is given by 
$f^*F(U):=F(f_\sharp U)$, compare Lemma \cite[2.7]{bss}.
In  addition, the functor 
 $f^*:\Sh\bH\to\Sh\bG$ is exact (see \cite[2.5.9]{bss})
and preserves flat sheaves of abelian groups.

\begin{lem}\label{derivedadj}
If  $f:X\to Y$ is a morphism between locally compact stacks which has local sections, then we have the derived
adjunction
$$f^*\colon D^+(\Sh_\Ab\bY)\leftrightarrows D^+(\Sh_\Ab\bX)\colon Rf_*\ .$$
\end{lem}
\begin{proof}
 Since $f^*$ is exact its right adjoint $f_*$ preserves injectives.
If $G\in C^+(\Sh_\Ab\bX)$ is a complex of injectives and $F\in C^+(\Sh_\Ab\bY)$, then we have
\begin{multline*}
  R\Hom_{\Sh_\Ab\bY}(F,Rf_*(G))\cong \Hom_{\Sh_\Ab\bY}(F,f_*(G))\\
\cong
  \Hom_{\Sh_\Ab\bX}(f^*(F),G)\cong R\Hom_{\Sh_\Ab\bX}(f^*(F),G)\ .$$
\end{multline*}
This implies the assertion. \end{proof}


\subsection{} 

\begin{lem}\label{qwuidiuwqdwqdwqd}
Let $X$ be a locally compact stack. If $C,B\to X$ are maps from locally compact spaces, then  $C\times_XB$ is locally compact.
\end{lem}
\begin{proof} By assumption $X$ is locally compact so that we can chose an atlas  $A\to X$   such that
$A$ and $A\times_XA$ are locally compact. 
Since $A\to X$ is surjective and has local sections, there
exists an open covering $(B_i)$ of $B$ such that we have lifts
$$\xymatrix{&&A\ar[d]\\B_i\ar@{.>}[urr]\ar[r]&B\ar[r]&X}\ .$$
Then $(A\times_X B_i)$ is an open covering of $A\times_XB$.
In order to show that $A\times_XB$ is locally compact it suffices
to show that the space $A\times_XB_i$ is locally compact. By 
$A\times_XB_i\cong (A\times_XA)\times_AB_i\subseteq A\times_XA\times B_i$, this space is a closed
(note that $A$ is Hausdorff) subspace of a locally compact space and hence itself locally compact.

The same argument shows that $C\times_XA$ is locally compact.
We now write
$C\times_XB_i\cong (C\times_XA)\times_A B_i\subseteq  (C\times_XA) \times B_i$
in order to see that $C\times_XB_i$ is locally compact. Since
$(C\times_XB_i)$ is an open covering of $C\times_XB$ we conclude that
$C\times_XB$ is locally compact.
\end{proof}

\subsection{}

Let $f:X\to Y$ be a morphism between locally compact stacks.
\begin{lem}\label{obensharp}\label{iwdoi2d23d}
If $f$ is representable,  then it induces a morphism
of sites $f^\sharp:\bY\to \bX$ given by $f^\sharp(V\to Y):=(X\times_YV\to X)$.
\end{lem}
\begin{proof}
Let $B\to X$ be a locally compact atlas.
We consider $(V\to Y)\in \bY$ and form the   diagram of Cartesian squares  $$\xymatrix{V\times_YB\ar[r]\ar[d]&B\ar[d]\\U\ar[r]\ar[d]&X\ar[d]^f\\
V\ar[r]&Y}\ .$$
In order to check that $(U\to X)\in \bX$ we must show that $U$ is locally compact. 
Since $B\to X$ is surjective and has local sections we see  that
 $V\times_YB\to U$ is surjective and has local sections, too.
Since $Y$ is locally compact we see by Lemma \ref{qwuidiuwqdwqdwqd} that
$V\times_YB$ is locally compact.
Let $u\in U$ and $W\subseteq U$ be a neighborhood of $u$ such that there exists a section 
$$\xymatrix{&V\times_YB\ar[d]^\pi\\W\ar@{.>}[ur]^s\ar[r]&U}\ .$$
Let $K\subseteq \pi^{-1}(W)$ be a compact neighborhood of $s(u)$.
Then $s^{-1}(K)$ is a compact neighborhood of $u$.
Indeed, $s^{-1}(K)$ is a closed subset of the compact set $\pi(K)$.

It is easy to see  that $f^\sharp$ maps covering families to covering families and preserves the fiber products required for a morphism of sites, see
 \cite[1.2.2]{MR1317816}.
\end{proof}

If $f:X\to Y$ is a representable morphism between locally compact stacks, then
 we have the relations $f^*=(f^\sharp)_*:\Sh\bY\to \Sh\bX$
and $f_*=(f^\sharp)^*:\Sh \bX\to \Sh\bY$
, see \cite[Lemma 2.9]{bss}.

\subsection{}\label{obstr1a}\label{desc_sheaves_on_U}

Let $X$ be a topological stack and $(U\to X)\in \bX$. Let $(U)$ denote the site whose objects 
and morphisms are the open subsets of $U$ and inclusions, and whose coverings are coverings by families of open subsets.
We have restriction functors
$\nu_U:\Sh\bX  \to \Sh(U)$ and 
${}^p\nu_U:\Pr\bX\to \Pr(U)$. For $F\in \Sh\bX$ we also  write $\nu_U(F)=:F_U$. We have the following assertions,
most of which are straightforward to prove. 

\begin{enumerate}
\item Let $i^\sharp$ and $i^\sharp_U$ denote the sheafification functors on the sites $\bX$ and $(U)$. Then we have a natural isomorphism $$i_U^\sharp\circ {}^p\nu_U\cong \nu_U\circ i^\sharp\ ,$$ see
\cite[Lemma 2.4.7]{bss} 
\item Let $F\in \Sh\bX  $. If $f\colon U\to V$ is a morphism (\ref{gsdhasddas}) 
in $\bX  $, then we have a natural map $f^*F_V\to F_U$.
\item There is a one-to one correspondence of sheaves $F\in \Sh\bX  $ on the one hand,  and of collections $(F_U)_{(U\to X)\in \bX  }$ of sheaves $F_U\in \Sh(U)$ together with functorial maps $f^*F_V\to F_U$ for all morphisms $f\colon U\to V$ in $\bX$ on the other hand.
\item Let $F,G\in \Sh\bX  $. There is a one-to-one correspondence between compatible collections of morphisms $g_U\colon F_U\to G_U$ for  
all $(U\to X)\in \bX $ and maps $g\colon F\to G$.
\item If $F,G\in \Sh\bX  $ or
$F,G\in D^+(\Sh_\Ab\bX  )$, then
a map $F\to G$ is an isomorphism if and only if the induced map $F_U\to G_U$ is an isomorphism for all
$(U\to X)\in \bX  $.
\item Let $f\colon X\to Y$ be a representable map of locally compact stacks, $(A\to
  Y)\in \bY  $  and $(B:=A\times_YX\to X)\in \bX$.  Let $g\colon B\to A$ be the projection onto the first factor and $g_*:\Sh (B)\to \Sh(A)$. Then
  we have for $F\in \Sh\bX  $ or $G\in D^+(\Sh_\Ab\bX)$
$$(f_*F)_A\cong g_*(F_B)\ ,\quad (Rf_*G)_A\cong Rg_*(G_B)\ .$$
The second isomorphism follows from the first using the fact that the restriction $\nu_B$ preserves flabby or even injective sheaves (see Lemma \ref{hwfwefewfw}).
\item If $f\colon X\to Y$ is a map of topological stacks which has local sections, $(B\to X)\in \bX  $, then we have $(B\to X\to Y)\in \bY  $ and for $F\in\Sh\bY  $ 
$$(f^*F)_B\cong F_B\ .$$ 
\item
The collection of restriction functors $(\nu_U)_{(U\to X)\in \bX}$ detects flabby (flasque, flat) sheaves (see Definition \ref{def:flabby}), i.e.
a sheaf $F\in \Sh_\Ab\bX$ is flabby (flasque, flat)  if and only if $F_U\in \Sh_\Ab(U)$ is flabby for all $(U\to X) \in \bX$ (compare \ref{flatdetect} for the flat case). 
\item The collection of restriction functors $(\nu_U)_{(U\to X)\in \bX}$ detects exact sequences, i.e.
a sequence $F\to G\to H$ of sheaves of abelian groups on $\bX$ is exact if and only if
$F_U\to G_U\to H_U$ is exact for all $(U\to X)\in \bX$. 
\end{enumerate}

\begin{lem}\label{hwfwefewfw}
Let $(U\to X)\in \bX$.
The functor $\nu_U:\Sh_\Ab\bX\to \Sh_\Ab(U)$ preserves injective sheaves.
\end{lem}
\begin{proof}
We show that $\nu_U$ has an exact  left adjoint $\nu_\Z^U:\Sh_\Ab(U)\to \Sh_\Ab\bX$. 
We first show that the restriction functor ${}^p\nu_U:\Pr_\Ab\bX\to \Pr_\Ab(U)$ fits into an adjoint pair
$${}^p\nu_\Z^U:\Pr_\Ab(U)\leftrightarrows \Pr_\Ab\bX:{}^p\nu_U\ .$$
The left-adjoint is given by
 $${}^p\nu_\Z^U(F)(A\to X):= \colim F(V)\ ,$$
where the colimit is taken over the category of diagrams
$$\xymatrix{V\ar[d]&A\ar@{.>}[dl]^\phi\ar[d]\ar[l]\\U\ar [r]&X}\ ,$$
where $V\to U$ is the embedding of an open subset.
As explained in \cite[II.3.18]{MR559531} we have a decomposition of this category into a union of  categories $S(\phi)$ with  $\phi\in \Hom_\bX((A\to X),(U\to X))$. The category $S(\phi)$ is the category of open neighborhoods of $\phi(A)$ and their inclusions.  It is cofiltered.
Therefore  $F\mapsto \colim_{S(\phi)} F(V)$ preserves finite limits and is in particular left exact.
 This implies that  
${}^p\nu_\Z^U$ 
given by
$${}^p\nu_\Z^U(F)(A\to X)\cong \bigoplus_{\phi} \colim_{S(\phi)} F(V)  $$
 is left-exact, too.
We now get $\nu^U_\Z:=i^\sharp\circ {}^p\nu^U_\Z\circ i_U$.
As a left-adjoint it is right-exact. Since $i_U$ is left exact and $i^\sharp$ is exact, this composition is also left-exact.
\end{proof}

%
%
%

\subsection{}

%

\begin{lemma}\label{lem:pullpush}
  Consider the following Cartesian diagram in locally compact topological stacks
  \begin{equation*}
    \begin{CD}
     H @>v>> G\\
      @VVgV @VVfV\\
      Y @>u>> X     \end{CD}
  \end{equation*}
  In this situation the two canonical ways to define a natural transformation
$$u^*f_*\to g_*v^*\colon \Sh_\Ab(\bG)\to \Sh_\Ab(\bY)$$ 
give the same result, i.e. the diagram
\begin{equation}\label{hfiewufuwef}
\xymatrix{u^*f_*\ar@{=}[d]\ar[r]^{unit}&g_*g^*u^*f_*\ar[r]^{u g=f v}&g_*v^*f^*f_*\ar[r]^{\hspace{0.1cm}counit}&g_*v^*\ar@{=}[d]\\u^*f_*\ar[r]^{unit}&u^*f_*v_*v^*\ar[r]^{u  g=f v}&u^*u_*g_*v^*\ar[r]^{\hspace{0.2cm}counit}&g_*v^*}
\end{equation}
commutes. 
This transformation is functorial with respect to composition of Cartesian diagrams.

Moreover, 
if $u$ has local sections, then this transformation induces isomorphisms
  \begin{gather}
    u^*f_*\iso g_*v^*\colon \Sh_\Ab(\bG)\to \Sh_\Ab(\bY),\label{fewfwfwefewf}\\
    u^*Rf_* \iso Rg_* v^* \colon D^+\Sh_\Ab(\bG)\to D^+\Sh_\Ab(\bY).\label{wqdwqdqw}
  \end{gather}
If $u$ and $f$ have local sections, then we get 
  commutative diagrams
 $$
\xymatrix{&u_*\ar[rd]^{unit}\ar[dl]^{unit}&\\
u_*g_*g^*\ar[r]^\sim&f_*v_*g^*&f_*f^*u_*\ar[l]_\sim}, \xymatrix{&v_* &\\
f^*f_*v_*\ar[ur]^{counit}\ar[r]^\sim&f^*u_*g_*\ar[r]^\sim&v_*g^*g_*\ar[ul]^{counit}}$$
    
$$
\xymatrix{&u^*\ar[rd]^{unit}\ar[dl]^{unit}&\\
u^*f_*f^*\ar[r]^\sim&g_*v^*f^*&g_*g^*u^*\ar[l]_\sim}, \xymatrix{&v^* &\\
v^*f^*f_*\ar[ur]^{counit}\ar[r]^\sim&g^*u^*f_*\ar[r]^\sim&g^*g_*v^*\ar[ul]^{counit}}$$
   and
their derived versions, e.g.

\begin{equation}\label{ddwqzwqidqwdww}
 \xymatrix{&u^*\ar[rd]^{unit}\ar[dl]^{unit}&\\
u^*Rf_*f^*\ar[r]^\sim&Rg_*v^*f^*&Rg_*g^*u^*\ar[l]_\sim}\ ,
  \end{equation}
and also
   \begin{equation}\label{wdqidwqdwdwdewdwedwed}
\xymatrix{&&Ru_*u^*\ar[drr]^{unit}\ar[dll]_{unit}&&\\
Ru_*u^*Rf_*f^*\ar[r]^\sim&Ru_*Rg_*v^*f^*\ar[r]^\sim&Rf_*Rv_*v^*f^*\ar[r]^\sim&Rf_*Rv_*g^*u^*&Rf_*f^*Ru_*u^*\ar[l]^\sim}  
   \end{equation}
\end{lemma}
\begin{proof}
\textit{Most of the following arguments and the large diagrams  were supplied by \textbf{A. Schneider}.}
For convenience we present a proof of (\ref{hfiewufuwef}), see also \cite[Expose XVII, Proposition 2.1.3]{deligne}. 
%
We first observe that
\begin{equation}\label{jdewkdwed}
\xymatrix{
v^*f^*f_*v_*\ar[r]^{counit}\ar[d]^\sim&v^*v_*\ar[d]_{counit}\\
(fv)^*(fv)_*\ar[r]^{counit}&\id
}
\end{equation}
commutes. Using this in addition to standard functorial properties 
we check that all squares in the following diagram commute:
$$\hspace{-0.8cm}
 \xymatrix{
u^*f_*\ar[r]^{unit}\ar@{=}[ddddd]&
g_*g^* u^*f_*\ar[r]^\sim\ar[d]^{unit}&g_*(ug)^*f_*\ar[r]^=\ar[d]^{unit}&
g_*(fv)^*f_*\ar[r]^\sim\ar[d]^{unit}&g_*v^* f^*f_*\ar[r]^{counit}\ar[d]^{unit}&
g_*v^* \ar[d]^{unit}\ar@/^1cm/[dd]^{\rm id}\\
&
g_*g^* u^*f_*v_*v^*\ar[r]^\sim\ar@{=}[d]&g_*(ug)^*f_*v_*v^*\ar[r]^=\ar[d]^\sim&
g_*(fv)^*f_*v_*v^*\ar[r]^\sim\ar[d]^\sim&g_*v^* f^*f_*v_*v^*\ar[r]^{counit}\ar[d]^\sim&
g_*v^* v_*v^*\ar[d]^{counit}\\
&
g_*g^* u^*f_*v_*v^*\ar[r]^\sim&g_*(ug)^*(fv)_*v^*\ar[r]^=&
g_*(fv)^*(fv)_*v^*\ar@{=}[r]&g_* (fv)^*(fv)_*v^*\ar[r]^{counit}&g_*v^*\\
&
g_*g^* u^*f_*v_*v^*\ar[r]^\sim\ar@{=}[u]&g_*(ug)^*(fv)_*v^*\ar[r]^=\ar@{=}[u]&
g_*(ug)^*(ug)_*v^*\ar@{=}[r]\ar[u]_=&g_* (ug)^*(ug)_*v^*\ar[r]^{counit}\ar[u]_=&g_*v^*\ar@{=}[u]\\
&
g_*g^* u^*f_*v_*v^*\ar[r]^\sim\ar@{=}[u]&g_*g^*u^*(fv)_*v^*\ar[r]^=\ar[u]_\sim&
g_*g^*u^*(ug)_*v^*\ar[r]^\sim\ar[u]_\sim&g_*g^* u^*u_*g^*v^*\ar[r]^{counit}\ar[u]_\sim&
g_*g^* g_*v^*\ar[u]_{counit}\\
u^*f_*\ar[r]^{unit}&
u^*f_*v_*v^*\ar[r]^\sim\ar[u]_{unit}&u^*(fv)_*v^*\ar[r]^=\ar[u]_{unit}&
u^*(ug)_*v^*\ar[r]^\sim\ar[u]_{unit}&u^*u_*g^*v^*\ar[r]^{counit}\ar[u]_{unit}&
g_*v^*\ar[u]_{unit}\ar@/_1cm/[uu]_{\rm id}.
}
$$ 
The two ways to go along the boundary from the upper left to lower right corner give the two maps $u^*f_*\to g_*v^*$ in question.

The isomorphism (\ref{fewfwfwefewf}) can be shown as in \cite[Lemma 2.16]{bss}, where the
assumption of smoothness of $u$ in  \cite{bss} corresponds to the assumption of local sections in the present setting. The derived version (\ref{wqdwqdqw}) can be shown using the simplicial models as in \cite [Lemma 2.43]{bss}.
Alternatively one can use the commutativity of the diagram asserted in Lemma \ref{uiqehewqdqwdwqdqd}
and the isomorphism (\ref{uiidwqdwqdwqd}).

We now show the compatibility of the units and counits with Cartesian diagrams. The arguments are purely formal and only use
that the functors involved occur as parts of adjoint pairs. We will only give the details for the two triangles involving derived functors. 
 If in addition to $u$ also  $f$ has local sections, then so has $g$. In this case we have the adjoint pairs $(f^*,Rf_*)$ and $(g^*,Rg_*)$.
In order to see (\ref{ddwqzwqidqwdww}) we must show that  
$$
\xymatrix{
u^*\ar[r]^{unit}\ar@/_1cm/[rrrrr]_{unit}&
u^*Rf_*f^*\ar[r]^{\Psi}&Rg_*v^*f^*\ar[r]^\sim&Rg_*(fv)^*\ar[r]^=&
Rg_*(ug)^*\ar[r]^\sim&Rg_*g^*u^*,
}
$$
commutes, 
where $\Psi:u^*Rf_*f^*\to Rg_*v^*f^*$ is induced by (\ref{wqdwqdqw}).
This is a consequence of the  commutativity of 
$$
\xymatrix{
u^*\ar[r]^{unit}&
u^*Rf_*f^*\ar[d]^{unit}\ar[rrr]^{\Psi}&&&
Rg_*v^*f^* \ar@/^2cm/[dd]\\
&
Rg_*g^* u^*f_*f^*\ar[r]^\sim&Rg_*(ug)^*Rf_*f^*\ar[r]^=&
Rg_*(fv)^*Rf_*f^*\ar[r]^\sim&Rg_*v^* f^*Rf_*f^*\ar[u]_{counit}\\
u^*\ar[r]^{unit}\ar@{=}[uu]&
Rg_*g^* u^*\ar[r]^\sim\ar[u]_{unit}&Rg_*(ug)^*\ar[u]_{ unit}\ar[r]^=\ar@/^0.5cm/[l]&
Rg_*(fv)^*\ar[u]_{unit}\ar[r]^\sim\ar@/^0.5cm/[l]&Rg_*v^* f^*\ar@/^0.5cm/[l]\ar[u]_{ unit}\ar@/_1cm/[uu]_{\rm id}
}
$$
which follows from standard functorial properties of units and counits.

The same properties are used in the proof of (\ref{wdqidwqdwdwdewdwedwed}) which is
represented by the boundary of the following big array of small commutative squares and triangles
\newpage
\pagestyle{empty}
\hspace{13cm}
{\tiny \begin{rotate}{270}
$
\hspace{-2cm}
\xymatrix{
&Rf_*f^*u_*u^*\ar[r]^{unit}\ar@/^0.5cm/[rrrrr]^\Phi\ar[d]^{unit}	&Rf_*f^*u_*Rg_*g^*u^*\ar[r]^\sim\ar[d]^{unit}		&Rf_*f^*R(ug)_*g^*u^*\ar[r]^=	\ar[d]^{unit}		&Rf_*f^*R(fv)_*g^*u^*\ar[r]^\sim\ar[d]^{unit}		&Rf_*f^*Rf_*v_*g^*u^*\ar[r]^{counit}\ar[d]^{unit}		&Rf_*v_*g^*u^*\ar@{=}[rd]\ar[d]^{unit}&\\
&Rf_*f^*u_*u^*Rf_*f^*\ar[r]^{unit}\ar@{=}[d]&Rf_*f^*u_*Rg_*g^*u^*Rf_*f^*\ar[r]^\sim\ar@{=}[d]	&Rf_*f^*R(ug)_*g^*u^*Rf_*f^*\ar[r]^=	\ar[d]^\sim	&Rf_*f^*R(fv)_*g^*u^*Rf_*f^*\ar[r]^\sim\ar[d]^\sim	&Rf_*f^*Rf_*v_*g^*u^*Rf_*f^*\ar[r]^{counit}\ar[d]^\sim	&Rf_*v_*g^*u^*Rf_*f^*\ar[d]^\sim
&Rf_*v_*g^*u^*\ar[d]^\sim\ar[l]_{unit}\\
&Rf_*f^*u_*u^*Rf_*f^*\ar[r]^{unit}\ar@{=}[d]&Rf_*f^*u_*Rg_*g^*u^*Rf_*f^*\ar[r]^\sim\ar@{=}[d]	&Rf_*f^*R(ug)_*(ug)^*Rf_*f^*\ar[r]^=	\ar@{=}[d]		&Rf_*f^*R(fv)_*(ug)^*Rf_*f^*\ar[r]^\sim\ar[d]^=		&Rf_*f^*Rf_*v_*(ug)^*Rf_*f^*\ar[r]^{counit}\ar[d]^=		&Rf_*v_*(ug)^*Rf_*f^*\ar[d]^=
&Rf_*v_*(ug)^*\ar[d]^=\ar[l]_{unit}\\
&Rf_*f^*u_*u^*Rf_*f^*\ar[r]^{unit}\ar@{=}[d]&Rf_*f^*u_*Rg_*g^*u^*Rf_*f^*\ar[r]^\sim\ar@{=}[d]	&Rf_*f^*R(ug)_*(ug)^*Rf_*f^*\ar[r]^=\ar@{=}[d]		&Rf_*f^*R(fv)_*(fv)^*Rf_*f^*\ar[r]^\sim\ar@{=}[d]	&Rf_*f^*Rf_*v_*(fv)^*Rf_*f^*\ar[r]^{counit}\ar[d]^\sim	&Rf_*v_*(fv)^*Rf_*f^*\ar[d]^\sim
&Rf_*v_*(fv)^*\ar[l]_{unit}\ar[d]^\sim\\
&Rf_*f^*u_*u^*Rf_*f^*\ar[r]^{unit}\ar@{=}[dd]&Rf_*f^*u_*Rg_*g^*u^*Rf_*f^*\ar[r]^\sim\ar@{=}[dd]	&Rf_*f^*R(ug)_*(ug)^*Rf_*f^*\ar[r]^=\ar@{=}[dd]		&Rf_*f^*R(fv)_*(fv)^*Rf_*f^*\ar[r]^\sim\ar@{=}[dd]	&Rf_*f^*Rf_*v_*v^*f^*Rf_*f^*\ar[r]^{counit}\ar@{=}[dd]	&Rf_*v_*v^*f^*Rf_*f^*\ar[dr]^{counit}
&Rf_*v_*v^*f^*\ar[l]_{unit}\ar[d]^{\rm id}\\
u_*u^*\ar@/^0.7cm/[uuuuur]^{unit}\ar@/_0.7cm/[dddddr]_{unit}&&&&&&&Rf_*v_*v^*f^*\\
&Rf_*f^*u_*u^*Rf_*f^*\ar[r]^{unit}\ar@{=}[d]&Rf_*f^*u_*Rg_*g^*u^*Rf_*f^*\ar[r]^\sim\ar@{=}[d]	&Rf_*f^*R(ug)_*(ug)^*Rf_*f^*\ar[r]^=\ar@{=}[d]		&Rf_*f^*R(fv)_*(fv)^*Rf_*f^*\ar[r]^\sim\ar@{=}[d]	&Rf_*f^*Rf_*v_*v^*f^*Rf_*f^*\ar[r]^{counit}			&Rf_*f^*Rf_*v_*v^*f^*\ar[ur]^{counit}
&Rf_*v_*v^*f^*\ar[l]_{unit}\ar[u]_{\rm id}\\
&Rf_*f^*u_*u^*Rf_*f^*\ar[r]^{unit}\ar@{=}[d]&Rf_*f^*u_*Rg_*g^*u^*Rf_*f^*\ar[r]^\sim\ar@{=}[d]	&Rf_*f^*R(ug)_*(ug)^*Rf_*f^*\ar[r]^=\ar@{=}[d]		&Rf_*f^*R(fv)_*(fv)^*Rf_*f^*\ar[r]^\sim			&Rf_*f^*R(fv)_*v^*f^*Rf_*f^*\ar[r]^{counit}\ar[u]_\sim	&Rf_*f^*R(fv)_*v^*f^*\ar[u]_\sim
&R(fv)_*v^*f^*\ar[u]_\sim\ar[l]_{unit}\\
&Rf_*f^*u_*u^*Rf_*f^*\ar[r]^{unit}\ar@{=}[d]&Rf_*f^*u_*Rg_*g^*u^*Rf_*f^*\ar[r]^\sim\ar@{=}[d]	&Rf_*f^*R(ug)_*(ug)^*Rf_*f^*\ar[r]^=				&Rf_*f^*R(ug)_*(fv)^*Rf_*f^*\ar[r]^\sim\ar[u]_=		&Rf_*f^*R(ug)_*v^*f^*Rf_*f^*\ar[r]^{counit}\ar[u]_=		&Rf_*f^*R(ug)_*v^*f^*\ar[u]_=
&R(ug)_*v^*f^*\ar[u]_=\ar[l]_{unit}\\
&Rf_*f^*u_*u^*Rf_*f^*\ar[r]^{unit}		&Rf_*f^*u_*Rg_*g^*u^*Rf_*f^*\ar[r]^\sim			&Rf_*f^*u_*Rg_*(ug)^*Rf_*f^*\ar[r]^=\ar[u]_\sim	&Rf_*f^*u_*Rg_*(fv)^*Rf_*f^*\ar[r]^\sim\ar[u]_\sim	&Rf_*f^*u_*Rg_*v^*f^*Rf_*f^*\ar[r]^{counit}\ar[u]_\sim	&Rf_*f^*u_*Rg_*v^*f^*\ar[u]_\sim
&u_*Rg_*v^*f^*\ar[u]_\sim\ar[l]_{unit}\\
&u_*u^*Rf_*f^*\ar[r]^{unit}\ar@/_0.5cm/[rrrrr]^\Psi\ar[u]_{unit}	&u_*Rg_*g^*u^*Rf_*f^*\ar[r]^\sim	\ar[u]_{unit}	&u_*Rg_*(ug)^*Rf_*f^*\ar[r]^=\ar[u]_{unit}		&u_*Rg_*(fv)^*Rf_*f^*\ar[r]^\sim\ar[u]_{unit}		&u_*Rg_*v^*f^*Rf_*f^*\ar[r]^{counit}\ar[u]_{unit}		
&u_*Rg_*v^*f^*\ar@{=}[ru]\ar[u]_{unit}&\\
}
$
\end{rotate}}
\newpage
\pagestyle{plain}


\end{proof}

\section{Tensor products and the projection formula}

\subsection{}

We consider a Grothendieck site $\bX$ and a commutative ring $R$. 
The goal of the present Subsection is to discuss aspects of the closed monoidal structures on
the categories of presheaves $\Pr_{R-\Mod}\bX$ and sheaves $\Sh_{R-\Mod}\bX$ of $R$-modules on $\bX$. The material is standard, but we need to understand in detail the relation between the sheaf and presheaf versions in order to show the compatibility with the operations
induced by a morphism of stacks.

\subsection{}

Let $F,G\in \Pr_{R-\Mod}\bX$ be presheaves of $R$-modules.
The tensor product $F\otimes^p G\in \Pr_{R-\Mod}\bX$ is defined as the presheaf which associates to $(U\to X)$ the $R$-module
$F(U)\otimes^p_R G(U)$. In this way $\Pr_{R-\Mod}\bX$ becomes a symmetric monoidal category.

Since colimits of presheaves are defined objectwise
we have for a diagram of presheaves of $R$-modules $(F_i)_{i\in I}$ that
$$\colim_{i\in I} (F_i\otimes^p_RG)\cong (\colim_{i\in I}F_i)\otimes^p_RG\ .$$

\newcommand{\uHom}{\underline{\tt Hom}}
\subsection{}\label{system202}

For $U\in \bX$ let $h_U\in \Pr\bX$ denote the presheaf represented by $U$ and
$h_U^R\in \Pr_{R-\Mod}\bX$ be the presheaf of $R$-modules generated by $h_U$.
Let $F,G\in \Pr_{R-\Mod}\bX$. We define the presheaf
$$\uHom^p(F,G)\in {\Pr}_{R-\Mod}\bX$$ by
$$\uHom^p(F,G)(U):=\Hom_{\Pr_{R-\Mod}\bX}(h_U^R\otimes^p F,G)\ .$$

The topology of the site of a locally compact stack is sub-canonical. Hence, in this case $h_U$ is actually a sheaf.
But even in the case of a sub-canonical topology $h^R_U$ is only a presheaf, in general.

If $U\to V$ is a morphism in $\bX$, then 
$\uHom^p(F,G)(V)\to \uHom^p(F,G)(U)$ is induced by the morphism $h_U\to h_V$.
If $H\in  \Pr_{R-\Mod}\bX$, then we have
\begin{eqnarray*}
\Hom_{ \Pr_{R-\Mod}\bX}(H,\uHom^p(F,G))&\cong&\Hom_{ \Pr_{R-\Mod}\bX}(\colim_{h^R_V\to H}h_V^R, \uHom^p(F,G))\\
&\cong&
\lim_{h^R_V\to H}\Hom_{ \Pr_{R-\Mod}\bX}(h_V^R, \uHom^p(F,G))\\
&\cong&\lim_{h^R_V\to H} \uHom^p(F,G)(V)\\
&=&\lim_{h^R_V\to H} \Hom_{ \Pr_{R-\Mod}\bX}(h_V^R \otimes^p F,G)\\
&\cong&\Hom_{ \Pr_{R-\Mod}\bX}(\colim_{h^R_V\to H} (h_V^R \otimes^p F),G)\\
&\cong&\Hom_{ \Pr_{R-\Mod}\bX}((\colim_{h^R_V\to H} h_V^R) \otimes^p F,G)\\
&\cong&\Hom_{ \Pr_{R-\Mod}\bX}(H\otimes^p F,G)\end{eqnarray*}
In other words, the pair  $(\otimes^p,\uHom^p)$ together with this natural isomorphism defines a closed symmetric monoidal structure on $\Pr_{R-\Mod}\bX$.
In particular, if $(F_i)_{i\in I}$ is a diagram of presheaves, then we have
\begin{equation}\label{mb12}\uHom^p(\colim_{i\in I} F_i,G)\cong \lim_{i\in I}\uHom^p(F_i,G)\ .\end{equation} 

\subsection{}

An element of 
 $$\uHom(F,G)(U)=\Hom_{ \Pr_{R-\Mod}\bX}(h_U^R\otimes^p F,G) $$
is given by a collection of $R$-linear maps $(\phi_{V\to U}:F(V)\to G(V))_{(V\to U)\in \bX/U}$ such that for a morphism
$(W\to U)\mapsto (V\to U)$ in $\bX/U$ the diagram
$$\xymatrix{F(V)\ar[r]\ar[d]^{\phi_{V\to U}}&F(W)\ar[d]^{\phi_{W\to U}}\\G(V)\ar[r]&G(W)}$$ commutes. Therefore
$$\uHom(F,G)(U)\cong \Hom_{ \Pr_{R-\Mod}\bX/U}(F_{|U},G_{|U})\ .$$
\begin{lem}\label{vorh1}
If $G$ is a sheaf, then $\uHom(F,G)$ is a sheaf.
\end{lem}
\begin{proof}
Let $U\in \bX$ and $(U_i\to U)_{i\in I}$ be a covering.
In order to simplify the notation we consider $V:=\sqcup_{i\in I} U_i$.
We must show that the sequence
$$0\to \uHom(F,G)(U)\to  \uHom(F,G)(V)\to  \uHom(F,G)(V\times_UV)$$
is exact.

Let $\psi\in  \Hom_{\Pr_{R-\Mod}\bX/U}(F_{|U},G_{|U})$ be such that
its restriction to $V$ vanishes.
If $(W\to U)\in \bX/U$, then $W\times_UV\to W$ is a covering of $W$,
and $\pr_W^*:G(W)\to G(W\times_UV)$ is injective since $G$ is a sheaf.
In view of the commutative diagram 
$$\xymatrix{F(W)\ar[r]^{\pr_W^*}\ar[d]^{\psi_W}&F(W\times_UV)\ar[d]^{(\psi_{|V})_{W\times_UV}}\\G(W)\ar[r]^{\pr_W^*}&G(W\times_UV)}$$
we see that $\psi_W=0$.

Let now $\phi\in  \Hom_{\Pr_{R-\Mod}\bX/V}(F_{|V},G_{|V})$
be such that the induced map 
$$\Phi\in \Hom_{\Pr_{R-\Mod}\bX/(V\times_UV)}(F_{|V\times_UV},G_{|V\times_UV})$$
vanishes. We will construct
$\psi\in  \Hom_{\Pr_{R-\Mod}\bX/U}(F_{|U},G_{|U})$ such that $\psi_{|V}=\phi$.
Let $(W\to U)\in \bX/U$ and $f\in F_{|U}(W\to U)=F(W)$.
Then $W\times_UV\to W$ is a covering of $W$ and $\pr_W^*f\in F_{|V}(W\times_UV\to V)=F(W\times_UV)$. We get an element
$$\phi_{W\times_UV\to V} (\pr_W^*(f))\in G(W\times_UV)= G_{|V}(W\times_UV\to V)\ .$$
Note that
$(W\times_UV)\times_W(W\times_UV)\cong W\times_U(V\times_UV)$. The difference
of the pull-backs of $\phi_{W\times_UV\to V} (\pr_W^*(f))$ with respect to the two projections to $W\times_UV$ induces
$$\Phi_{W\times_U(V\times_UV)}(\pr_W^*(f))=0\in G((W\times_UV)\times_W(W\times_UV))\ .$$
Again, since $G$ is a sheaf
there is a unique element $\psi_W(f)\in G(W)$ such that
$$\psi_W(f)_{|W\times_UV}=\phi_{W\times_UV\to V} (\pr_W^*(f))\ .$$
The morphism $\psi$ is now given by the collection $(\psi_W)_{(W\to U)\in \bX/U}$.
\end{proof}

\subsection{}\label{system101}

If $F,G\in \Sh_{R-\Mod}\bX$, then we define
$F\otimes G\in  \Sh_{R-\Mod}\bX$ to be
$$F\otimes G :=i^\sharp(i(F)\otimes^p i(G))\ .$$ 
We furthermore define
$$\uHom(F,G):=i^\sharp\uHom^p(i(F),i(G))\ .$$
Using the fact \ref{vorh1} that $\uHom^p(i(F),i(G))$ is a sheaf at the isomorphism marked by $!$ we get for every $H\in \Sh_{R-\Mod}\bX$ that 
\begin{eqnarray*}
\Hom_{\Sh_{R-\Mod}\bX}(H\otimes F,G)&\cong&\Hom_{\Sh_{R-\Mod}\bX}(i^\sharp(i(H)\otimes^p i(F),G)\\
&\cong&\Hom_{\Pr_{R-\Mod}\bX}(i(H)\otimes^p i(F),i(G))\\
&\cong&\Hom_{\Pr_{R-\Mod}\bX}(i(H),\uHom^p(i(F),i(G))\\
&\stackrel{!}{\cong}&\Hom_{\Pr_{R-\Mod}\bX}(i(H),i\circ i^\sharp (\uHom^p(i(F),i(G))))\\
&\cong&\Hom_{\Sh_{R-\Mod}\bX}(i^\sharp\circ i(H),\uHom(F,G))\\
&\cong&\Hom_{\Sh_{R-\Mod}\bX}(H,\uHom(F,G))\ .\end{eqnarray*}
In other words, the pair 
$(\otimes,\uHom)$ together with this natural isomorphism make $\Sh_{R-\Mod}\bX$ into a closed symmetric monoidal category.

\subsection{}\label{flatdetect}

Let $F,G\in \Sh_{R-\Mod}\bX$ and $(U\to X)\in \bX$. Then we have
$$(F\otimes G)_U\cong F_U\otimes G_U\ .$$
Indeed, this follows from the fact that sheafification commutes with the restriction from the site $\bX$ to the site $(U)$, see \ref{obstr1a}.
Since the collection of functors $(\nu_U)_{(U\to X)\in \bX}$ detects exact sequences
it now follows that a sheaf $F\in  \Sh_{R-\Mod}\bX$ is flat if and only if $F_U\in \Sh_{R-\Mod}(U)$ is flat for all $(U\to X)\in \bX$. This fact was claimed in \ref{obstr1a}.

\subsection{}

\begin{lem}\label{kion}
For $F,G\in \Pr_{R-\Mod}\bX$ we  have
$i^\sharp(F\otimes^p G)\cong i^\sharp(F)\otimes i^\sharp(G)$.
\end{lem}
\begin{proof}
This follows from (we omit the functor $i$ at various places in order to simplify the formula)
\begin{eqnarray*}
\Hom_{\Sh_{R-\Mod}\bX}( i^\sharp(F\otimes^p G),H)&\cong&
\Hom_{\Pr_{R-\Mod}\bX}( F\otimes^p G,H)\\
&\cong&\Hom_{\Pr_{R-\Mod}\bX}(F,\uHom^p(G,H))\\
&\stackrel{!}{\cong}&\Hom_{\Pr_{R-\Mod}\bX}(i^\sharp(F),\uHom^p(G,H))\\
&\cong&\Hom_{\Pr_{R-\Mod}\bX}(i^\sharp(F)\otimes^p G,H)\\
&\cong&\Hom_{\Pr_{R-\Mod}\bX}(G,\uHom^p(i^\sharp F,H))\\
&\stackrel{!}{\cong}&\Hom_{\Pr_{R-\Mod}\bX}(i^\sharp G,\uHom^p(i^\sharp F,H))\\
&\cong&\Hom_{\Pr_{R-\Mod}\bX}(i^\sharp G\otimes^p i^\sharp F,H)\\
&\cong&\Hom_{\Sh_{R-\Mod}\bX}(i^\sharp G\otimes i^\sharp F,H) 
 \end{eqnarray*}
for arbitrary $H\in \Sh_{R-\Mod}\bX$, where we use Lemma \ref{vorh1} at the isomorphisms marked by $!$.
\end{proof}
\subsection{}

Let $f\colon X\to Y$ be a morphism of  locally compact stacks. 
Let $\bX$ and $\bY$ be the sites associated to $X$ and $Y$.  Consider the adjoint pair of functors
$${}^pf^*\colon \Pr_{R-\Mod}\bY\leftrightarrows\Pr_{R-\Mod}\bX\colon {}^pf_*\ .$$
The proof of the following Lemma uses the product in $\bY$ 
described in \cite[Lemma 3.1]{bss}
in a specific way.
\begin{lem}
For $F,G\in \Pr_{R-\Mod}\bY$ we have a natural isomorphism
$${}^pf^*(F\otimes^pG)\cong {}^pf^*F\otimes^p{}^p f^*G\ .$$
\end{lem}
\begin{proof}
We use the notation introduced in \cite[2.1.4]{bss}.
For  $(U\to X)\in \bX$ we consider the category $U/\bY$ of diagrams $$\xymatrix{U\ar[r]\ar[d]&X\ar[d]\\V\ar[r]&Y}\ .$$
The functor ${}^pf^*$ is defined in \cite[Definition 2.3]{bss}
as a colimit over this category.

We consider  the diagonal functor
$U/\bY\to U/\bY\times U/\bY$ given on objects by
$$\xymatrix{U\ar[r]\ar[d]&X\ar[d]\\V\ar[r]&Y}\hspace{0.5cm}\mapsto\hspace{0.5cm} (\xymatrix{U\ar[r]\ar[d]&X\ar[d]\\V\ar[r]&Y},\xymatrix{U\ar[r]\ar[d]&X\ar[d]\\V\ar[r]&Y})\ .$$
In view of the definition of ${}^pf^*$ by colimits
it induces a map
$${}^pf^*(F\otimes^pG)\to {}^pf^*F\otimes^p {}^pf^*G\ .$$
In the other direction we have the functor
$U/\bY\times U/\bY\to U/\bY$ given by
$$(\xymatrix{U\ar[r]\ar[d]&X\ar[d]\\V\ar[r]&Y},\xymatrix{U\ar[r]\ar[d]&X\ar[d]\\V^\prime\ar[r]&Y})\hspace{0.5cm}\mapsto\hspace{0.5cm} \xymatrix{U\ar[r]\ar[d]&X\ar[d]\\V\times_Y V^\prime\ar[r]&Y}\ .$$
This together with the projections
$V\times_Y V^\prime\to V$ and $V\times_Y V^\prime\to V^\prime$ it induces the inverse map
$${}^pf^*F\otimes^p {}^pf^*G \to {}^pf^*(F\otimes^pG)\ .$$
\end{proof}

\subsection{}

Let $f\colon X\to Y$ be a morphism of locally compact stacks.
\begin{lem}\label{tens-pres}
For $F,G\in \Sh_{R-\Mod}\bY$ we have a natural isomorphism
$$f^*(F\otimes G)\cong f^*F\otimes f^*G\ .$$
\end{lem}
\begin{proof}
For $H\in \Sh_{R-\Mod} \bX$, using the fact that ${}^pf_*$ preserves sheaves (see \ref{sec:morph_of_sheaves}) and Lemma \ref{kion}, we have
\begin{eqnarray*}
\Hom_{\Sh_{R-\Mod}\bX}(f^*(F\otimes G),H)&\cong&\Hom_{\Sh_{R-\Mod}\bX}(F\otimes G,f_*(H))\\
&\cong&\Hom_{\Sh_{R-\Mod}\bY}(i^\sharp(i(F)\otimes^p i(G)),i^\sharp \circ f_*^p\circ i(H))\\
&\cong&\Hom_{\Pr_{R-\Mod}\bY}((i(F)\otimes^p i(G)), f_*^p\circ i(H))\\
&\cong&\Hom_{\Pr_{R-\Mod}\bX}({}^pf^*(i(F)\otimes^p i(G)), i(H))\\
&\cong&\Hom_{\Pr_{R-\Mod}\bX}({}^pf^*\circ i(F)\otimes^p {}^pf^*\circ i(G), i(H))\\
&\cong&\Hom_{\Sh_{R-\Mod}\bX}(i^\sharp({}^pf^*\circ i(F)\otimes^p {}^pf^*\circ i(G)), H)\\
&\cong&\Hom_{\Sh_{R-\Mod}\bX}(f^*(F)\otimes f^*(G), H)
\end{eqnarray*}
\end{proof}

\subsection{}

For a derived version of Lemma \ref{tens-pres} we assume that the morphism  $f:X\to Y$ of locally compact stacks has local sections. For simplicity we only consider the case
$R=\Z$, i.e. sheaves of abelian groups (finite cohomological dimension of $R$ would suffice).
Then the exact functor $f^*=(f_\sharp)^*$ preserves torsion-free sheaves of abelian groups. Since the derived tensor product can be calculated using torsion-free resolutions we get the corollary 
\begin{kor} \label{tens-pres1} If $f:X\to Y$ has local sections, then
for $F,G\in D^+(\Sh_\Ab\bY)$ we have a natural isomorphism
$$f^*(F\otimes^L G)\cong f^*F\otimes^L f^*G\ .$$ 
\end{kor}
of Lemma \ref{tens-pres}.

\subsection{}

Let $f\colon X\to Y$ be a morphism of locally compact stacks.
\begin{lem}
For $F\in \Sh_{R-\Mod}\bY$ and $G\in \Sh_{R-\Mod}\bX$ we have a natural isomorphism
$$\uHom(F,f_*G)\cong f_*\uHom(f^*F,G)$$
in  $\Sh_{R-\Mod}\bY$
\end{lem}
\begin{proof}
For any $T\in \Sh_{R-\Mod}\bY$ we calculate
\begin{eqnarray*}
\Hom_{\Sh_{R-\Mod}\bY}(T,f_*\uHom(f^*F,G))&\cong&\Hom_{\Sh_{R-\Mod}\bX}(f^*T,\uHom(f^*F,G))\\
&\cong& \Hom_{\Sh_{R-\Mod}\bX}(f^*T\otimes f^*F,G)\\
&\cong& \Hom_{\Sh_{R-\Mod}\bX}(f^*(T\otimes F),G)\\
&\cong&\Hom_{\Sh_{R-\Mod}\bY}(T\otimes F,f_*G)\\
&\cong&\Hom_{\Sh_{R-\Mod}\bY}(T,\uHom(F,f_*G))
\end{eqnarray*}
\end{proof}

\subsection{}

Let $f\colon X\to Y$ be a morphism of locally compact stacks.
\begin{lem}\label{system80}
For $F\in \Sh_{R-\Mod}\bY$ and $G\in \Sh_{R-\Mod}\bX$ we have a natural morphism
$$f_*G\otimes F\to f_*(G\otimes f^* F)\ .$$
\end{lem}
\begin{proof}
The transformation is the image of the identity under the following chain of maps, where the first is induced by the counit $f^*\circ f_*\to \id $ of the adjoint pair $(f^*,f_*)$, and the second isomorphism is given by Lemma \ref{tens-pres}.
\begin{eqnarray*}
\Hom_{\Sh_{R-\Mod}\bX}(G\otimes f^*F,G\otimes f^*F)&\to&
\Hom_{\Sh_{R-\Mod}\bX}(f^* f_*G\otimes f^*F,G\otimes f^*F)\\
&\cong&\Hom_{\Sh_{R-\Mod}\bX}(f^*(f_*G\otimes F),G\otimes f^*F)\\
&\cong&\Hom_{\Sh_{R-\Mod}\bY}(f_*G\otimes F,f_*(G\otimes f^*F))\ .
\end{eqnarray*}
\end{proof}

\begin{lem}
If $f$ has local sections, then for $F\in \Sh_{\Ab}\bY$ and $G\in \Sh_{\Ab}\bX$ we have a natural morphism
$$f_*G\otimes^L F\to f_*(G\otimes^L f^* F)\ .$$
 \end{lem}
\begin{proof}
We use the same argument as for Lemma \ref{system80} based on the adjoint pair $(f^*,Rf_*)$ and Lemma \ref{tens-pres1}. \end{proof}

\subsection{}

Let $f\colon X\to Y$ be a morphism of locally compact stacks.  
\begin{lem}\label{gr213}
Let $F\in \Sh_{R-\Mod}\bY$ be sheaf which is locally isomorphic to $\uRR_{\bY}$, i.e. there exist an atlas
$a\colon U\to Y$ such that $a^*F\cong  \uRR_{\bU}$.
In this case we have the projection formula:
For all $G\in \Sh_{R-\Mod}\bX$ or $H\in D^+(\Sh_\Ab\bX)$ the natural morphism
$$f_*G\otimes F\to f_*(G\otimes f^* F)\ , \quad Rf_* H\tensor^L F\to Rf_*(H\tensor^Lf^*F)$$ are isomorphisms.
 \end{lem}
\begin{proof}
This can be checked locally on the atlas $U\to Y$.
We consider the pull-back
$$\xymatrix{V\ar[r]^b\ar[d]^g&X\ar[d]^f\\
U\ar@{:>}[ur]\ar[r]^a&Y}\ .$$
 We must check that
$$a^*\circ (f_*G\otimes F)\to a^*\circ f_*(G\otimes f^* F)$$
is an isomorphism.
This map is equivalent to

\begin{eqnarray*}
             a^*(f_*G\tensor F) &\cong& a^*f_*G\tensor a^*F\\& \cong& a^*f_*G\tensor
        \underline{R}_U\\
        &\cong& a^*f_*G\\&\cong& g_*b^*G\\
        &\cong& g_*b^*(G\tensor\underline{R}_X)\\& \cong& g_*(b^*G\tensor
        b^*f^*\underline{R}_Y) \\
        &\cong& g_*(b^*G\tensor g^*a^*\underline{R}_Y)\\& \cong& g_*(b^*G\tensor
        g^*a^*F)\\
        &\cong& g_*b^*(G\tensor f^*F)\\& \cong &a^*f_*(G\tensor f^*F)\ .
     \end{eqnarray*}
The derived version is shown in  similar manner. 
\end{proof} 

\subsection{}

We will also need the projection formula with different assumptions.
Let $f\colon X\to Y$ be a map of locally compact stacks.
We consider $F\in \Sh_{R-\Mod}\bY$ and $G\in \Sh_{R-\Mod}\bX$.

\begin{lem}\label{system81}
Assume that $f$ is proper and representable, and that  $F$ is flat.
Then the natural transformation
$$f_* G\otimes F\to f_*(G\otimes f^*F)$$
of \ref{system80}
is an isomorphism.
\end{lem}
\begin{proof}
Using the observations \ref{obstr1a} we see that
it suffices to show that for all $(U\to Y)\in \bY  $
the induced morphism
\begin{equation}\label{kassp}
g_*G_V\otimes F_U\to g_*(G_V\otimes g^*F_U)
\end{equation}
is an isomorphism. Here $g\colon V\to U$ is the proper map of locally compact spaces defined by the Cartesian diagram
$$\xymatrix{V\ar[d]^g\ar[r]&X\ar[d]^f\\U\ar[r]&Y}\ .$$
The map (\ref{kassp}) is an isomorphism by
\cite[Prop. 2.5.13]{MR1299726}. \end{proof}

\subsection{}
 
We also have  a derived version of the projection formula in the case of sheaves of abelian groups.
The main point is that the ring $\Z$ has a  finite cohomological dimension (in fact equal to $1$). 
Let $f\colon X\to Y$ be a morphism of locally compact stacks.
\begin{lem}\label{projefoa}Assume that $f$ is proper and representable. 
If $G\in D^+(\Sh_\Ab\bY  )$ and $F\in D^+(\Sh_\Ab \bX  )$,
then we have
$$Rf_*G\otimes^LF\stackrel{\sim}{\to} Rf_*(G\otimes^L f^* F)$$
in $D^{+ }(\Sh_{\Ab} \bY  )$.
\end{lem}
\begin{proof}
As in the proof of Lemma \ref{system81} we can reduce to the small sites $(U)$ for all objects $(U\to Y)\in \bY$. After this reduction we apply
\cite[Prop. 2.6.6]{MR1299726} and the fact that the cohomological dimension of $\Z$ is $1$, hence finite. 
\end{proof}

\subsection{}

The following derived adjunction again uses the finiteness of the cohomological dimension of $\Z$.

\begin{lem}\label{khkqdqwwqc}
For $F,G,H\in D^+(\Sh_\Ab\bX)$ we have a natural isomorphism
$$R\Hom_{\Sh_\Ab\bX}(F\otimes^LG,H)\cong R\Hom_{\Sh_\Ab\bX}(F,R\uHom(G,H))\ .$$
\end{lem}
\begin{proof}
If $G\in \Sh_\Ab\bX$ is flat and $H\in \Sh_\Ab\bX$ is injective,  then the functor $$\Sh_\Ab\bX\ni F\mapsto \Hom_{\Sh_\Ab\bX}(F,\uHom(G,H))\cong \Hom_{\Sh_\Ab\bX}(F\otimes G,H)\in \Ab$$ is, as a composition of exact functors, exact. It follows that $\uHom(G,H)$ is again injective.
We now show the Lemma.
We can assume that $H$ is a complex of injectives. Furthermore, since the cohomological dimension of $\Z$ is one, hence in particular finite, we can assume  that $G$  is a complex of flat sheaves.
Then we have
\begin{eqnarray*}
R\Hom_{\Sh_\Ab\bX}(F\otimes^LG,H)&\cong&
\Hom_{\Sh_\Ab\bX}(F\otimes G,H)\\
&\cong&\Hom_{\Sh_\Ab\bX}(F,\uHom(G,H))\\
&\cong&R\Hom_{\Sh_\Ab\bX}(F,\uHom(G,H))\ .
 \end{eqnarray*}
\end{proof}

\section{Verdier duality for locally compact stacks in detail}

\subsection{}

Let $f\colon X\to Y$ be a map of locally compact stacks.
\begin{ddd}
We say that the cohomological dimension of $f_*$ is not greater  than $n\in \nat$ if the derived functor $R^if_*:\Sh_\Ab\bX\to \Sh_\Ab\bY$ vanishes for all $i>n$. 
\end{ddd}

The main theorem of the present subsection is
\begin{theorem}\label{main123}
Assume that $f:X\to Y$ is a representable and proper map between locally compact stacks such that
 $f_*$  has finite cohomological dimension. Then the functor $Rf_*\colon D^+(\Sh_\Ab\bX  )\to D^+(\Sh_\Ab\bY  )$ admits a right adjoint
$f^!\colon  D^+(\Sh_\Ab\bY  )\to  D^+(\Sh_\Ab\bX  )$.
\end{theorem}

The proof of Theorem \ref{main123}  will be finished in \ref{uidbhqwodqwdwqddqw}.
The main idea is to
transfer the construction of $f^!$ from \cite[Section 3.1]{MR1299726} to the
present situation.

\subsection{}\label{canokcioa}
 
  We consider the functorial flabby resolution (see \ref{system100}) of the
  constant sheaf $\uZ_{\bX}\to \Fl(\uZ_{\bX})$ and form the truncated complex
  $K:=\tau^{\le n} \Fl(\uZ_{\bX})$ so that in particular
  $K^n=\ker(\Fl^n(\uZ_{\bX})\to \Fl^{n+1}(\uZ_{\bX}))$.  
\begin{lemma}\label{shortrest12}
Assume that $f$ is representable and that $f_*$ has cohomological dimension  not greater than $n$.
 Then the complex
\begin{equation}\label{reconl}0\to \uZ_{\bX }\to K^0\to K^1\to\dots\to K^n \to 0\end{equation}
is a flat and $f_*$-acyclic resolution of $\uZ_{\bX}$.
\end{lemma}
\begin{proof}
The sheaf $\ker(K^n\to K^{n+1})$ is a torsion-free subsheaf of a torsion-free sheaf
and therefore flat (compare \cite[Lemma 3.1.4]{MR1299726}). 
By Lemma \ref{hdqoiwdwqw} the flabby sheaves $K^i$ for $i=0,\dots,n-1$ are $f_*$-acyclic.
In order to see that $K^n$ is $f_*$-acyclic, it
 suffices to show that
$R^if_* (\ker(K^n\to K^{n+1}))\cong 0$ for $i\ge 1$.
In fact, we have 
$R^if_* (\ker(K^n\to K^{n+1}))\cong R^{i+n}f_* \uZ_{\bX }\cong 0$.
\end{proof}

\subsection{}

 The fibers of a representable and proper morphism of topological stacks are compact. This is explicitly used in the proof of the following Lemma.

\begin{lem}\label{sumpre45}
If $f:X\to Y$ is  a representable and proper morphism of locally compact stacks, then
the functor $f_*\colon \Sh_\Ab\bX  \to \Sh_\Ab\bY  $ preserves direct sums.
\end{lem}
\begin{proof}
 Let $(G_i)_{i\in I}$ be a family of sheaves in 
$\Sh_\Ab\bX  $. Then we have a canonical map
$$\bigoplus_{i\in I} \circ f_*(G_i)\to f_*\circ \bigoplus_{i\in I}(G_i)\ .$$
In order to show that this map is an isomorphism we show that the induced map
$$(\bigoplus_{i\in I} \circ f_*(G_i))_U\to (f_*\circ \bigoplus_{i\in I}(G_i))_U$$ is an isomorphism for all
$(U\to Y)\in \bY  $. Choose such  $(U\to Y)$ and consider the Cartesian diagram
$$\xymatrix{V\ar[d]^g\ar[r]&X\ar[d]^f\\U\ar[r]&Y}\ .$$ It suffices to show that the induced map
$$\bigoplus_{i\in I} \circ g_*(G_i)_U\to g_*\circ \bigoplus_{i\in I}(G_i)_U$$
is an isomorphism. We consider the induced map on the stalk at $x\in U$. Since the restriction to $g^{-1}(x)$ commutes with the sum and $g^{-1}(x)$ is compact it is given by
$$\bigoplus_{i\in I} \circ \Gamma(g^{-1}(x),[(G_i)_U]_{|g^{-1}(x)})\to
\Gamma(g^{-1}(x),\bigoplus_{i\in I} [(G_i)_U]_{|g^{-1}(x)})$$
(see \cite[Proposition 2.5.2]{MR1299726}). But this last map is an isomorphism since the
global section functor on sheaves on a compact space commutes with sums.
\end{proof}

\subsection{}\label{uihiewfe}
Fix $j\in \{0,1,2\dots,n\}$ and set $K:=K^j$, see \ref{canokcioa}
\begin{lem}\label{ggtre1}
Let $f:X\to Y$ be a representable, proper morphism of locally compact stacks such that $f_*$ has cohomological dimension not greater than $n$.
Then the functor
$G\mapsto f_*(G\otimes K)$ is an exact functor $\Sh_\Ab \bX  \to \Sh_\Ab \bY  $.
Furthermore, $G\otimes K$ is $f_*$-acyclic.
\end{lem}
\begin{proof}
In the following proof we freely use the facts listed in \ref{desc_sheaves_on_U}.
Let $G^\cdot$ be an exact complex in $\Sh_\Ab \bX  $.
For $(U\to Y)\in \bY  $  consider the Cartesian  diagram
$$\xymatrix{V\ar[d]^g\ar[r]&X\ar[d]^f\\U\ar[r]&Y}\ .$$
Note that $(V\to X)\in \bX  $.
By construction (see \cite[Lemma 3.1.4]{MR1299726}) $K_V$ is flat and $g$-soft.
The complex $G_V^\cdot$ is exact.
By \cite[Lemma  3.1.2 (ii)]{MR1299726}
the complex $g_*(G_V^\cdot\otimes K_V)=(f_*(G^\cdot\otimes K))_U$ is exact.
Since this is true for all $(U\to Y)\in \bY  $ we conclude that
$f_*(G^\cdot\otimes K)$ is exact.

We now show that $G\otimes K$ is $f_*$-acyclic.
We must show that
$R^if_*(G\otimes K)\cong 0$ for all $i\ge 1$.
For $(U\to Y)\in \bY  $ as above we have
$(R^if_*(G\otimes K))_U\cong R^ig_*(G_U\otimes K_U)\cong 0$,
since $G_U\otimes K_U$ is $g$-soft by \cite[Lemma 3.1.2 (i)]{MR1299726} (note that $K_U$ is flasque and flat).
Since $(U\to Y)$ was arbitrary this implies that $R^if_*(G\otimes K)\cong0$
\end{proof}

\subsection{}

For $(V\to X)\in \bX$ let $\hat h_V^\Z$ denote the sheafification of the presheaf $h_V^\Z$, the presheaf of free abelian groups generated by the sheaf $h_V$ represented by $V$.
We consider the functor $f^!_K\colon \Sh_\Ab\bY  \to  \Pr_\Ab \bX  $ which associates to
a sheaf  $F\in \Sh_\Ab\bY  $ the presheaf
$f^!_{K}( F)\in \Pr_\Ab \bX $ given by
$$\bX  \ni (V\to X)\mapsto f^!_KF(V):=\Hom_{\Sh_\Ab\bY}(f_*(\hat h_V^\Z\otimes K),F)\in \Ab\ .$$
Note that $K\to f^!_K(F)$ is also a functor in $K$ (for fixed $F$).

\begin{lem}\label{injpw}
Let $K$ be as in \ref{uihiewfe} and $f:X\to Y$ be a representable, proper morphism of locally compact stacks such that $f_*$ has cohomological dimension not greater than $n$.
Assume that  $F\in \Sh_\Ab\bY  $ is an injective sheaf.   Then
$f^!_K(F)$ is an injective  sheaf.  Furthermore, for $G\in \Sh_\Ab\bX  $ there is a canonical isomorphism
\begin{equation}\label{ggtre}\Hom_{\Sh_\Ab\bY}(f_*(G\otimes K),F)\cong \Hom_{\Sh_\Ab\bX}(G,f^!_K(F))\ .\end{equation}
\end{lem}
\begin{proof}
We show that $f^!_KF$ is a sheaf by copying the corresponding argument in the proof of \cite[Lemma 3.1.3]{MR1299726}. The functor
$G\mapsto \Hom_{\Sh_\Ab\bY}(f_*(G\otimes K),F)$ is exact by Lemma \ref{ggtre1} and injectivity of $F$. If we establish the isomorphism (\ref{ggtre}), then we also have shown that $f^!_K(F)$ is injective.

For $(W\to X)\in \bX$ we have a canonical isomorphism
\begin{equation}\label{sus55}\Hom_{\Sh_\Ab\bY}(f_*(\hat h_W^\Z\otimes K),F)= f^!_K (F)(W)\cong\Hom_{\Sh_\Ab\bX}(\hat h_W^\Z,f^!_K(F))\ .\end{equation}
For a system $(G_i)_{i\in I}$ of sheaves we have a natural map $\colim_{i\in I}\circ f_*(G_i)\to f_*\circ \colim_{i\in I} (G_i)$.
For $G\in \Sh_\Ab\bX  $ we get
\begin{eqnarray*}
 \Hom_{\Sh_\Ab\bY}(f_*(G\otimes K),F)&\cong &\Hom_{\Sh_\Ab\bY}(f_*((\colim_{\hat h_W^\Z\to G} \hat h_W^\Z)\otimes K),F)\\
&\stackrel{!}{\cong}&\Hom_{\Sh_\Ab\bY}(f_*\circ \colim_{\hat h_W^\Z\to G} (\hat h_W^\Z\otimes K),F)\\
&\to&\Hom_{\Sh_\Ab\bY}(\colim_{\hat h_W^\Z\to G}\circ  f_*(\hat h_W^\Z\otimes K),F)\\
&\cong&\lim_{\hat h_W^\Z\to G}  \Hom_{\Sh_\Ab\bY}(f_*(\hat h_W^\Z\otimes K),F)\\
&\cong&\lim_{\hat h_W^\Z\to G}\Hom_{\Sh_\Ab\bX}(\hat h_W^\Z,f^!_K(F))\\
&\cong&\Hom_{\Sh_\Ab\bX}(\colim_{\hat h_W^\Z\to G}\hat h_W^\Z,f^!_K(F))\\
&\cong&\Hom_{\Sh_\Ab\bX}(G,f^!_K(F))\ .
\end{eqnarray*}
The marked isomorphism uses that the tensor product of sheaves commutes with colimits, a consequence of the fact  \ref{system101} that it is part of a closed monoidal structure.
It remains to show that this composition is an isomorphism.
If we write out the definition of the colimit in
$G\cong \colim_{\hat h_W^\Z\to G}\hat h_W^\Z$ we obtain an exact sequence of the form
\begin{equation}\label{sus54}\bigoplus_{j\in J}\hat h^\Z_{W_j}\to \bigoplus_{i\in I}\hat h^\Z_{V_i}\to G\to 0\ .\end{equation}
Now observe that for any collection $(G_i)_{i\in I}$ of sheaves in $\Sh_\Ab \bX  $ we have
$$\Hom_{\Sh_\Ab\bY}(f_*((\bigoplus_i G_i)\otimes K),F)\cong \prod_{i\in I} \Hom_{\Sh_\Ab\bY}(f_*(G_i\otimes K),F)$$
since $f_*$ (Lemma \ref{sumpre45})  and $\dots\otimes K$ commute with sums.
Clearly we also have
$$\Hom_{\Sh_\Ab\bX}(\bigoplus_i G_i,f^!_K(F))\cong \prod_{i\in I} 
\Hom_{\Sh_\Ab\bX}(G_i,f^!_K(F))\ .$$ 
{}From (\ref{sus54}) we thus get the diagram
\begin{equation*}
  \begin{CD}
    0 && 0\\
 @VVV @VVV \\
   \Hom_{\Sh_\Ab\bY}(f_*(G\otimes K),F) @>>> \Hom_{\Sh_\Ab\bX}(G,f^!_K(F)) \\
   @VVV @VVV\\
  \prod_{i\in I}\Hom_{\Sh_\Ab\bY}(f_*(\hat h_{V_i}^\Z\otimes K),F) @>\alpha>>\prod_{i\in I}\Hom_{\Sh_\Ab\bX}(\hat h_{V_i}^\Z,f^!_K(F)) \\
  @VVV @VVV\\
\prod_{j\in J}\Hom_{\Sh_\Ab\bY}(f_*(\hat h_{W_j}^\Z\otimes K),F)
@>\beta>>\prod_{j\in J}\Hom_{\Sh_\Ab\bX}(\hat h_{W_j}^\Z,f^!_K(F)) .\\
  \end{CD}
\end{equation*}
Because of the isomorphism (\ref{sus55}) the maps $\alpha$ and $\beta$ are isomorphisms.
The left vertical sequence is exact by Lemma \ref{ggtre1}.
The right vertical sequence is exact by the left-exactness of the $\Hom$-functor.
It follows from the five Lemma that (\ref{ggtre}) is an isomorphism.
\end{proof}

\subsection{}\label{uidbhqwodqwdwqddqw}

Let $I\Sh_{\Ab}\bX  \subset \Sh_{\Ab}\bX  $ denote the full subcategory of injective objects and $K^+(I\Sh_\Ab \bX  )$ be the category of complexes in $I\Sh_\Ab \bX $ which are bounded below, and whose morphisms are
homotopy classes of chain maps.
Then we have an equivalence of triangulated categories
$$K^+(I\Sh_\Ab \bX  )\cong D^+(\Sh_\Ab \bX  )\ .$$

Let $f:X\to Y$ be a representable, proper morphism of locally compact stacks such that $f_*$ has cohomological dimension not greater than $n$, and let $K^\cdot$ be as in \ref{canokcioa}.
We then define the functor
$f^!\colon K^+(I\Sh_\Ab \bY  )\to K^+(I\Sh_\Ab \bX  )$ by
$$f^!(F^\cdot)=(f^!_{K^\cdot}(F^\cdot))_{tot}\ ,$$
where $E^{\cdot,\cdot}_{tot}$ denotes the total complex of the double complex $E^{\cdot,\cdot}$.
Since $f^!_K$ preserves injective sheaves by Lemma \ref{injpw} this functor is well-defined.
Furthermore, for $F\in K^+(I\Sh_\Ab \bY  )$ and $G\in K^+(I\Sh_\Ab \bX  )$ we have by 
 Lemma \ref{injpw} a natural isomorphism between spaces of chain maps
$$\Hom_{C^+(\Sh_\Ab\bY  )}(f_*(G^\cdot\otimes K^\cdot)_{tot},F^\cdot)\cong  \Hom_{C^+(\Sh_\Ab\bX  )}(G^\cdot, f^!(F^\cdot))$$
which descends to an isomorphism on the level of homotopy classes
$$\Hom_{K^+(I\Sh_\Ab\bY  )}(f_*(G^\cdot\otimes K^\cdot)_{tot},F^\cdot)\cong  \Hom_{K^+(I\Sh_\Ab\bX  )}(G^\cdot, f^!(F^\cdot))\ .$$
Since $f^!(F^\cdot)$ is a complex of injective sheaves we have
$$\Hom_{K^+(I\Sh_\Ab\bX  )}(G^\cdot, f^!(F^\cdot))\cong \Hom_{D^+(\Sh_\Ab\bX  )}(G^\cdot, f^!(F^\cdot))\ .$$
Note that
$G^\cdot\cong G^\cdot\otimes \uZ_{\bX  }\to (G^\cdot\otimes K^\cdot)_{tot}$
is a quasi-isomorphism, and the complex $G^\cdot\otimes K^\cdot$ consists of $f_*$-acyclic sheaves by Lemma \ref{ggtre1}. Therefore $f_*(G^\cdot\otimes K^\cdot)_{tot}\cong Rf_*(G^\cdot)$.
Since $F^{\cdot}$ is injective we have
$$\Hom_{K^+(\Sh_\Ab\bY  )}(f_*(G^\cdot\otimes K^\cdot)_{tot},F^\cdot)\cong \Hom_{D^+(\Sh_\Ab\bY  )}(Rf_*(G^\cdot),F^\cdot)\ .$$
We conclude that
$$\Hom_{D^+(\Sh_\Ab\bY  )}(Rf_*(G^\cdot),F^\cdot)\cong  \Hom_{D^+(\Sh_\Ab\bX  )}(G^\cdot, f^!(F^\cdot))\ .$$
This finishes the proof of Theorem \ref{main123}. \hB

 \subsection{}

 We consider  morphisms $f\colon X\to Y$ and $u:U\to Y$ of locally compact stacks and form the Cartesian diagram $$\xymatrix{V\ar[r]^v\ar[d]^g&X\ar[d]^f\\U\ar[r]^u&Y}\ .$$

\begin{lem}\label{nattr121}
Assume the  $f$ is representable, proper and that $f_*$ has finite cohomological dimension.
Assume furthermore that $u$ has local sections. Then
we have a natural transformation
 $v^*\circ f^!\to g^!\circ u^*$.
\end{lem}
\begin{proof}
First note that $g$  is representable, proper and that $g_*$ has finite cohomological dimension.
Furthermore, $v$ has local sections.
 We apply $f^!$ to the unit
 $\id\to Ru_*\circ u^*$ and obtain a morphism
 \begin{equation}\label{plug1}f^!\to f^!\circ Ru_*\circ u^*\ .\end{equation}
 Since $f$ is representable and $u$ has local sections we have the isomorphism (see Lemma \ref{lem:pullpush} or \cite[Lemma 2.43]{bss}) 
 $$u^*\circ Rf_*\cong Rg_*\circ v^*\ .$$ Taking its right adjoint yields the isomorphism
 $$f^!\circ Ru_*\cong Rv_*\circ g^!\ .$$
 We plug this into (\ref{plug1}) and obtain a transformation
 $$f^!\to Rv_*\circ g^!\circ u^*\ .$$
Its adjoint is the desired transformation \end{proof}

\subsection{}

The following adjunction is a consequence of the derived projection formula Lemma \ref{projefoa}
and the derived adjunction Lemma \ref{khkqdqwwqc}
\begin{lem}
If  $f:X\to Y$ is a representable proper morphism of locally compact stacks which has local sections and is such  that $f_*$ has finite cohomological dimension, then for $G\in D^+(\Sh_\Ab\bX)$ and $F\in D^+(\Sh_\Ab\bX)$
we have a natural isomorphism
$$Rf_*R\uHom(G,f^!F)\cong R\uHom(Rf_*G,F)\ .$$
\end{lem}
\begin{proof}
Let $H\in D^+(\Sh_\Ab\bX)$ be arbitrary.
Then we calculate using Lemma \ref{derivedadj} and Lemma \ref{projefoa} that
\begin{eqnarray*}
R\Hom_{\Sh_\Ab\bY}(H,Rf_*R\uHom(G,f^!F))&\cong&R\Hom_{\Sh_\Ab\bX}(f^*H,R\uHom(G,f^!F))\\
&\cong&R\Hom_{\Sh_\Ab\bX}(f^*H\otimes^L G,f^!F)\\
&\cong&R\Hom_{\Sh_\Ab\bY}(Rf_*(f^*H\otimes^L G),F)\\
&\cong&R\Hom_{\Sh_\Ab\bY}(H\otimes^L Rf_*G,F)\\
&\cong&R\Hom_{\Sh_\Ab\bY}(H,R\uHom(Rf_*G,F))\ .
\end{eqnarray*}
\end{proof}

%
%
%
%
%
%
%
%
%
%
%
%
%
%

\subsection{}
\begin{ddd}
If  $f:X\to Y$ is a proper morphism of locally compact stacks such that $f_*$ has finite cohomological dimension, then we define the relative dualizing complex by $$\omega_{X/Y}:=f^!(\uZ_\bY)\ .$$
\end{ddd}
It would be interesting to know the structure of $\omega_{X/Y}$ for a topological submersion $f$ in the sense of \cite[Def. 3.3.1]{MR1299726}. 

\subsection{}

In a different setup  of Artin  stacks and  the lisse-\'etale site in \cite{laszlo-2005}
a six functor calculus was constructed. Starting with the observation that dualizing sheaves
on the small sites are sufficiently functorial the functors $Rf_!$ and $f^!$ are constructed on constructible sheaves by duality. In this approach one can relate the global $f^!$ 
with the local versions without any difficulty.  

A similar approach may work in the present topological context as well, but it is not clear how the
resulting $f^!$ will relate to the construction in the present paper.

\section{The integration map}\label{system103}

\subsection{}

Let $M$ be a closed connected orientable $n$-dimensional topological manifold.
\begin{ddd}
A map between locally compact stacks $f:X\to Y$ is a locally trivial fiber bundle with fiber $M$ if
for every space $U\to X$ the pull-back $U\times_YX\to U$ is a locally trivial
fiber bundle of spaces with fiber $M$.
\end{ddd}
Note that a locally trivial fiber bundle $f$ with fiber $M$ is representable, proper and has local sections, and $f_*$ has finite cohomological dimension.
 In order to see the last fact and to calculate $R^nf_*(\uZ_\bX)$ we consider $(U\to Y)\in \bY$ and let $V\to U$ be surjective and locally an open embedding such that we have a diagram with Cartesian squares
\begin{equation}\label{iudhqoiwdhqwidwqd}
\xymatrix{M\ar[d]^q&V\times_YX\ar[l]\ar[d]^h\ar[r]&U\times_YX\ar[d]^g\ar[r]&X\ar[d]^f\\{*}&V\ar[l]^p\ar[r]&U\ar[r]^u &Y}\ .
\end{equation}

The map $g$ is a topological submersion in the sense of  \cite[Def. 3.3.1]{MR1299726}.
As remarked in   \cite[Sec.~3.3]{MR1299726}   the cohomological dimension of $g_*$  is  not greater than $n$. This implies  $(R^if_*F)_U\cong R^ig_* (F_{U\times_YX})=0$ for all $i>n$.
Since this holds true for all $(U\to Y)\in \bY$ we conclude that $R^if_*F=0$ for all $i>n$.

We use the left part of the diagram (\ref{iudhqoiwdhqwidwqd})  in order to see that
$R^nf_*(\uZ_\bX)$ is locally isomorphic to $\uZ_\bY$.
In fact, we have
$$Rf_*(\uZ_\bX)_V\cong Rh_* \uZ_{(V\times_YX)}\cong p^* Rq_*\uZ_{(M)}\ .$$ 
A choice of an orientation of $M$ gives an isomorphism $R^nq_*\uZ_{(M)}\cong \uZ_{(*)}$ and therefore
$ R^nf_*(\uZ_\bX)_V\cong p^*\uZ_{(*)}\cong \uZ_{(V)}$.

\begin{ddd}\label{iquhuiqwdddqw}
A locally trivial fiber bundle $f\colon X\to Y$ with fiber $M$ is called orientable if
there exists an isomorphism $R^nf_*(\uZ_\bX)\cong \uZ_\bY$. An orientation of $f$ is a choice of such an isomorphism.
\end{ddd}

\subsection{}
Let $f:X\to Y$ be a locally trivial fiber bundle with fiber $M$, where $M$ is  a compact closed $n$-dimensional topological manifold.
We consider the $f_*$-acyclic and flat resolution $K$ defined in (\ref{reconl}).
The following was observed in \ref{uidbhqwodqwdwqddqw} 
\begin{corollary}\label{corol:calculate_Rf}
  The functor $Rf_*\colon D^+(\Sh_\Ab\bX)\to D^+(\Sh_\Ab\bY)$ is represented by
  $f_*\circ T_K$, where $T_K$ is tensor product with the complex $K$.
\end{corollary}

We now define a natural transformation
$$R\uHom(R^nf_*( \uZ_{\bX }),F)\to Rf_*\circ f^!(F)\ .$$
Let $F\in C^+(I\Sh_\Ab\bY )$ be a complex of injectives. 
We start from the observation that
$$R^nf_*( \uZ_{\bX })\cong f_*(K^n)/\im(f_*(K^{n-1})\to f_*(K^n))\ .$$
For $(U\to Y)\in \bY $ we thus obtain a chain of maps of complexes
\begin{eqnarray*}
\uHom(R^nf_* \uZ_{\bX },F)(U)
&\cong&\Hom_{\Sh_\Ab\bY}(\hat h_U^\Z,\uHom(R^nf_* \uZ_{\bX },F))\\&\cong&
\Hom_{\Sh_\Ab\bY}(\hat h^\Z_U\otimes R^nf_*( \uZ_{\bX }),F)\\&\cong&\Hom_{\Sh_\Ab\bY}(\hat  h_U^\Z\otimes f_*(K^n)/\im(f_*(K^{n-1})\to f_*(K^n)),F)\\
&\stackrel{!}{\to}&\Hom_{\Sh_\Ab\bY}(\hat h_U^\Z\otimes f_*(K),F)\\
&\stackrel{\ref{system81}}{\cong}&\Hom_{\Sh_\Ab\bY}(f_*(f^*\hat  h_U^\Z\otimes K),F)\\
&\stackrel{\ref{injpw}}{\cong}&\Hom_{\Sh_\Ab\bX}(f^*\hat   h_U^\Z ,f^!_K(F))\\
&\cong&\Hom_{\Sh_\Ab\bX}(\hat   h_U^\Z,f_*\circ f_K^!(F))\\
&\cong&f_*\circ f_K^!(F)(U)\ ,
\end{eqnarray*}
where the map marked by $!$ has degree $n$.
The projection formula Lemma \ref{system81} can be applied since $f^*\hat h^\Z_U$ is flat.
This transformation preserves homotopy classes of morphisms $F\to F^\prime$.
Since $F$ is injective we have
$$\uHom(R^nf_* \uZ_{\bX },F)\cong R\uHom(R^nf_* \uZ_{\bX },F)\ .$$
Further note that
$f^!_K(F)$ is still a complex of injectives by Lemma \ref{injpw}.
Therefore
$f_*\circ f_K^!(F)\cong Rf_*\circ f^!(F)$.
Hence this chain of maps of complexes induces a transformation
\begin{equation}\label{fetsch}R\uHom(R^nf_* \uZ_{\bX },F)\to Rf_*\circ f^!(F)\ .\end{equation}

\subsection{}

Its adjoint is a natural transformation
$$Rf_* f^* R\uHom(R^nf_* \uZ_{\bX },F)\to F\ .$$
Let us now assume that $f:X\to Y$ is in addition oriented by an isomorphism
$R^nf_*\uZ_\bX\cong \uZ_{\bY}$.
We precompose with this isomorphism  and get the integration map. 
\begin{ddd}\label{rdphi}
The integration map
$$\int_f\colon Rf_*\circ f^*\to \id$$ is the natural transformation of functors $D^+(\Sh_\Ab\bY)\to D^+(\Sh_\Ab\bY)$ of degree $-n$ defined as the composition
$$Rf_* f^*(F)\cong Rf_*f^*(\uHom( \uZ_{\bY },F))\cong 
Rf_* f^*(\uHom(R^nf_*( \uZ_{\bX }),F))\to F\ .$$ 
\end{ddd}

In Lemmas \ref{ttt12} and \ref{ttt13} we will verify in the more general case of unbounded derived categories that the integration map is functorial
with respect to compositions and compatible with pull-back diagrams.

\section{Operations with unbounded derived categories}\label{system4001}

%
%
%
%

\subsection{}

The category of sheaves $\Sh_\Ab\bX$ on a locally compact stack is a Grothendieck abelian category (see \ref{system104}). The category of complexes in a Grothendieck abelian category carries a model category structure (see \ref{ztoou6}). The unbounded derived category is the associated homotopy category.
The goal of the present subsection is to extend the sheaf theory operations $(f^*,f_*)$ and the integration map
$\int_f$ to the unbounded derived category.

Many results of the present subsection would continue to hold if one drops the assumption of local compactness in the definition of the site associated to stacks as well as for the stacks themselves.  
But the assumption of local compactness is important for the integration map
since it uses versions of the projection formula.

\subsection{}

Let $f\colon X\to Y$ be a morphism between locally compact  stacks.
Then we have an adjoint pair of functors
$$f^*\colon C(\Sh_\Ab\bY)\leftrightarrows C(\Sh_\Ab \bX):f_*\ .$$
In order to descend the functor $f_*$ to the bounded below derived category
it was sufficient to know that $f_*$ is left exact. In this case the idea is to apply 
$f_*$ to injective resolutions. The descent of the other functor $f^*$ is usually only considered if it exact, but see e.g. \cite{MR2312554} for more general constructions.  We know by  \ref{prexact} that  the functor $f^*$ is exact 
if $f$ has local sections.

It is not possible to show using the left exactness that $f_*$ preserves quasi-isomorphisms between unbounded complexes of injectives. Even worse, it is not clear how to resolve an unbounded complex by an injective complex. The method to descend $f_*$ to the derived category  uses abstract homotopy theory and works under the additional assumption that $f$ has local sections.

Recall that we use a model structure on the category $C(\Sh_\Ab\bX)$ of unbounded complexes of
sheaves for which the equivalences are the quasi-isomorphisms,  and the cofibrations are the
level-wise injections (see \ref{ztoou6}).  The inclusion
$C^+(\Sh_\Ab\bX)\hookrightarrow C(\Sh_\Ab\bX)$ of the full subcategory of
bounded below complexes induces an identification
$D^+(\Sh_\Ab\bX)\cong hC^+(\Sh_\Ab\bX)\hookrightarrow
hC(\Sh_\Ab\bX)=:D(\Sh_\Ab\bX)$ of the bounded below derived category as a full subcategory of the unbounded derived category.

The functor 
$Rf_*:D^+(\Sh_\Ab\bX)\to D^+(\Sh_\Ab\bY)$  is the adjoint of the restriction
of $f^*$ to the bounded below derived categories, and it is therefore the restriction of $Rf_*:D(\Sh_\Ab\bX)\to D(\Sh_\Ab\bY)$ to be defined below.

\begin{lem}
If the morphism $f\colon X\to Y$  of locally compact stacks  has local sections, then
$(f^*,f_*)$ is a Quillen adjoint pair.
\end{lem}
\begin{proof}
We use the  criterion \cite[Def. 1.3.1 (2)]{MR1650134} in order to show that $f^*$ is a left Quillen functor. We must show that it preserves cofibrations and trivial cofibrations. In other words, we must show that $f^*$ preserves injections and injections which induce isomorphisms on cohomology. Both properties follow from the exactness of $f^*\colon \Sh_\Ab\bY\to \Sh_\Ab \bX$. 
%
\end{proof}

\subsection{}

Let $f\colon X\to Y$ be a map between locally compact stacks which has local sections. Since $(f^*,f_*)$ is a Quillen adjoint pair
it induces a derived adjoint pair
$$Lf^*\colon hC(\Sh_\Ab\bY)\leftrightarrows hC(\Sh_\Ab \bX):Rf_*$$
(see Lemma \cite[Lemma 1.3.10]{MR1650134}).
Since $f^*$ is exact it directly  descends to the homotopy category.

%

\subsection{}

Let $g\colon Y\to Z$ be a second map of locally compact stacks which admits local sections.
Then we have adjoint  canonical isomorphisms
\begin{equation}\label{wjehfdajewd}
(g\circ f)^*\cong f^*\circ g^*\ ,\quad (g\circ f)_*\cong g_*\circ f_*\ .
\end{equation}
\begin{lem}\label{6.53}
We have a canonical isomorphism
$$R(g\circ f)_*\cong Rg_*\circ Rf_*\ .$$
\end{lem} 
\begin{proof}
Using \cite[Thm. 1.3.7]{MR1650134} we have a natural transformation
\begin{equation}\label{ajhwoirofg}
R(g\circ f)_*\cong R(g_*\circ f_*)\to Rg_*\circ Rf_*
\end{equation}
which is adjoint to
\begin{equation}\label{hjefdauzew}
Lf^*\circ  Lg^*\to L(f^*\circ g^*)\cong L(g\circ f)^*\ .
\end{equation}
Since $Lf^*$, $Lg^*$, and $L(g\circ f)^*$ are plain descents of
$f^*$, $g^*$, and $(g\circ f)^*$ to the homotopy category it follows from (\ref{wjehfdajewd}) that (\ref{hjefdauzew}) is an isomorphism. Therefore its adjoint
(\ref{ajhwoirofg}) is also an isomorphism. 
\end{proof} 

\subsection{}

Consider a Cartesian diagram of locally compact stacks 
$$\xymatrix{U\ar[d]^g\ar[r]^v&X\ar[d]^f\\V\ar[r]^u&Y}\ ,$$
where all maps have local sections.
Using the unit
$\id\to v_*\circ v^*$, the counit $u^*\circ u_*\to \id$,  and (\ref{wjehfdajewd}) we define (see Lemma \ref{lem:pullpush}) a transformation
$$u^*\circ f_*\to u^*\circ f_*\circ v_*\circ v^*\cong u^*\circ u_*\circ g_*\circ v^*\to g_*\circ v^*\ .$$ It is functorial with respect to compositions of such Cartesian diagrams.
By the same method we obtain a transformation
\begin{equation}\label{jfurfcscdde}
Lu^*\circ Rf_*\to Rg_*\circ Lv^*\ .
\end{equation}

\subsection{}
By Lemma \ref{lem:pullpush} we know that the transformation
$$u^*\circ f_*\to g_*\circ v^*$$ is in fact
an isomorphism. The derived version is more complicated
 and needs an additional assumption.
 
\begin{lem}\label{eiuwh}
Assume that $g$ is representable and  $g_*\colon \Sh_\Ab\bU\to \Sh_\Ab \bV$ has finite cohomological dimension. Then the transformation (\ref{jfurfcscdde}) is an isomorphism.
\end{lem}
\begin{proof}
We choose  fibrant replacement functors $$I_X\colon C(\Sh_\Ab\bX)\to
C(\Sh_\Ab\bX)\ ,\quad I_U: C(\Sh_\Ab\bU)\to
C(\Sh_\Ab\bU)\ .$$ 
In terms of these replacement functors we can write the compositions of derived functors as descents of quasi-isomorphism preserving functors on the level of chain complexes:
$$Lu^*\circ Rf_*\cong u^*\circ f_*\circ I_X\ ,\quad Rg_* \circ Lv^*\cong g_*\circ I_U\circ v^*\ .$$
Let $F\in  C(\Sh_\Ab\bX)$. We must show that the marked arrows  (induced by $\id\to I_U$ and $\id\to I_X$) in the following sequence are quasi-isomorphisms
$$u^*f_*I_X(F)\cong g_* v^*I_X(F)\stackrel{(*)}{\to} g_*I_Uv^*I_X(F)\stackrel{(**)}{\leftarrow} g_*I_Uv^*(F)\ .$$
The arrow marked by $(**)$ is a quasi-isomorphism since the functors
$g_*I_U$  and $v^*$ preserve quasi-isomorphisms, and  $F \to I_X(F)$ is a quasi-isomorphism.

The morphism $(*)$ is more complicated, and it is here where we need the assumption.
It is a property of the injective model structure on the chain complexes of a Grothendieck  abelian category that a fibrant complex consists of  injective objects.
An injective sheaf is in particular flabby. Since $v$ has local sections $v^*$  preserves flabby sheaves (Lemma \ref{flabbypres}). We conclude that $v^*I_X(F)$ is a complex of flabby sheaves.

Let $G\in C(\Sh_\Ab \bU)$ be a complex of flabby sheaves. We must show that
$g_*(G)\to g_*I_U(G)$ is a quasi-isomorphism. Since $g_*$ is an additive functor this assertion is equivalent to the assertion that $g_*(C)$ is exact, where $C$ is the mapping cone
of $G\to I_U(G)$. Note that $C$ is an exact complex of flabby sheaves. It decomposes into short exact sequences
$$0\to Z^n\to C^n\to Z^{n+1}\to 0\ ,$$
where $Z^n:=\ker(C^n\to C^{n+1})$.
Since $g$ is representable we know by Lemma \ref{hdqoiwdwqw}  that flabby sheaves are $g_*$-acyclic.
Therefore we obtain the exact sequence
$$0\to g_*(Z^n)\to g_*(C^n)\to g_*(Z^{n+1})\to R^1g_*(Z^{n})\to 0$$
and the isomorphisms
$$R^kg_*(Z^{n+1})\cong R^{k+1}g_*(Z^{n})$$
for all $k\ge 1$.
By induction we show that for $k\ge 1$ and all $l\in \nat$ we have
$$R^kg_*(Z^{n})\cong R^{k+l}g_*(Z^{n-l})\ .$$
Since we assume that $g_*$ has bounded cohomological dimension we conclude that
$R^k(Z^n)\cong 0$ for all $n\in \Z$ and $k\ge 1$.
In particular the sequences 
$$0\to g_*(Z^n)\to g_*(C^n)\to g_*(Z^{n+1})\to 0$$
are exact for all $n\in \Z$.  This shows the exactness of $g_*(C)$. 
\end{proof}
%

\subsection{}

Let now $f\colon X\to Y$ be a representable map between locally compact stacks which is an oriented locally trivial fiber bundle of closed  oriented manifolds of dimension $n$. In particular, $f$ has local sections and is proper, and $f_*$ has cohomological dimension $\le n$. We consider the canonical flabby resolution (see \ref{system100})
$$0\to \uZ_\bX \to \Fl^0(\uZ_\bX)\to \Fl^1( \uZ_\bX)\to \dots
\ .$$   
Then we know that $f_*\Fl( \uZ_\bX)$ is exact above degree $n$.
We let
$K $ denote the truncation (\ref{reconl}) of this resolution at $n$. 
Then the orientation of the bundle (see \ref{iquhuiqwdddqw}) gives the last isomorphism in the following composition
$$f_*K^n\to f_*K^n/\im(f_* K^{n-1}\to f_*K^n)\cong R^nf_*\uZ_\bX\cong  \uZ_\bY\ .$$ We let 
$T_{K}\colon C(\Sh_\Ab\bX)\to C(\Sh_\Ab\bX)$ denote the functor which associates to the complex
$F$ the total complex $T_{K}(F)$ of $F\otimes K$.
The  projection formula Lemma \ref{system81} for the proper representable map $f$ gives an isomorphism
$$f_*\circ T_{K}\circ f^*(F)\cong T_{f_*K}(F)$$
for complexes of flat sheaves $F\in C(\Sh_\Ab\bY)$.
The inclusion $ \uZ_\bX\to K$ and the projection $f_*K\to \uZ_\bY[-n]$ induces transformations
\begin{equation}\label{iofejoiwfwefwef}
\id\to T_{K}\ ,\quad T_{f_*K^\cdot}\to \id[-n]\ .
\end{equation}

\subsection{}

We know by Lemma \ref{ggtre1} that the functor $$f_*\circ T_K\colon \Sh_\Ab \bX\to \Sh_\Ab \bY$$
is exact. We choose a functorial fibrant replacement functor $\id\to I$ on $C(\Sh_\Ab\bX)$. Let $R:C(\Sh_\Ab\bY)\to C(\Sh_\Ab\bY)$ be the functorial flat
resolution functor of \ref{hdwidhwqdiwqdwdw}, extended to unbounded complexes.
Then we consider sequence
\begin{equation}\label{jhwjhdkqews}
f_*\circ I\circ f^*\to f_*\circ T_K\circ I\circ f^*\stackrel{!}{\leftarrow} f_*\circ T_K\circ f^*\stackrel{!}{\leftarrow} 
f_*\circ T_K\circ f^*\circ R \cong T_{f_* K}\circ R \to R[-n]\to  \id[-n]\ .
\end{equation} 
All functors in this sequence preserve quasi-isomorphisms and therefore descend plainly to the homotopy category $hC(\Sh_\Ab \bX)$. Since $f_*\circ T_K$ is exact the arrows marked by $!$  induce isomorphisms of functors on the homotopy category.
Now observe that  the descent of $f_*\circ I\circ f^*$ to the homotopy category is isomorphic to $Rf_*\circ Lf^*$.  Therefore (\ref{jhwjhdkqews}) induces a transformation
\begin{equation}\label{jkqwjhewuqieq3w}
\int_f\colon Rf_*\circ Lf^*\to \id[-n]\ .
\end{equation}

\begin{ddd}
The transformation (\ref{jkqwjhewuqieq3w}) is called the integration map.
\end{ddd}
It generalizes  Definition \ref{rdphi} from the bounded below to the unbounded derived category.

\subsection{}

In order to have a simple definition we have defined the integration map using a canonical resolution of $\uZ_\bX$ of length $n$. In fact, we can use more general resolutions. This will turn out to be useful for the verification of functorial properties of the integration map.

\subsection{}\label{system303}

Let us first recall some notation.
An object  $(U\to X)\in \bX$ represents the presheaf $h_U\in \Pr\bX$ (see also \ref{system202}).  We let $h_U^\Z\in \Pr_\Ab \bX$ be the free abelian presheaf generated by $h_U$ and form $\hat h_U^\Z := i^\sharp h_U^\Z \in \Sh_\Ab\bX$. 
\begin{ddd}
Let $f:X\to Y$ be a map of locally compact stacks.
A sheaf $F\in \Sh_\Ab\bX$ is called locally $f_*$-acyclic, if for every $(U\to X)\in \bX$ and $k\ge 1$ we have $R^k f_*(\hat h_U^\Z\otimes F)\cong 0$.
\end{ddd}

\subsection{}

 Let  $f:X\to Y$ be a map of locally compact stacks.
\begin{lem}\label{9odclkmfwef}
 Assume that the cohomological dimension of $f_*$ is bounded by $n$.
If $$L^0\to L^1\to \dots \to L^{n-1}\to L^n\to0$$ is an exact complex such that the
$L^i$ are $f_*$-acyclic (or locally $f_*$-acyclic) for $i=0,\dots,n-1$, then $L^n$ is $f_*$-acyclic (or locally $f_*$-acyclic, respectively). 
\end{lem}
This can be shown by a similar induction argument as used in the proof of Lemma \ref{eiuwh}. 
\hB

\subsection{}

Let $f:X\to Y$ be a map of locally compact stacks.
\begin{lem}\label{wgdqwededed}
Let $(V\to X)\in \bX$ and $F$ be locally $f_*$-acyclic. Then $\hat h_V^\Z\otimes F$ is locally $f_*$-acyclic.
\end{lem}
\begin{proof}
 Indeed, let $(U\to  X)\in \bX$.
Then we have 
$$\hat h_U^\Z\otimes(\hat h_V^\Z\otimes F)\cong (\hat h_U^\Z\otimes \hat h_V^\Z)\otimes F\ .$$
Furthermore we have
$$\hat h_U^\Z\otimes \hat h_V^\Z\stackrel{Lemma\: \ref{kion}}{\cong} i^\sharp (h^\Z_U\otimes^p h^\Z_V)\cong
i^\sharp (h_U\times h_V)^\Z\cong i^\sharp h_{U\times_\bX V}^\Z\cong 
  \hat h^\Z_{U\times_X V}\ ,$$
where we use the fact, that the absolute product in $\bX$ is given by the fiber product spaces over $X$ (\cite[Lemma 2.3.3]{bss}).
It follows that
$$R^kf_*(\hat h_U^\Z\otimes(\hat h_V^\Z\otimes F))\cong R^kf_*(\hat h_{U\times_X V}^\Z\otimes F)\cong 0$$
for all $k\ge 1$. 
\end{proof}

\subsection{}

Let $f:X\to Y$ be a map of locally compact stacks.
\begin{lem}\label{wejjqklxskj}
Assume that $f$ is proper, representable, and  that the cohomological dimension of $f_*$ is bounded.
If $F\in \Sh_\Ab\bX$ is flat and locally $f_*$-acyclic, then for any sheaf
$G\in \Sh_\Ab\bX$ the tensor product $G\otimes F$ is  $f_*$-acyclic and locally $f_*$-acyclic).
\end{lem}
\begin{proof}
We construct a resolution
$\dots\to G_j\to G_{j-1}\to \dots\to G_0\to G$, where all
$G_j$ are coproducts of sheaves of the form $\hat h^\Z_U$.
In fact, we have a surjection
$$\bigoplus_{\hat h_U^\Z\to G}\hat h_U^\Z \to G\ .$$
If we have already constructed $G_j\to\dots\to G_0\to G$, then we extend
this complex by
$$\bigoplus_{\hat h_U^\Z\to \ker(G_j\to G_{j-1})}\hat h^U_\Z \to G_j\ .$$
Since $F$ is flat, the complex
$$\dots \to G_j\otimes F\to \dots\to G_0\otimes F\to G\otimes F$$
is exact. The tensor product commutes with direct sums.
Therefore $G_j\otimes F$ is a sum of  $f_*$-acyclic sheaves, and by Lemma 
\ref{wgdqwededed} also of locally $f_*$-acyclic sheaves.  Since $f_*$ commutes with direct sums (Lemma \ref{sumpre45}) the sheaves $G_j\otimes F$ are themselves  $f_*$-acyclic and locally $f_*$-acyclic. With Lemma \ref{9odclkmfwef} we conclude that
$G\otimes F$ is  $f_*$-acyclic and locally $f_*$-acyclic. 
\end{proof}

\subsection{}

Let $f:X\to Y$ be a map of locally compact stack.
\begin{lem}
If $f$ is representable, then a flasque sheaf is locally $f_*$-acyclic.
\end{lem}
\begin{proof}
Let $F\in \Sh_\Ab\bX$ be flasque. We consider
 $(U\to Y)\in \bY$ and form the Cartesian diagram
$$\xymatrix{V\ar[r]\ar[d]^g&X\ar[d]^f\\U\ar[r]&Y}\ .$$
Then $(V\to X)\in \bX$ and we have $Rf_*(F)_U\cong Rg_*(F_V)$.
The restriction $F_V\in \Sh_\Ab(V)$ is still flasque.
A flasque sheaf on $(V)$ is $g$-soft (see  \cite[Definition 3.1.1]{MR1299726}).
But this implies that
$R^kg_*(F_V)=0$ for $k\ge 1$.
Since $U\to Y$ was arbitrary we see that
$R^kf_*(F)=0$ for $k\ge 1$.
\end{proof}

\subsection{}

Let us from now on until the end of this subsection assume that $f:X\to Y$ is a proper representable map of locally compact stacks which is an oriented locally trivial fiber bundle with fiber a closed connected topological manifold of dimension $n$.

Since a flat and flasque sheaf is locally $f_*$-acyclic and $K$ is a truncation of a flat and flasque resolution of $ \uZ_{\bX}$ we see by Lemma \ref{9odclkmfwef} that $K$ is a complex of flat and locally $f_*$-acyclic sheaves.  These are the two properties which make the theory work.

Let $L\to M$ be a quasi-isomorphism between upper bounded complexes of locally $f_*$-acyclic and flat sheaves. 

\begin{lem}\label{uuusaaassq}
For every  complex $F\in C(\Sh_\Ab\bX)$ the induced map
$$f_*(F\otimes L)\to f_*(F\otimes M)$$
is a quasi-isomorphism.
\end{lem}
\begin{proof}
We form the mapping cone $C$ of $L\to M$. It is an exact complex of locally $f_*$-acyclic and flat sheaves. Since the tensor product and $g_*$ commute with the formation of a mapping cone it suffices to show that $f_*(F\otimes C)$ is exact.
 
We know by Lemma \ref{wejjqklxskj} that $F\otimes C$ is a complex of $f_*$-acyclic sheaves. We claim  that $F\otimes C$ is exact. 

To this end we first show that $H\otimes C$ is exact for an arbitrary sheaf $H\in \Sh_\Ab\bX$.
We decompose the exact complex $C$ into short exact sequences   
$$S(k)\colon 0\to Z^k\to C^k\to Z^{k+1}\to 0$$ where
$Z^k:=\ker(C^k\to C^{k+1})$. 
Using the fact that the sheaves $C^k$ are flat we obtain
$$0\to \Tor_1(H,Z^{k+1})\to H\otimes Z^k\to H\otimes C^k \to H\otimes Z^{k+1}\to 0$$ and the isomorphisms $\Tor_{m+1}(H,Z^{k+1})\cong \Tor_{m}(H,Z^k)$ for all $m\ge 1$.
Since $\Z$ is one-dimensional we know that $\Tor_{m}\cong 0$ for $m\ge 2$. Inductively
we conclude that $\Tor_1(H,Z^k)\cong 0$ for all $k\in \Z$. 
It follows that $H\otimes S(k)$ is exact for all $k\in \Z$. 
This implies that $H\otimes C$ is exact.

Let now $F$ be a complex. Using the previous result and a spectral sequence argument we conclude that  the total complex associated to the double complex $F\otimes C$ is exact. 


Let now $C\in C(\Sh_\Ab\bX)$ be an exact  complex of $f_*$-acyclic sheaves.
We show that this implies that $f_*(C)$ is exact.
 The complex $C$  decomposes into short exact sequences
$$0\to Z^n\to C^n\to Z^{n+1}\to 0\ ,$$
where $Z^n:=\ker(C^n\to C^{n+1})$.
Using the fact that $C^n$ is $f_*$-acyclic we
obtain the exact sequence
$$0\to f_*(Z^n)\to f_*(C^n)\to f_*(Z^{n+1})\to R^1f_*(Z^{n})\to 0$$
and the isomorphisms
$$R^kf_*(Z^{n+1})\cong R^{k+1}f_*(Z^{n})$$
for all $k\ge 1$.
By induction we show that for $k\ge 1$ and all $l\in \nat$ we have
$$R^kf_*(Z^{n})\cong R^{k+l}f_*(Z^{n-l})\ .$$
Since  $f_*$ has bounded cohomological dimension we conclude that
$R^kf_*(Z^n)\cong 0$ for all $n\in \Z$ and $k\ge 1$.
In particular the sequences 
$$0\to f_*(Z^n)\to f_*(C^n)\to f_*(Z^{n+1})\to 0$$
are exact for all $n\in \Z$.  This shows the exactness of $f_*(C)$. 
\end{proof}

\subsection{}

\begin{lem}\label{edwiuhed}
The integration map is independent of the choice of a flat locally $f_*$-acyclic resolution
$K$ of $ \uZ_{\bX}$ of length $n$.
\end{lem}
\begin{proof}
Let $K,L$ are two such resolutions. Assume that there exists a quasi-isomorphism
$K\to L$.
The identification
$$\coker (f_*L^{n-1}\to f_*L^n)\cong \coker (f_*K^{n-1}\to f_*K^n)\cong R^nf_*(\uZ_\bX)\cong \uZ_\bY$$
gives a map $f_*L\to  \uZ_\bY[-n]$  which induces
the transformation
$T_{f_*L}\to \id$ of degree $-n$.

It induces a commutative diagram
$$\xymatrix{f_*If^*\ar[r]\ar@{=}[d]&f_*T_KIf^*\ar[d]&f_*T_Kf^*\ar[d]\ar[l]_\sim &f_*T_Kf^*R\ar[l]_\sim\ar[r]^\cong\ar[d]&T_{f_*K}R\ar[r]\ar[d] &R\ar[r]\ar@{=}[d]&\id\ar@{=}[d]\\
f_*If^*\ar[r]&f_*T_LIf^*&f_*T_Lf^*\ar[l]_\sim&f_*T_Lf^*R\ar[l]_\sim\ar[r]^\cong&T_{f_*L}R\ar[r]&R\ar[r] &\id}$$

The upper horizontal composition is the integration map defined using $K$ (see \ref{jhwjhdkqews}), and the lower horizontal composition is the integration map defined using $L$. We see that both maps are equal.

Let now $K,L$ again be flat and locally $f_*$-acyclic resolutions of $\uZ_\bX$ of length $n$.
We complete the proof of the Lemma by showing that there exists a third such resolution $M$
together with quasi-isomorphisms $K\stackrel{\sim}{\to}M\stackrel{\sim}{\leftarrow} L$.

The maps $\uZ_{\bX}\to K$ and $\uZ_{\bX}\to L$, respectively, induce maps
$K\to K\otimes L$ and $L\to K\otimes L$ which are quasi-isomorphisms.  We further get induced quasi-isomorphisms
\begin{equation}\label{sswuebh}
K\to \Fl(K\otimes L)\ ,\quad L\to \Fl(K\otimes L)\ .
\end{equation}
We let $M:=\tau^{\le n} \Fl(K\otimes L)$. 
The maps (\ref{sswuebh}) factorize over $M$.
Note that $K\otimes L$ is flat. Since $\Fl$ and truncation preserve flatness (see Lemma \ref{flat-preserv}), we see that $M$ is flat. Since $\Fl$ in fact produces flasque and hence locally $f_*$-acyclic resolutions, and the cohomological dimension of $f_*$ is bounded by $n$ we conclude by Lemma \ref{9odclkmfwef} that $M$ is locally $f_*$-acyclic.
\end{proof}

\subsection{}
In this paragraph we show that the integration map is functorial. 
Let $g\colon Y\to Z$ be a second proper and representable map of locally compact stacks which is
an oriented  locally trivial fiber bundle of closed   $m$-dimensional manifolds.

\begin{lem}\label{ttt12}
We have a commutative diagram
$$\xymatrix{Rg_*\circ Rf_*\circ Lf^*\circ Lg^*\ar[d]^{Rg_*(\int_f)}\ar[r]^{\cong}& R(g\circ f)_*\circ L(g\circ f)^*\ar[d]^{\int_{g\circ f}}\\
Rg_*\circ Lg^*[-n]\ar[r]^{\int_g}&\id[-n-m]}\ .$$
 \end{lem}
\begin{proof}
The following sequence of modifications transforms the down-right composition into the right-down composition. 
\begin{multline}
  g_*If_* If^*g^*\to g_*I f_*T_K If^*g^*\stackrel{\sim}{\leftarrow} g_*
  If_*T_K f^*g^* R
  \to g_*Ig^*R\\ 
  \to g_*T_L I g^*R\stackrel{\sim}{\leftarrow}
  g_*T_L g^*R \to\id \label{riocsn1}
\end{multline}
\begin{multline}
  g_*I f_* If^*g^*\to g_* T_L I f_* If^*g^*\to g_* T_L I f_*T_K
  If^*g^*
\stackrel{\sim}{\leftarrow} g_* T_L f_* T_K I f^*g^*\\
  \stackrel{\sim}{\leftarrow} g_* T_L f_*T_K f^*g^*R \to g_*T_L g^*R
  \to\id \label{riocsn2}
\end{multline}
\begin{multline}
  g_*I f_* If^*g^* \to g_* T_L I f_* If^*g^*\stackrel{\sim}{\leftarrow} g_*
  T_L f_* If^*g^*
 \to g_* T_L f_*T_K If^*g^*\\
\stackrel{\sim}{\leftarrow} g_*
  T_L f_* T_K f^*g^*R \to g_*T_L g^* R\to\id \label{riocsn3}
\end{multline}
\begin{multline}
  g_*I f_* If^*g^* \stackrel{\sim}{\leftarrow}g_* f_* If^*g^* \to g_* T_L f_*
  If^*g^*
 \to g_* T_L f_*T_K If^*g^*\\
\stackrel{\sim}{\leftarrow} g_* T_L f_*
  T_K f^*g^*R \to g_*T_L g^*R \to\id \label{riocsn4}
\end{multline}
\begin{multline}
  g_* f_* If^*g^* \to g_* T_L f_* If^*g^*\to g_* T_L f_*T_K If^*g^*
  \stackrel{\sim}{\leftarrow} g_* T_L f_*T_KR  If^*g^*\\
\stackrel{\sim}{\leftarrow} g_* T_L f_* T_KR f^*g^*R \to g_*T_L g^*R
  \to\id \label{riocsn44}
\end{multline}
\begin{multline}
  g_* f_* If^*g^* \to g_* T_L f_*T_K If^*g^*\stackrel{\sim}{\leftarrow}g_*f_*
  T_{f^*L\otimes K} R If^*g^*\\
\stackrel{\sim}{\leftarrow} g_*f_* T_{f^*L\otimes
    K}R f^*g^* R \to g_*T_L g^* R\to\id \label{riocsn5}
\end{multline}
\begin{multline}
  (g\circ f)_* I (g\circ f)^* \to (g\circ f)_* T_{M} I (g\circ
  f)^*
\stackrel{\sim}{\leftarrow} (g\circ f)_* T_{M} (g\circ f)^*R
  \to\id \label{riocsn6}
\end{multline}
The transition from (\ref{riocsn1}) to (\ref{riocsn2}) uses the fact  that tensoring  with $L$ and the map $\id\to T_L$ can be commuted with the intermediate operations.
In order to go from (\ref{riocsn2}) to (\ref{riocsn3}) we use the fact that $g_*T_L$ preserves quasi-isomorphisms. The same reason and the fact that $f_*$ preserves fibrant objects is behind the transition from (\ref{riocsn3}) to (\ref{riocsn4}). We use e.g. the isomorphism
$g_*f_*If^*g^*\stackrel{\sim}{\to} g_*I f_* I f^*g^*$. 
There is a vertical quasi-isomorphism from (\ref{riocsn44}) to (\ref{riocsn4}).
The step from (\ref{riocsn44}) to (\ref{riocsn5}) uses the isomorphism $T_Lf_*T_KR\stackrel{\sim}{\to}
 f_*T_{f^*L\otimes K}R$ given by the  projection formula.
 The weak equivalence in (\ref{riocsn5}) is not obvious (since $f^*L\otimes K$ is not obviously $g_*f_*$-acyclic), but follows from the fact, that this line is isomorphic to the previous
(\ref{riocsn44}).
In the last step from (\ref{riocsn5}) to (\ref{riocsn6}) we use the map $f^*L\otimes K\to M$ given by a  truncated flabby resolution of $f^*L\otimes K$ and the fact that the integration map is independent of the choice of the resolution.
\end{proof}

\subsection{}

Consider a cartesian diagram of locally compact stacks 
\begin{equation}\label{udgqwuidqwdqwdqwdw}\xymatrix{V\ar[d]^g\ar[r]^v&X\ar[d]^f\\U\ar[r]^u&Y}\ .\end{equation}
We assume that $f$ and $u$, and hence also $g$ and $v$ have local sections.
Furthermore we assume that $f$ is representable and a locally trivial oriented fiber bundle
with a closed manifold as fiber. Then $g$ has these properties, too.
The orientation of $g$ is induced by
$$R^ng_*\uZ_{\bV}\cong R^ng_*v^*\uZ_{\bX}\cong  u^*R^nf_*\uZ_{\bX}\cong u^*\uZ_{\bY}\cong \uZ_\bU$$
We get   diagrams
\begin{equation}\label{ieuwfhwefjkp}\xymatrix{u^*Rf_*f^*\ar[d]_{u^*\int_f}\ar[r]^{(\ref{jfurfcscdde})}&Rg_*v^*f^*\ar[d]^{(\ref{hjefdauzew})}\\u^*&Rg_*g^*u^*\ar[l]^{\int_g}} \end{equation}\begin{equation}\label{eq:comm_u_lower}
    \xymatrix{
      Ru_* Rg_* g^* \ar[d]^{Ru_*\int_g}\ar[r]&Rf_* Rv_* g^*\\
       Ru_* & Rf_*f^* Ru_*\ar[l]^{\int_f Ru_*}\ar[u]^\sim
    }
\end{equation}

 For the upper horizontal transformation in (\ref{ieuwfhwefjkp}) we use \ref{6.53}, and for the
  right vertical one (\ref{wqdwqdqw}) or \ref{eiuwh}. Note that
  only in the bounded below derived 
  category the right vertical morphism is an
  equivalence for general $u$ (which is anyway the situation in which we will
  apply the assertion). 

\begin{lem}\label{system300}\label{ttt13}\label{lem:funct_int} \label{lem:funct_int1}
The diagrams (\ref{ieuwfhwefjkp}) and  (\ref{eq:comm_u_lower})  commutes.
 \end{lem}


To prove Lemma \ref{system300}, we start with the following two technical
lemmas.
 \begin{lem}       \label{lem:two_projections}
    Given a Cartesian diagram (\ref{udgqwuidqwdqwdqwdw}) of locally compact stacks such that $f$ and $u$ have local sections, then  for 
     sheaves $K\in \Sh_\Ab\bX$  and $F\in \Sh_\Ab\bU$  the following diagram commutes:
    \begin{equation*}
      \begin{CD}
        f_*K\tensor u_* F @>=>> f_*K\tensor u_*F\\
        @VV{\ref{system80}}V @VV{\ref{system80}}V\\
          f_*(K\tensor f^*u_*F) &&  u_*(u^*f_*K\tensor F)\\
          @V{\sim}V{\ref{lem:pullpush}}V @V{\sim}V{\ref{lem:pullpush}}V\\
          f_*(K\tensor v_*g^*F) &&  u_*(g_*v^*K\tensor F)\\
          @VV{\ref{system80}}V @VV{\ref{system80}}V\\
          f_*v_*(v^*K\tensor g^*F)&&  u_*g_*(v^*K\tensor g^*F)\\
          @V{\sim}V{\ref{feriueewwwwzzz}}V  @V{\sim}V{\ref{feriueewwwwzzz}}V\\
          h_*(v^*K\tensor g^*F)@>=>> h_*(v^*K\tensor g^*F)
      \end{CD}\ ,
    \end{equation*}
    where $h:=f\circ v=u\circ g$.
  \end{lem}
  \begin{proof}
    By Definition \ref{system80}, the left vertical morphism is the
    image of the identity under the following sequence of maps
    \begin{multline*}
      \Hom(v^*K\tensor g^*K,v^*K\tensor g^*K) \to \Hom(v^*f^*f_*K\tensor
      v^*v_*g^*K,v^*K\tensor g^*K)\\
      \to \Hom(v^*(f^*f_*K\tensor
      f^*u_*K),v^*K\tensor g^*K)
      \to \Hom(f^*(f_*K\tensor u_*K)
      ,v_*(v^*K\tensor g^*K))\\
 \to \Hom(f_*K\tensor u_*K, f_*v_*(v^*K\tensor
      g^*K) ) \to \Hom(f_*K\tensor u_*K, h_*(v^*K\tensor g^*F)).
    \end{multline*}
    The right vertical morphism, on the other hand, is given by
    \begin{multline*}
       \Hom(v^*K\tensor g^*K,v^*K\tensor g^*K) \to \Hom(g^*g_*v^*K\tensor
      g^*u^*u_*K,v^*K\tensor g^*K)\\
      \to \Hom(g^*(u^*f_*K\tensor
      u^*u_*K),v^*K\tensor g^*K)
      \to \Hom(u^*(f_*K\tensor u_*K)
      ,g_*(v^*K\tensor g^*K))\\
 \to \Hom(f_*K\tensor u_*K, u_*g_*(v^*K\tensor
      g^*K) ) \to \Hom(f_*K\tensor u_*K, h_*(v^*K\tensor g^*F)).
    \end{multline*}
    In both cases, we first use the counit, then ``commute'' pushdown and
    pullback using Lemma \ref{lem:pullpush} and finally use adjunctions. By
    Lemma \ref{lem:pullpush}, the two ways to apply the counit and the push-pull isomorphism
    commute. This implies commutativity of the diagram of homomorphism sets,
    and therefore the commutativity of the original diagram.
  \end{proof}


  \begin{lemma}\label{lem:projection_and_pull}
In the situation of Lemma \ref{lem:two_projections}
for $K\in \Sh_\Ab\bX$ and $F\in \Sh_\Ab \bY$ the following diagram commutes:
    \begin{equation*}
      \begin{CD}
        u^*(f_*K\tensor F) @>{\ref{system80}}>> u^*f_*(K\tensor f^*F)\\
        @VV{\ref{tens-pres}}V @VV{\ref{lem:pullpush}}V\\
        u^*f_*K\tensor u^*F && g_*v^*(K\tensor f^*F)\\
        @VV{\ref{lem:pullpush}}V @VV{\ref{tens-pres}}V\\
        g_*v^*K\tensor u^*F && g_*(v^*K\tensor v^*f^*F)\\
        @VV=V @VV{\ref{uefhewiufuwefzzz}}V\\
         g_*v^*K\tensor u^*F @>>{\ref{system80}}> g_*(v^*K\tensor
        g^*u^*F) 
      \end{CD}\ .
    \end{equation*}
\end{lemma}
\begin{proof}
  The left vertical and lower composition is by definition the image of the
  identity under the sequence of maps
  \begin{multline*}
    \Hom(K\tensor f^*F,K\tensor f^*F) \xrightarrow{unit} \Hom(K\tensor
    f^*F,v_*v^* (K\tensor f^*F))\\
    \xrightarrow{adj} \Hom(v^*(K\tensor f^*F),v^*(K\tensor f^*F)) \to \Hom
    (v^*K\tensor g^*u^*F,v^*K\tensor g^*u^*F)\\
    \xrightarrow{counit} \Hom( g^*g_*v^*K \tensor g^*u^*F,v^*K\tensor g^*u^*F)
    \xrightarrow{adj} \Hom(g_*v^*K\tensor u^*F,g_*(v^*K\tensor g^*u^*F))\\
    \to \Hom(u^*(f_*K\tensor F),g_*(v^*K\tensor g^*u^*F)).
  \end{multline*}
  The upper and right vertical composition is the image of the identity under
  the sequence of maps
  \begin{multline*}
    \Hom(K\tensor f^*F,K\tensor f^*F) \xrightarrow{counit} \Hom(f^*f_*K\tensor
    f^*F,K\tensor f^*F)\\
    \xrightarrow{adj} \Hom(f_*K\tensor F,f_*(K\tensor f^*F))
    \xrightarrow{unit} \Hom(f_*K\tensor F,u_*u^*f_*(K\tensor f^*F))\\
    \xrightarrow{adj} \Hom(u^*(f_*K\tensor F),u^*f_*(K\tensor f^*F)) \to
    \Hom(u^*(f_*K\tensor F),g_*v^*(K\tensor f^*F))\\
    \to \Hom(u^*(f_*K\tensor F),g_*(v^*K\tensor v^*f^*F)) \to
    \Hom(u^*(f_*K\tensor F),g_*(v^*K\tensor g^*u^*F)).
  \end{multline*}

  These two maps coincide, as follows from the fact that units and counits
  commute (in the appropriate sense) with $\alpha_*$ and $\beta^*$.
\end{proof}
\subsection{}

We now show that (\ref{ieuwfhwefjkp}) commutes.
We  simplify the definition of the integration map which is represented by all horizontal compositions in the following diagram.
$$\xymatrix{f_*If^*\ar[r]\ar[d]^\sim&f_*T_KIf^*\ar[d]&f_*T_Kf^*R\ar[l]_{\sim}\ar[r]\ar[d]&\id\ar[d]^\sim\\
f_*If^*I\ar[r]^\sim&f_*T_KIf^*I&f_*T_Kf^*RI\ar[l]_{\sim}\ar[r]&I\ar@{=}[d]\\
f_*f^*I\ar[u]^\sim\ar[r]\ar[d]^\sim&f_*T_Kf^*I\ar[d]^\sim\ar[u]^\sim&f_*T_Kf^*RI\ar[l]^\sim\ar[d]^\sim\ar@{=}[u]\ar[r]&I\ar[d]^\sim\\
f_*f^*I\Fl\ar[r]&f_*T_Kf^*I\Fl&f_*T_Kf^*RI\Fl\ar[l]^\sim\ar[r]&I\Fl\\
f_*f^*\Fl\ar[u]^\sim\ar[r]&f_*T_Kf^*\Fl\ar[u]^\sim&f_*T_Kf^*R\ar[l]^\sim\Fl\ar[u]^\sim\ar[r]&\Fl\ar[u]^\sim
}$$
Let us comment about the isomorphisms in the first column. Let $F\in C(\Sh_\Ab\bX)$.
Then $f_*If^*(F)\to f_*If^*I(F)$ is a quasi-isomorphism since $f_*If^*$ preserves quasi-isomorphisms and $F\to I(F)$ is a quasi-isomorphism.
The map $f_*f^*I(F)\to f_*If^*I(F)$ is a quasi-isomorphism since $I(F)$ is a complex of injective, hence flabby sheaves, the functor $f^*$ preserves flabby sheaves, and therefore the
acyclic mapping cone of 
$C:=C(f^*I(F)\to If^*I(F))$ is an exact complex of flabby sheaves. In particular it is 
an exact complex of $f_*$-acyclic sheaves. Since $f_*$ has bounded cohomological dimension this implies that $f_*(C)$ is exact (see the argument in the proof of Lemma \ref{uuusaaassq}), and therefore $f_*f^*I(F)\to f_*If^*I(F)$ is a quasi-isomorphism.
The map $f_*f^*I(F)\to f_*f^*I\Fl(F)$ is a quasi-isomorphism by a similar argument.
In fact, $f^*\Fl(F)\to f^*I\Fl(F)$ is a quasi-isomorphism of $f_*$-acyclic sheaves.  This implies again by the  mapping cone argument, that
$f_*f^*\Fl(F)\to f_*f^*I\Fl(F)$ is a quasi-isomorphism.

The lower line of the diagram  (\ref{ieuwfhwefjkp}) expresses the integration map in terms of 
the flabby resolution functor $\Fl$. Since we know that $\Fl$ preserves flat sheaves (we do not know this for $I$) we can drop the flat resolution functor $R$ from the construction of the integration
by adopting the convention that  the functors are applied  to complexes of 
flat sheaves.

We get the following commutative diagram
\begin{equation}
  \label{eq:6.78}
  \begin{CD}
    u^* Rf_* f^* @>{\sim}>> u^* Rf_*f^* @>{u^*\int_f}>> u^*\\
    @VV{\sim}V @VV{\sim}V @VV{\sim}V\\
    u^*f_*T_K f^*\Fl @<{\sim}<< u^* T_{f_*K} \Fl @>>> u^*\Fl\\
    @VV{\sim}V @VV{\sim}V @VV{\sim}V\\
    g_*v^* T_K f^*\Fl &&  T_{u^*f_*K} u^*\Fl @>>> T_{u^*\Z}u^*\Fl\\
    @VV{\sim}V @VV{\sim}V @VV{\sim}V\\
    g_* T_{v^*K} v^*f^*\Fl && T_{g^*v_*K} u^*\Fl @>>> u^*\Fl\\
    @VV{\sim}V @VV=V @VV=V\\
    g_* T_{v^*K}g^* u^*\Fl @<{\sim}<< T_{g_*v^*K} u^*\Fl @>>> u^*\Fl\\
    @VV{\sim}V @VV{\sim}V @VV{=}V \\
    g_*T_{v^*K}g^*\Fl u^* @<{\sim}<< T_{g_*v^*K} \Fl u^* @>>> \Fl u^*\\ 
    @VV{\sim}V @VV{\sim}V @VV{\sim}V\\
    Rg_* g^* u^* @>=>> Rg_*g^*u^*  @>{\int_g u^*}>> u^*\\
  \end{CD}
\end{equation}
The commutativity of all the small squares is evident. The commutativity of
the large rectangle relies on the fact that the projection formula is
compatible with pullbacks, this is the statement of Lemma
\ref{lem:projection_and_pull}. The commutativity of the boundary of this diagram gives
 (\ref{ieuwfhwefjkp}).


\subsection{}

In order to show that (\ref{eq:comm_u_lower})  commutes we start with the following observation.

  \begin{lemma}\label{lem:two_integrations}
    Assume, in the situation of Lemma \ref{lem:two_projections}, that $K$ is a
    flat locally $f_*$-acyclic resolution of $\uZ_X$ of length $n$, and that $f$ is a
    projection of a locally trivial orientable fiber bundle of $n$-dimensional
    closed manifolds. Assume that $f_*K\to \uZ_Y$ is an orientation. Let
    $g_*v^*K\to \uZ_U$ be the induced orientation of the pullback bundle
    $g$. Then the following diagram commutes, where all the horizontal maps
    are given by the orientations.
    \begin{equation*}
      \begin{CD}
        f_*K\tensor u_*F @>>> \uZ_Y\tensor u_*F\\
        @VVV      @VVV\\
        u_*(u^*f_*K\tensor F) @>>> u_*(u^*\uZ_Y\tensor F)\\
        @VV{\sim}V   @VV{\sim}V\\
        u_*(g_*v^*K\tensor F)  @>>> u_*(\uZ_U \otimes F)
      \end{CD}
    \end{equation*}
  \end{lemma}
  \begin{proof}
    The upper diagram commutes because of the naturality of the homomorphism
    of the projection formula, the lower diagram commutes by the definition of
    the induced orientation of $g$.
  \end{proof}

To understand the relation between derived pushdown along a non-representable
map and 
integration we need to use an explicit model of the derived pushdown. If
$u\colon U\to Y$ is a morphism between locally compact stacks which has local
sections, then $Ru_*$ is given by $C_A\circ \Fl$, where $\Fl$ is the
functorial flabby resolution functor, and $C_A$ is defined in
Section \ref{iowefefwewqfqfefewf}, using an atlas $A\to U$. Note that $C_A$ indeed can be decomposed as the
composition of a functor $L_A$ on sheaves on $U$ and $u_*$. Here $L_A$ is the
sheafification of the functor on presheaves given by $${}^pL^k_AF(W\to 
U):= F(\underbrace{A\times_U\dots\times_U A}_{\text{$k+1$ factors}}\times_U
W\to U)\ .$$ i.e.~${}^pL^k_A={p_k}_*p_k^*$, with $p_k\colon
\underbrace{A\times_U\dots\times_U A}_{\text{$k+1$ factors}} \to U$.

  \begin{lemma}\label{lem:commute_pushdowns}
    In the situation of Lemma \ref{lem:two_integrations}, we obtain a
    commutative diagram
    \begin{equation*}
      \begin{CD}
        f_* T_K f^* u_* L_A\Fl @>=>>   f_* T_K f^* u_* L_A\Fl  @<{\sim}<<  T_{f_*K} u_* L_A \Fl @>>>
        u_*         L_A\Fl\\ 
     @VV{\sim}V   @VV{\sim}V   @VVV @VV=V\\
  f_*T_Kv_*L_{g^*A}g^*\Fl @>{\ref{pulcomghjdf}}>{\sim}>      f_*T_K v_*g^* L_A\Fl && u_*T_{u^*f_*K} L_A\Fl @>>> u_* L_A\Fl\\
&&        @VVV @VV{\sim}V  @VV=V\\
  &&      f_*v_*T_{v^*K} g^* L_A \Fl && u_*T_{g_*v^*K} L_A\Fl @>>> u_* L_A\Fl\\
    &&    @VV{\sim}V @VV=V @VV=V\\
    &&    u_*g_* T_{v^*K}g^* L_A\Fl @<<< u_*T_{g_*v^*K} L_A\Fl @>>> u_* L_A\Fl.
      \end{CD}
    \end{equation*}
    Here, the right horizontal maps are given by the orientations $f_*K\to\uZ_Y$
    and $g_*v^*K\to\uZ_U$.
  \end{lemma}
  \begin{proof}
    This is the direct translation of Lemma \ref{lem:two_projections} and
    Lemma \ref{lem:two_integrations}. 
  \end{proof}

Note that the upper composition is a representation (when applied to flat
sheaves) of
\begin{equation*}
  Rf_*f^*Ru_* \xrightarrow{\int_f} Ru_*.
\end{equation*}

The leftmost vertical arrow represents the morphism
\begin{equation}\label{eq:rightmap}
  Rf_*f^*Ru_* \to Rf_*Rv_*g^*,
\end{equation}
since $g^*$ preserves flabby sheaves, and $v_*L_{g^*A}$ indeed is a model for
$C_{g^*A}$, which can be used to calculate $Rv_*$.

Therefore the diagram in Lemma \ref{lem:commute_pushdowns}
contains one part (lower right-up) of the diagram (\ref{eq:comm_u_lower}).

\subsection{}

To represent the other composition of the diagram (\ref{eq:comm_u_lower}), we have to  commute not only $u_*$
but also $L_A$ with the other operations. Recall that $L_A$ provides some kind of a resolution, i.e.~we have a
  canonical map $\id\to L_A$, which is  used in the Lemma below.

\begin{lemma}\label{lem:LA_and_int}
  In the situation of Lemma \ref{lem:two_integrations}, the following diagram
  commutes, where the horizontal maps are induced by the orientation of
  $g$. 
  \begin{equation*}
    \begin{CD}
      u_* T_{g_*v^*K}L_A\Fl @>>>  u_* T_{\Z} L_A \Fl\\
        @VVV   @VVV\\
        u_* T_{L_A g_*v^*K} L_A \Fl @>>> u_* T_{L_A\Z} L_A\Fl\\
        @VVV @VVV\\
        u_*L_A T_{g_*v^*K} \Fl @>>> u_* L_A T_\Z \Fl
    \end{CD}
  \end{equation*}
  The second vertical map in each column follows from a variant of the
  projection 
  formula, using that $L_A$ is given by application of $(p_k)_*p_k^*$ (or by
  directly inspecting the definitions).
\end{lemma}
\begin{proof}
If $G\to H$ is a morphism of sheaves, then we get a natural transformation of functors $T_G\to T_H$. This naturality 
implies the commutativity of the first square.
  The second square
  is commutative by the naturality of the morphism in the projection formula.
\end{proof}

Observe that we have a natural isomorphism $g^*L_A\iso L_{g^*A}g^*$. 

\begin{lemma}\label{lem:move_LA}
  In the situation of Lemma \ref{lem:two_integrations}, we obtain the
  following commutative diagram
  \begin{equation*}
    \begin{CD}
      u_*g_* T_{v^*K}g^*L_A\Fl @<<< u_*T_{g_*v^*K} L_A\Fl\\
      @VVV @VVV\\
      u_*g_*T_{L_{g^*A}v^*K}g^*L_A\Fl @<<< u_* T_{g_*L_{g^*A}v^*K} L_A \Fl\\
      @V{\ref{pulcomghjdf}}V{\sim}V @VV{\sim}V\\
      u_*g_* T_{L_{g^*A}v^*K} L_{g^*A}g^*\Fl && u_*T_{L_Ag_*v^*K} L_A\Fl\\
      @VV{}V @VV{}V\\
      u_*g_* L_{g^*A} T_{v^*K} g^*\Fl && u_* L_A T_{g_*v^*K}\Fl\\
      @V{\ref{pulcomghjdf}}V{\sim}V   @VV=V\\
      u_* L_A g_*T_{g_*v^*K}g^*\Fl @<<< u_* L_A T_{g_*v^*K}\Fl
       \end{CD}
  \end{equation*}
\end{lemma}
\begin{proof}
  The upper square is commutative because of the naturality of the morphism in the projection
  formula. The commutativity of the lower rectangle follows from Lemma
  \ref{lem:two_projections}, as we basically have to commute two different
  applications of the projection formula.
\end{proof}

We now prove the commutativity of \eqref{eq:comm_u_lower}. Using explicit
representatives of the maps in question, we obtain (applied to flat sheaves)
\begin{equation*}
  \begin{CD}
    Rf_*f^* Ru_*   @>>=>  Rf_*f^*Ru_* @>{\int_fRu_*}>> Ru_*\\
    @VV{\sim}V  @VV{\sim}V @VV{\sim}V\\
        f_* T_K f^* u_* L_A\Fl  @<{\sim}<<  T_{f_*K} u_* L_A \Fl @>>>
        u_*         L_A\Fl\\
        @VVV @VVV @VV=V\\
        u_*g_* T_{v^*K}g^* L_A\Fl @<<< u_*T_{g_*v^*K} L_A\Fl @>>> u_* L_A\Fl\\
        @VVV @VVV @VV=V\\
      u_* L_A g_*T_{g_*v^*K}g^*\Fl @<<< u_* L_A T_{g_*v^*K}\Fl @>>> u_* L_A T_\Z \Fl\\
      @VVV @VVV @VV=V\\
      u_* L_A \Fl g_*T_{g_*v^*K}g^*\Fl @<<< u_* L_A \Fl T_{g_*v^*K}\Fl @>>>
    u_* L_A  \Fl\\
    @VV{\sim}V @VV{\sim}V @VV{\sim}V\\
    Ru_*Rg_*g^* @>>=> Ru_*Rg_*g^* @>{Ru_*\int_g}>> Ru_*
    \end{CD} 
\end{equation*}

Here, the first and the last rows are just added as illustration what the next
or preceding line, respectively, computes in the derived category. The map
from the third-last to the second-last row is induced by the inclusion into
the flabby resolution. This step is necessary because we don't know that the
functors in question are $u_*$-acyclic, and explains why one can directly
define only the map $f^*Ru_*\to Rv_*g^*$, and why it is hard to show that this
is an equivalence. The other vertical maps, and the commutativity of the
remaining four squares, is given by Lemmas \ref{lem:commute_pushdowns},
\ref{lem:LA_and_int}, \ref{lem:move_LA}.

Note that the left vertical composition is the composition 
\begin{equation*}
Rf_*f^* Ru_* \to Rf_* Rv_* g^* \to Ru_*Rg_*g^*,
\end{equation*}
as shown in the reasoning for \eqref{eq:rightmap}. The assertion follows.
\hB

\subsection{}

Compared with the simplicity of its statement the proof of Lemma \ref{system300} seems to be too long. But
let us mention that the proof of a similar result in the algebraic context is
quite involved, too. The book \cite{MR1804902} is devoted to this problem.

\section{Extended sites}\label{system3000}

\subsection{}

We consider the lower right Cartesian square of the  diagram
$$\xymatrix{&U\times_YB\ar@{.>}[r]\ar@{.>}[d]&B\ar@{.>}[d]\\A\times_YX\ar@{.>}[d]\ar@{.>}[r]&U\times_YX\ar[d]\ar[r]&X\ar[d]^f\\A\ar@{.>}[r]&U\ar[r]&Y}$$
in stacks where $U,X,Y$ are locally compact.
\begin{lem}\label{zqwduwdwdwdwqdwdzzz}
If $U$ is a space or $f$ is representable, then
$U\times_YX$ is a locally compact stack.
 \end{lem}
\begin{proof}
We first assume that $U$ is a locally compact space.
Let $B \to X$ be a locally compact  atlas. Then
$U\times_YB\to U\times_YX$ is an atlas. Indeed, 
surjectivity, representability, and local sections for this map are implied by the corresponding properties of the map $B\to X$.  The stack $U\times_YB$ is a space since $U\to Y$ is representable by 
Proposition \ref{lem:representability}. By Lemma \ref{qwuidiuwqdwqdwqd} the space $U\times_YB$
is locally compact.
Furthermore, again by Lemma \ref{qwuidiuwqdwqdwqd},
$$(U\times_YB)\times_{(U\times_YX)} (U\times_YB)\cong U\times_Y(B\times_XB)$$ is locally compact
since $B\times_XB$ is locally compact.
Hence the atlas $U\times_YB\to U\times_YX$ has the properties required in Definition \ref{qwuidiuwqdwqdwqd1fwefw} so that  $U\times_YX$ is a locally compact stack. 

We now assume that $f$ is representable. 
Let $A\to U$ be a locally compact atlas such that $A\times_UA$ is locally compact.
Then $A\times_YX\cong A\times_U(U\times_YX)\to U\times_YX$ is an atlas of $U\times_YX$.
We again verify  the properties required in Definition \ref{qwuidiuwqdwqdwqd1fwefw}.
By the special case of the Lemma already shown this atlas is locally compact. Moreover
$[A\times_U(U\times_YX)]\times_{U\times_YX}[A\times_U(U\times_YX)]\cong
(A\times_UA)\times_YX$ is locally compact.
\end{proof}

\subsection{}

If $f\colon X\to Y$ is a representable map with local sections between   locally compact stacks, then
for $(U\to Y)\in \bY$ we have ${}^pf^* h_U\cong h_{U\times_XY}$ (see the proof of Lemma \ref{hhdf36434746zzz} below). If we drop the assumption that $f$ is representable, then in general ${}^pf^*h_U$ is not representable. In order to overcome this defect we enlarge the site $\bX$ to $\tilde \bX$ so that it contains the stacks $U\times_XY\to X$ over $X$.

We consider the $2$-category $\Stacks^{top, lc}/_{ls,rep}X$ of locally compact stacks  $U\to X$ over $X$, where the structure map is representable and has local sections. A morphism in this category is a diagram
$$\xymatrix{U\ar[dr]\ar[rr]&\ar@{:>}[d]&V\ar[dl]\\&X&}$$
consisting of a one-morphism  and a two-morphism. The composition is defined in the obvious way.
If there is a two-morphism  between two such one-morphisms, then it is unique by the representability of the structure maps. Therefore $\Stacks^{top, lc}/_{ls,rep}X$ is equivalent in two-categories to the one-category obtained by identifying all isomorphic one-morphisms.

\subsection{}

Let $f:X\to Y$ be a map between locally compact stacks.
\begin{ddd}
We let $\tilde \bX$  be the category obtained from   $\Stacks^{top, lc}/_{ls,rep}X$   by identifying all isomorphic one-morphisms. 
\end{ddd}

%
%
%
%
%

 We now define the topology on $\tilde \bX$.
 A covering family $(U_i\to U)$ of $(U\to X)\in \tilde \bX$  is a family of  locally compact  stacks over $U$ such that $U_i\to U$ is representable, has local sections and  $\sqcup_{i\in  I} U_i\to U$ is surjective\footnote{These maps are actually equivalence classes, but in order to simplify the language we will not mention this explicitly in the following}. Using Lemma \ref{zqwduwdwdwdwqdwdzzz}
one easily checks the axioms listed in \cite[1.2.1]{MR1317816}.

Let $\hat \bX$ be  the site with the same underlying category as $\tilde \bX$, but with the topology generated by the   covering families of $(U\to X)$ given by families $(U_i\to U)\in \Stacks^{top, lc}/X$ such that $U_i\to U$ is a map from a locally compact space with local sections and $\sqcup_i U_i\to U$ is surjective.
 
\begin{lem}\label{wuiefhewffewfewfzzz}
We have a canonical isomorphism
$$\Sh \tilde \bX\cong \Sh\hat \bX \ .$$
\end{lem}
\begin{proof}
The covering families of $\hat \bX$ are covering families in
$\tilde \bX$. Here we use Proposition \ref{lem:representability} in order to see that the maps $U_i\to U$ from spaces $U_i$ are representable.
On the other hand, every covering family  $(U_i\to U)$ of $(U\to X)$ in $\tilde \bX$  can be refined
to a covering family in $\hat \bX$  by  choosing a locally compact   atlas $A_i\to U_i$ for each $U_i$.  
This implies the lemma. \end{proof}

\subsection{}

The natural functor $\Top^{lc}/X\to \Stacks^{top,lc}/X$  from locally compact spaces over $X$ to locally compact stacks over $X$ induces a map of sites $j\colon \bX\to \tilde \bX$.
\begin{lem}\label{k7783nndnbdzzz}
The restriction functor
$$j^*\colon \Sh\tilde \bX\to \Sh \bX $$ is  an equivalence of
categories. \end{lem}
\begin{proof}
The inverse of $j^*$ is the functor $j_*$ given by
$$j_* F(U):=\lim_{(V\to U)\in \bX//U} F(V)$$
for all $(U\to X)\in \tilde \bX$, where
$\bX//U$ is the category of all pairs $(V\in \bX, j(V)\to U\in  \Mor(\tilde \bX))$ such that the map
$j(V)\to U$ has local sections.

If $U\in j(\bX)$, then $(U,\id_{j(U)}\colon j(U)\to j(U))$ it is the final object of $\bX//U$.
This gives a natural isomorphism $j^* j_*(F)(U)\cong F(U)$.

We now define a natural isomorphism
$j_*j^*(F)\to F$ for all $F\in \Sh\tilde \bX$.
Let $(U\to X)\in \tilde \bX$.
The family $(V\to U)_{\bX//U}$ is a covering family
of $U\to X$ in $\hat \bX$. 
Since $F$ is also a sheaf on $\hat \bX$  by Lemma \ref{wuiefhewffewfewfzzz} we get an isomorphism
$$j_*j^*(F)(U)\cong \lim_{(V\to U)\in \bX//U} j^*(F)(V)\cong F(U)\ .$$
\end{proof}

\subsection{}

\begin{lem}\label{wkuhwqeddqwzzz}
A map $f:X\to Y$  between locally compact  
stacks induces a map of sites
$$\tilde f^\sharp\colon \bY\to \tilde \bX$$ by
$$\tilde f^\sharp(U\to Y):=U\times_YX\to X\ .$$
\end{lem}
\begin{proof}
Indeed, if $U\to Y$ is a map from a locally compact  space, then the stack
$U\times_YX$ is locally compact by Lemma \ref{zqwduwdwdwdwqdwdzzz}.
If $(U_i\to U)$ is a covering family of $(U\to Y)\in \bY$ by open subspaces, then
$(U_i\times_YX\to U\times_YX)$ is a covering family in $\tilde \bX$
by open substacks.

Furthermore it is easy to see that $\tilde f^\sharp$  preserves fiber products, i.e.
if $(U_i\to U)$ is a covering family and $V\to U$ is a morphism in $\bY$, then
$\tilde f^\sharp(U_i\times_UV)\cong \tilde f^\sharp(U_i)\times_{\tilde f^\sharp(U)}\tilde f^\sharp(V)$. 
\end{proof}

\subsection{}
We consider a map $f:X\to Y$  between locally compact  
stacks. Then we have  an adjoint pair of functors
$$\tilde f^\sharp_* \colon  \Sh \bY\leftrightarrows \Sh\tilde \bX : (\tilde f^\sharp)^* \ .$$
\begin{lem}\label{hhdf36434746zzz}
We have an isomorphism of functors
$j^*\circ \tilde f^\sharp_* \cong f^*\colon \Sh\bY\to \Sh\bX$ 
\end{lem}
\begin{proof}
 The map $j\colon \bX \to \tilde \bX $ induces a map ${}^pj^*\colon \Pr\tilde \bX \to \Pr\bX $.
We show the relation first on representable presheaves.
Let $(U\to Y)\in \bY $ and observe that $(U\times_YX\to X)\in \tilde \bX $
by Lemma \ref{zqwduwdwdwdwqdwdzzz}. 
The following chain of natural isomorphisms (for arbitrary $F\in \Pr\tilde \bX$) shows that
$\tilde f^\sharp_* h_U\cong h_{U\times_YX}$:
\begin{eqnarray*}
\Hom_{\Pr\tilde \bX}(\tilde f^\sharp_* h_U,F)&\cong&\Hom_{\Pr\bY}(h_U,(\tilde f^\sharp)^*F)\\&\cong&
(\tilde f^\sharp)^*F(U)\\&\cong&F(\tilde f^\sharp(U))\\&\cong&F(U\times_YX)\\&\cong&
\Hom_{\Pr\tilde \bX}(h_{U\times_YX},F)\ .
\end{eqnarray*}
For $(U\to Y)\in \bY $
we have $ {}^pf^* h_U\cong {}^pj^* h_{U\times_YX}$. Indeed,
for $(V\to X)\in \bX $ we have
$${}^pj^* h_{U\times_YX}(V)\cong\Hom_{\tilde \bX}(j(V), U\times_YX)
\stackrel{!}{\cong} {}^pf^* h_U(V)\ ,$$
where the marked isomorphism can be seen by making the definition of ${}^pf^*$ explicit.
Since ${}^pj^*\circ {}^p\tilde  f^\sharp_*$ and  ${}^p f^*$ commute with colimits
the equation ${}^pj^*\circ {}^p\tilde f^\sharp_*\cong {}^pf^*$ holds on all presheaves.
The restriction to sheaves (note that all functors preserve sheaves) gives
$j^*\circ \tilde f^\sharp_*\cong f^*$.
\end{proof}
By adjointness we get
\begin{equation}\label{uqiwdiqwdqwdwqdwqdwdwqd}
(\tilde f^\sharp)^*\circ j_*\cong f_*\ .
\end{equation}

\subsection{}

Consider two composeable  maps between locally compact stacks.
$$X\stackrel{f}{\to} Y\stackrel{g}{\to} Z\ . $$
The following lemma generalizes \cite[Lemma 2.23]{bss}
by dropping the unnecessary additional assumptions that $f$ has local sections or $g$ is representable.
\begin{lem}\label{feriueewwwwzzz}
We have an isomorphism of functors
$g_*\circ f_*\cong (g\circ f)_*\colon \Sh\bX\to \Sh \bZ$.
\end{lem}
\begin{proof}

We consider the following diagram:
$$\xymatrix{\Sh\bX\ar@/_2cm/[dd]^{(g\circ f)_*}\ar[r]^{j^X_*}\ar[d]^{f_*}&\Sh\tilde \bX\ar@/_-2cm/[dd]^{(\widetilde{(g\circ f)}^\sharp)^*}\ar[d]^{(\tilde f^\sharp)^*}\\
\Sh\bY\ar[r]^{j^Y_*}\ar[d]^{g_*}&\Sh\tilde \bY\ar[d]^{(\tilde g^\sharp)^*}\\
\Sh\bZ\ar[r]^{j^Z_*}&\Sh\tilde \bZ}\ .$$
We know that the squares commute (Equation (\ref{uqiwdiqwdqwdwqdwqdwdwqd})), and that the horizontal arrows are isomorphisms
(Lemma \ref{k7783nndnbdzzz}).
It follows from the constructions that
$$\tilde f^\sharp\circ \tilde g^\sharp=\widetilde{(g\circ f)}^\sharp$$
on the level of sites. Hence the right triangle commutes, too. This implies
commutativity of the left triangle. 
\end{proof}

Taking adjoints we get:
\begin{kor}\label{uefhewiufuwefzzz}
We have an isomorphism $f^*\circ g^*\cong (g\circ f)^*\colon \Sh\bZ\to \Sh\bX$.
\end{kor}

%

\subsection{}

We consider a topological stack $X$ and the inclusion $j\colon \bX\to \tilde \bX$ which induces by Lemma \ref{k7783nndnbdzzz}  an equivalence of categories of sheaves
$$j^*\colon \Sh\tilde\bX\leftrightarrows \Sh \bX\colon j_*\ .$$
Note that the notion of flabbiness depends on the site. 
\begin{ddd}
We call a sheaf $F\in \Sh_\Ab\bX$ strongly flabby if $j_*(F)$ is flabby.
\end{ddd}
Since flabbiness is a condition to be checked for all covering families
and since all covering families in $\bX$ induce covering families in $\tilde \bX$ it follows that a strongly flabby sheaf is flabby. Since injective sheaves are strongly flabby
each sheaf admits a strongly flabby resolution.

\subsection{}

Let $f\colon X\to Y$ be a morphism of locally compact stacks. 
\begin{lem}
Strongly flabby sheaves are $f_*$-acyclic.
\end{lem}
\begin{proof}
 In view of Lemma \ref{hhdf36434746zzz}
it suffices to show that flabby sheaves in $\Sh_\Ab\tilde\bX$
are $\tilde f_*$-acyclic. We now can write
$\tilde f_*=\tilde i^\sharp\circ {}^p\tilde f_*\circ \tilde i$, where
$\tilde i^\sharp$ and $\tilde i$ are the sheafification
functor and the inclusion of sheaves into presheaves for the tilded sites,
and ${}^p\tilde f_*={}^p(\tilde f^\sharp)^*\colon \Pr\tilde \bX\to \Pr\tilde \bY$.
Since
${}^p\tilde f_*(F)(V\to Y)=F(V\times_YX\to X)$ we see that
${}^p \tilde f_*$ is exact. Since strongly flabby sheaves are $\tilde i$-acyclic, and $\tilde i^\sharp$ is exact, it follows that strongly flabby sheaves are $\tilde f_*$-acyclic.
\end{proof}

\begin{lem}
The functor $$f_*\colon \Sh_\Ab\bX\to \Sh_\Ab\bY$$
 preserves strongly flabby sheaves. 
\end{lem}
\begin{proof}
We must show that $\tilde f_*$ preserves flabby sheaves.
Let $F\in \Sh_\Ab \tilde \bX$ and $\tau=(U_i\to U)$ be a covering family of $(U\to Y)$ in $\bY$. We must show that the \v{C}ech complex $C(\tau,\tilde f_*F)$ is acyclic. Note that
$\tilde f_*F(V)=F(V\times_YX)$. The family $f^\sharp(\tau):=(U_i\times_YX\to U\times_YX)$ is a covering family of $U\times_YX$ in $\tilde \bX$. We see that
$C(\tau,\tilde  f_*F)\cong C(f^\sharp(\tau),F)$. Since $F$ is strongly flabby, the complex $
C(f^\sharp\tau,F)$ is acyclic. 
\end{proof}

\subsection{}

Consider again a sequence of composeable  maps between locally compact stacks.
$$X\stackrel{f}{\to} Y\stackrel{g}{\to} Z\ . $$
The following Lemma generalizes \cite[Lemma 2.26]{bss}, again by dropping the unnecessary assumptions that $f$ has local sections or $g$ is representable.
\begin{lem}\label{keykey}
We have an isomorphism of functors
$$Rg_*\circ Rf_*\cong R(g\circ f)_*\colon D^+(\Sh_\Ab\bX)\to D^+(\Sh_\Ab \bZ).$$
\end{lem}
\begin{proof}
The isomorphism $(g\circ f)_*\to g_*\circ f_*$ induces a transformation
$R(g\circ f)_*\to Rg_*\circ Rf_*$.
Since injective sheaves are strongly flabby, $f_*$ preserves strongly flabby sheaves, and strongly flabby sheaves are $g_*$-acyclic, this transformation is indeed an isomorphism.
\end{proof}

\backmatter


\end{document}